\documentclass[a4paper,11pt]{article}
\usepackage[top=21mm,left=22mm]{geometry}\textwidth16.5cm\textheight24.3cm
\usepackage{amsmath,amssymb}
\usepackage{amsthm}
\usepackage[trick,hustyle]{optional}
\usepackage{url} 
\usepackage{graphicx}
\usepackage{pstricks}
\usepackage{amsfonts}
\usepackage{footmisc}
\usepackage{multirow}
\usepackage{array}
\usepackage{paralist}
\newcommand{\bpm}{\begin{pmatrix}}\newcommand{\epm}{\end{pmatrix}}
\def\pa{{\partial}}\def\lam{\lambda}
\def\huga#1{\begin{gather} #1 \end{gather}}

\newcommand{\reff}[1]{(\ref{#1})}
\def\Lh{\hat{L}}\def\CO{{\cal O}}
\def\ov{\overline}\def\er{{\rm e}}\def\ri{{\rm i}}
\newcommand{\hs}[1]{{\hspace{#1}}}\newcommand{\vs}[1]{{\vspace{#1}}}
\newcommand{\bce}{\begin{center}}\newcommand{\ece}{\end{center}}
\def\Uhet{U_{\text{het}}}\def\Uhom{U_{\text{hom}}}
\def\Ahom{A_{\text{hom}}}
\def\Afront{A_{\text{front}}}\def\Aback{A_{\text{back}}}
\def\ra{\to}\def\CS{{\cal S}}
\def\hot{\text{h.o.t.}}\def\eps{\varepsilon}\def\cc{{\rm c.c.}}
\def\ccf{{\rm c.c.f.}}
\newtheorem{theorem}{Theorem}[section]
\theoremstyle{definition}
\newtheorem{remark}[theorem]{Remark}

\def\brem{\begin{remark}}\def\erem{\end{remark}}
\def\eex{\hfill\mbox{$\rfloor$}}
\def\Om{\Omega}\def\al{\alpha}
\newcommand{\bi}{\begin{itemize}}\newcommand{\ei}{\end{itemize}}
\newcommand{\bci}{\begin{compactitem}}\newcommand{\eci}{\end{compactitem}}
\newcommand{\bcen}{\begin{compactenum}}\newcommand{\ecen}{\end{compactenum}}

\newcommand{\N}{\mathbb{N}}
\newcommand{\R}{\mathbb{R}}

\newcommand{\fettk}{\textbf{k}}
\newcommand{\kx}{m}\newcommand{\ky}{n}

\newcommand{\trh}{\text{tr} \hat{L}(k,\ld) }

\newcommand{\detlh}{\text{det}   \hat{L}(k,\ld)}

\newcommand{\etotal}{E_{\text{total}}}
\newcommand{\epot}{E_{\text{pot}}}\newcommand{\ekin}{E_{\text{kin}}}

\newcommand{\ld}{\lambda}
\newcommand{\sss}{\sum_{i=1}^3}
\newcommand{\rkl}{\right)}\newcommand{\lkl}{\left(}

\providecommand{\norm}[1]{\left\Vert#1\right\Vert}

\providecommand{\ssmatrix}[1]{\bigl(\begin{smallmatrix}#1 \end{smallmatrix} \bigr)}

\providecommand{\ali}[1]{\begin{align}#1\end{align}}
\providecommand{\alinon}[1]{\begin{align*}#1\end{align*}}

\def\ig{\includegraphics}\def\pdep{{\tt pde2path}}

\def\ds{\displaystyle}\def\del{\delta}\def\sig{\sigma}
\opt{hust}{
\renewcommand{\arraystretch}{1.1}\renewcommand{\baselinestretch}{1.05}
\addtolength{\belowcaptionskip}{-2mm}\addtolength{\abovecaptionskip}{-2mm}
}
\begin{document}
\author{Hannes Uecker, Daniel Wetzel 
\\ \small Institut f\"ur Mathematik, Universit\"at Oldenburg, D26111 Oldenburg
\\ \small \texttt{hannes.uecker@uni-oldenburg.de}\quad 
\small \texttt{daniel.wetzel@uni-oldenburg.de} 
}
\title{Numerical results for snaking of patterns over patterns 
in some 2D Selkov--Schnakenberg Reaction--Diffusion systems}
\maketitle
\begin{abstract}
  \noindent For a Selkov--Schnakenberg model as 
a prototype reaction-diffusion system  on 
  two dimensional domains we use the continuation and bifurcation 
software {\tt pde2path} to 
 numerically calculate branches of patterns embedded in patterns, for 
instance hexagons embedded in stripes and vice versa, with a planar 
interface between the two patterns. 
We use the 
Ginzburg-Landau reduction to approximate the locations 
of these branches by Maxwell points for the associated Ginzburg--Landau 
system.  For our basic model, some but not all of these 
branches show a snaking behaviour in parameter space, over the given 
computational domains. The (numerical) non--snaking 
behaviour appears to be related to too narrow bistable ranges with 
rather small Ginzburg-Landau energy differences. This claim 
is illustrated by a suitable generalized model. Besides the 
localized patterns with planar interfaces we also give a number of examples 
of fully localized patterns over patterns, for instance hexagon patches 
embedded in radial stripes, and fully localized hexagon patches over straight 
stripes. 
\end{abstract}

\tableofcontents 
\section{Introduction}  
Homoclinic snaking refers to the back and forth oscillation in parameter space 
of a branch of stationary localized patterns for some pattern forming 
partial differential equation (PDE). 
Two standard 
models are the quadratic--cubic Swift--Hohenberg equation (SHe) 
\huga{\label{qcshe}
\pa_t u=-(1+\Delta)^2u+\lam u+\nu u^2-u^3, \quad u=u(t,x)\in\R,\quad 
x\in\Om\subset\R^d, 
}
and the cubic--quintic SHe 
\huga{\label{cqshe}
\pa_t u=-(1+\Delta)^2u+\lam u+\nu u^3-u^5, \quad u=u(t,x)\in\R,\quad 
x\in\Om\subset\R^d, 
}
with suitable boundary conditions if $\Om\ne\R^d$, and 
where $\lam\in\R$ is the linear instability parameter, and $\nu>0$. 

In both equations the trivial solution $u\equiv 0$ is stable 
for $\lam<\lam_c:=0$, where in the 1D case $d=1$ we have 
a pitchfork bifurcation of periodic solutions with period near $2\pi$ 
(``stripes'', by trivial 
extension to 2D, also called ``rolls'' due to their connection to 
Rayleigh--B\'enard convection rolls). 
For $\nu>\nu_0\ge 0$, 
with $\nu_0=\sqrt{27/38}$ for \reff{qcshe} and $\nu_0=0$ for 
\reff{cqshe}, the bifurcation is subcritical and 
the periodic branch starts with unstable 
small solutions $r_-$ and turns around in a fold 
at $\lam=\lam_0(\nu)<0$ to yield 
$\CO(1)$ amplitude stable periodic solutions $r_+$. 
Thus, for $\lam_0<\lam<\lam_c$ there is a bistable regime 
of the trivial solution and $\CO(1)$ amplitude stripes. 

\begin{figure}[!ht]
\bce
{\small 
\begin{minipage}{60mm}(a)\\
\ig[width=55mm]{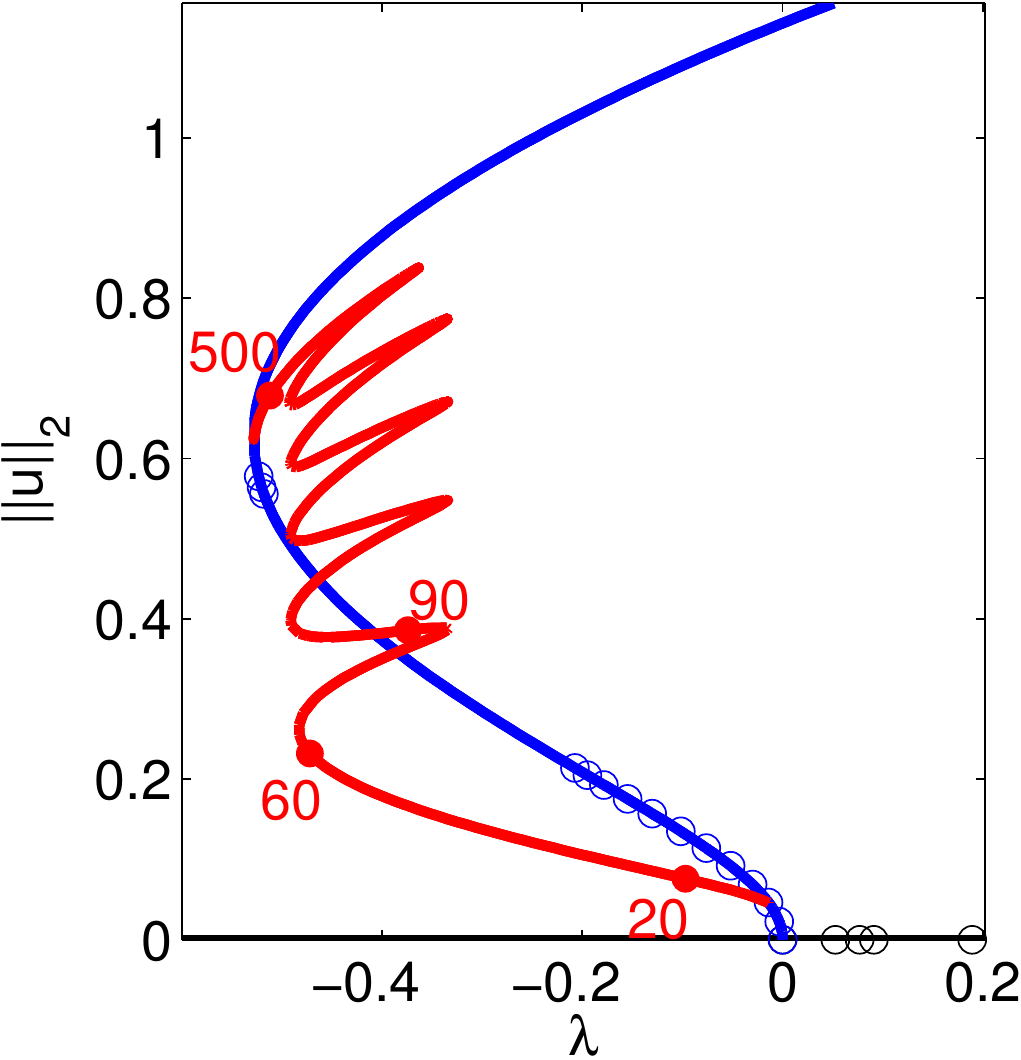}
\end{minipage}
\begin{minipage}{100mm}(b) $u$ at points 20,60,90 and 500 as indicated in (a)\\
\ig[width=100mm,height=13mm]{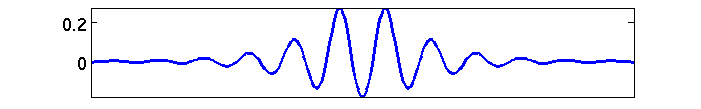}\\
\ig[width=100mm,height=13mm]{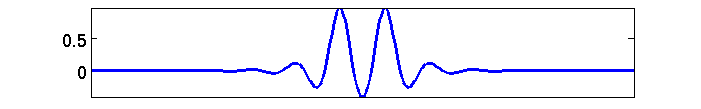}\\
\ig[width=100mm,height=13mm]{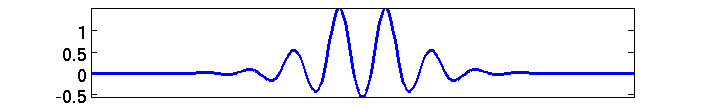}\\
\ig[width=100mm]{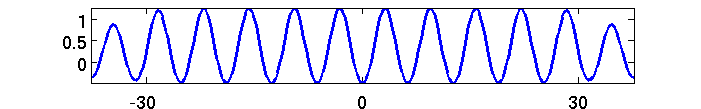}
\end{minipage}
}
\ece 

\vs{-5mm}
\caption{{\small Illustration of homoclinic snaking 
in \reff{qcshe} 
over a 1D domain $x\in(-12\pi,12\pi)$ with Neumann type boundary conditions 
$\pa_x u|_{x=\pm12\pi}=\pa_x^3 u|_{x=\pm12\pi}=0$, $\nu=2$. 
Here $\|u\|_2:=\left(\frac 1{|\Om|}\int_\Om 
|u(x,y)|^2{\rm d}(x,y)\right)^{1/2}$. 
The branch $r$ of primary roll solutions (blue) 
bifurcates subcritically at $\lam=0$ and becomes stable at $\lam\approx
-0.53$. Bifurcation points on the trivial branch and the roll branch 
are indicated by $\circ$, but omitted on the 
homoclinic snaking branch (red), which bifurcates 
from the second bifurcation point on $r$  and reconnects to $r$ 
just below the fold.\label{ipic1}}}
\end{figure}

In the 
simplest case the localized patterns then consist of 1D stripes over the 
homogeneous background $u=0$, and in each pair of turns in the snake 
the localized pattern grows by adding a stripe on both sides, 
and this continues for ever over the infinite line. 
See, e.g., \cite{burke, bukno2007, BKLS09} for seminal results in this setting, 
and \cite{chapk09,dean11} for detailed analysis using a Ginzburg--Landau 
formalism and beyond all order asymptotics. 
In finite domains snaking cannot continue for ever, and instead 
branches typically connect primary stripe branches, with either the same 
wave number or wave numbers near each other, see, e.g.,  \cite{hokno2009, 
bbkm2008,dawes1,KAC09} for detailed results, and Fig.~\ref{ipic1} for 
illustration. 

In 2D, by rotational invariance, depending on the domain and 
boundary conditions multiple patterns may bifurcate from $u\equiv 0$ 
at $\lam=0$, for instance straight stripes and so called 
hexagons. There are a number of studies of localized
patterns for \reff{qcshe} and \reff{cqshe} over two dimensional domains 
\cite{hexsnake,strsnake}, often combining analysis and numerics, 
and in more complicated systems like 
fluid convection \cite{kno2011}. See also \cite{lloyd13} for recent results 
on snaking of localized hexagon patches in a 2-component reaction--diffusion 
system in 2D. However, all these works essentially 
consider patterns over a homogeneous background. 
Already in \cite{pomeau} it is 
pointed out that ``pinned'' fronts connecting stripes and hexagons 
may exist in reaction-diffusion systems as a codimension 0 phenomenon 
in parameter space, i.e., for a whole interval of parameters. 
Similar ideas were also put forward 
in \cite{MNT90} with the 2D quadratic--cubic SHe  \reff{cqshe} 
as an example, see in particular \cite[Appendix C]{MNT90}. 
Such a pinned front is observed in \cite{hilaly} for 
\reff{cqshe} by time integration, but so far 
no studies of branches of stationary solutions involving 
different 2D patterns seem available. 

Here we start with a standard model problem for predator ($u$) prey ($v$) 
reaction diffusion systems, namely 
\huga{ \label{rdt} 
\pa_t U= D \Delta U+N(U,\lam), \qquad 
N(U,\lam)=\bpm-u+u^2v \\\lambda -u^2v\epm,
}
with $U=(u,v)(t,x,y)\in\R^2$, diffusion matrix 
$D=\ssmatrix{1 & 0\\ 0 & d}$, 
$d$ fixed to $d=60$, and  bifurcation parameter $\lambda\in \R_{+}$. 
The reaction term $N$ of \eqref{rdt}
is a special version of  $\ds 
N(U,\lam)=\bpm-u+u^2v+b+av \\\lambda -u^2v-av\epm$, known as 
the Selkov \cite{selkov} ($a\ge 0, b=0$) and Schnakenberg \cite{schnakenberg} 
 ($a=0, b\ge 0$) model. For simplicity, here we set $a=b=0$. 

In particular we consider the stationary system 
\begin{equation}\label{rd}
D \Delta U+N(U,\lam)=0. 
\end{equation} 
 The unique spatially homogeneous solution of \eqref{rd} is  
 $w^*=(\ld,1/\ld)$.  We write \eqref{rd} in the form
 \begin{align}
   \partial_t w = L(\Delta) w +G(w), \label{dglredu0}
 \end{align}
 with $w=U-w^*$ and $L(\Delta,\ld)=J(\ld)+D \Delta$, where $J$ is the 
Jacobian of $N$ in $w^*$. For $\fettk=(m,n)\in\R^2$ we
 have $L(\Delta,\ld) e^{i (\kx x+ \ky y)}=\hat{L}(\fettk,\ld)e^{i (\kx
   x+ \ky y)}$, where $\hat{L} (\fettk,\ld)=J(\ld)-Dk^2$, 
$k:=\sqrt{m^2+n^2}$, and thus 
we also write $\Lh(k,\lam)$. The eigenvalues of $\hat{L}$ are given by 
$\mu_\pm(\fettk,\ld)
=\mu_\pm(k,\ld)=\frac{\trh}{2} \pm \sqrt{\left(\frac{\trh}{2}\right)^2-\detlh}$. 
Following \cite[Chapter 14]{murray} we find that in
 $\lambda_c=\sqrt{d}\sqrt{3-\sqrt{8}}\approx 3.2085$ we have
 $\mu_+(\fettk,\ld_c)=0$ for all vectors $\fettk\in \R^2$ of  
 length $k_c=\sqrt{\sqrt2-1}\approx 0.6436$, and 
 all other $\mu_\pm (\fettk,\ld_c) < 0$, i.e., there is a Turing
 bifurcation at $\lambda_c$, with critical wave vectors $\fettk\in\R^2$ 
in the circle $|\fettk|=k_c$.  The most prominent Turing patterns 
near bifurcation are stripes and hexagons, which 
modulo rotational invariance can be expanded as 
\begin{align}\label{cos}
  \nonumber U&=w^*+ 2\lkl A\cos(k_c x)+ B\cos\lkl\frac{k_c}{2}\lkl-x+\sqrt{3}y\rkl\rkl+ B\cos\lkl\frac{k_c}{2}\lkl -x-\sqrt{3}y \rkl \rkl \rkl\Phi+\text{h.o.t.}\\
  &=w^*+ 2\lkl A\cos(k_cx)+ 2 B \cos\lkl\frac{k_c}{2}x\rkl
  \cos\lkl\frac{k_c}{2}\sqrt{3} y\rkl\rkl\Phi+\text{h.o.t.},
\end{align}
where $\Phi\in \R^2$ is
the critical eigenvector of $\Lh(k_c,\lam_c)$, and \text{h.o.t.} stands 
for higher order terms. 
The amplitudes $2A, 2B\in \R$ (where the factor 2 has been introduced 
for consistency with \S\ref{AnaRes}) of the corresponding
Turing pattern depend on $\lam$, with $A=B=0$ at bifurcation. 
The stripes bifurcate in a supercritical pitchfork, the hexagons 
bifurcate trancritically, and the subcritical part of the hexagons 
is usually called ``cold'' branch as here $u$ has a minimum in the center 
of each hexagon.   

\begin{figure}[!ht]
\bce
{\small 
\begin{minipage}{60mm}(a)\\
\ig[width=55mm]{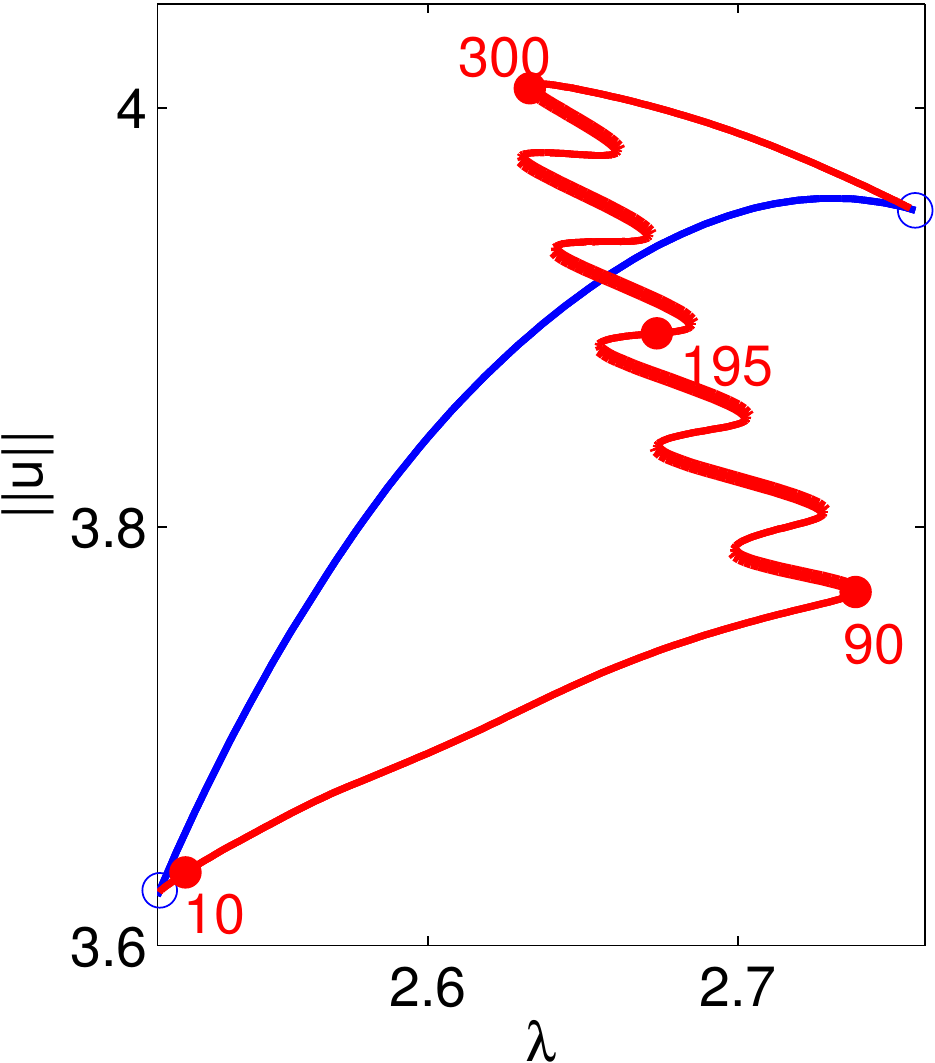}
\end{minipage}
\begin{minipage}{100mm}(b) $u$ at points 10,90,195 and 300 as indicated in (a)\\
\ig[width=100mm]{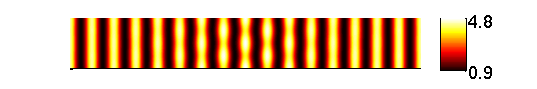}\\[-2mm]
\ig[width=100mm]{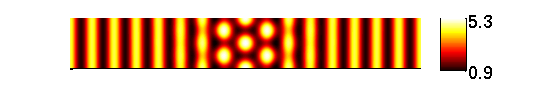}\\[-2mm]
\ig[width=100mm]{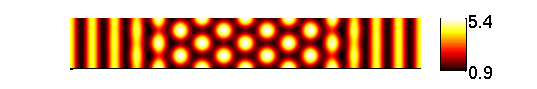}\\[-2mm]
\ig[width=100mm]{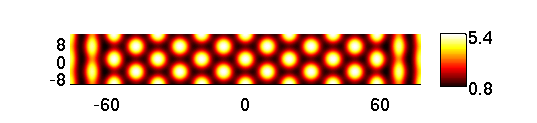}
\end{minipage}

\begin{tabular}{lll}
(c)&(d)&(e)\\[-0mm]
\ig[width=50mm]{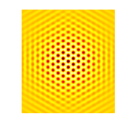}&
\ig[width=50mm]{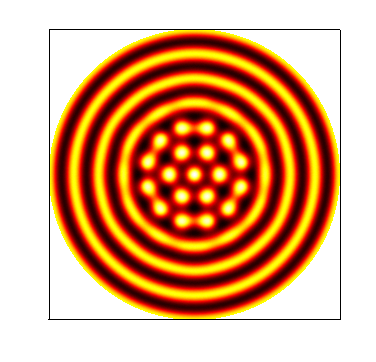}&
\ig[width=50mm]{./BB/fl200}
\end{tabular}
}
\ece 

\vs{-5mm}
\caption{{\small (a), (b) Mixed mode branch (blue) in the RD system \reff{rd} 
which connects stripes and hexagons, and a
branch of hexagons on a background of stripes. Here 
$\norm{u}=\norm{u}_{L^8}:=\left(\frac 1{|\Om|}\int_\Om 
|u(x,y)|^8{\rm d}(x,y)\right)^{1/8}$ was chosen for purely graphical reasons.  
(In Fig.~\ref{ipic1} any $L^p$ norm can be used to illustrate the snaking, 
but for \reff{rd} the stripes and hexagons have rather similar $L^2$ norms.) 
This example will be discussed in much more detail below, and here we 
mainly want to point out the similarities with Fig.~\ref{ipic1}. 
In Fig.~\ref{ipic1} the snaking branch is in the bistable range of (large) 
amplitude 
rolls and the homogeneous background, here it is in the 
bistable range of stripes and hexagons. In both cases the snaking 
branch bifurcates from an unstable branch connecting two different 
stable branches, and after bifurcation first needs a long transverse 
to enter the ``snaking region''. Finally, in both cases, the snaking 
branches consist of stable solutions after every second fold, 
indicated here by thicker lines. \newline
(c)-(e) $u$ for fully localized spot patterns 
over homogeneous background (c), and over radial (d) and straight (e) stripes. 
See below for colorbars and spatial scales. 
(d) and (e) again lie on some snaking branches, see \S\ref{radsec}
and \S\ref{flsec}, while the branch for (c) numerically does not snake, see 
\S\ref{crsec}. However, a snaking branch of patterns as in (c) can be obtained 
for the modified system \reff{mod1} below, see \S\ref{sigsec}. 
}  \label{ipic2}}
\end{figure}

In \S\ref{numrel} and \S\ref{addsec} we present some numerically calculated 
Turing patterns for \reff{rd}, including so called mixed modes, and, 
moreover, some 
branches of stationary solutions which involve different patterns, namely 
\bci
\item[(a)] planar fronts between stripes and hexagons, and associated localized 
patterns, e.g., hexagons localized in one direction 
on a background of stripes (see Fig.\ref{ipic2}(a),(b) for an example), 
and vice versa; 
\item[(b)] fully 2D localized patches of hexagons over a homogeneous 
background, and over radial and straight stripes 
(see Fig.\ref{ipic2}(c),(d),(e) for examples).
\eci 

Related to the transcritical bifurcation of the hexagons, the patterns in 
(a) come in two regimes: one ``hot'' which is rather 
far from the primary bifurcation at $\lam_c$; one ``cold'', with 
$\lam$ closer to $\lam_c$. Of these, only the ``hot'' branches show snaking  
behaviour in our numerical simulations for \reff{rd}. 
For \reff{qcshe} and \reff{cqshe} it is known (e.g., 
\cite{KC06}, \cite{dean11}), that the ``snaking width'' is exponentially 
small in the properly defined subcriticality parameter $\eps$. We expect 
a similar result here, which explains why snaking in the cold regimes might 
be below the numerical resolution. To make this quantitative we modify 
\reff{rdt} to make the  system ``more subcritical'' in the cold regime, 
and in this way 
we can switch on cold snaking, i.e., find it numerically. 

First, however, in \S\ref{AnaRes} we relate our numerical results to arguments
derived from the Ginzburg--Landau reduction. One result is the 
calculation of Ginzburg--Landau Maxwell points, which give
predictions for the $\lam$--ranges of the fronts and 
localized patterns in the full system \reff{rd}. In fact, for 
a consistent Ginzburg--Landau reduction of \reff{rdt} 
in principle one should include the parameters $a$,$b$ from $\ds 
N(U,\lam)=\bpm-u+u^2v+b+av \\\lambda -u^2v-av\epm$, to locate a co--dimension 
2 point where the quadratic terms in \reff{dglredu0} scale suitably. 
However, then the algebra for the bifurcation analysis will become even 
more involved. Moreover,  
here we are explicitly interested in phenomena $\CO(1)$ away 
from $(w,\lam){=}(0,\lam_c)$, and thus rather explore how far the 
Ginzburg--Landau reduction can take us.  

Then, motivated by the 
Ginzburg--Landau reduction, we modify the system \reff{rdt} to 
\huga{\label{mod1}\pa_t U=
D\Delta U+N(u,\lam)+\sigma \left(u-\frac 1 v\right)^2\bpm 1\\-1\epm. 
}
We do not claim any biological meaning for the term multiplied by $\sigma$, 
but the advantage of this modification is that the homogeneous solution and 
the linearization around it 
are unchanged, the Turing bifurcation is still at $\lam_c$, 
and only the nonlinear terms in 
\reff{dglredu0} are affected. Our first goal is to increase 
the subcriticality of the cold hexagons to obtain homoclinic snaking 
for branches connecting these with $w=0$. Similarly, via 
$\sigma$ we can increase the bistability range between cold stripes 
and cold hexagons and thus we can also switch on snaking in this regime. 
Finally, while \reff{rd} has no bistability between stripes 
and $w=0$ and thus no homoclinic snaking in 1D can be expected, 
in \reff{mod1} we can also turn the bifurcation of stripes from 
supercritical to subcritical to obtain 1D homoclinic snaking 
between stripes and $w=0$, see \S\ref{1dsec}. 

Sections \ref{addsec} contains additional numerical results 
about fully localized patterns over patterns, namely (hexagonal) spots 
over radial stripes, and fullly localized (hexagonal) spots over 
straight stripes, respectively, see also Fig.~\ref{ipic2}(d),(e). 
The numerics 
for such patterns are quite demanding, and therefore we restrict to 
some illustrative examples.

Thus, the remainder of this paper is organized in a slightly 
unconventional way, which is partly due to 
our own ``experimental mathematics'' 
approach, as newcomers to the field of homoclinic snaking: First we 
numerically calculated 
the Turing bifurcation diagram and spotted the hot snaking branches 
and the non-snaking cold localized branches, 
then we started the Ginzburg--Landau reduction to {\em a posteriori} 
get some analytical understanding, mainly via Maxwell points. 
Then, however, we used 
the Ginzburg--Landau reduction as a {\em predictive tool} to find 
ranges for the modified 
systems where we could switch on the cold snaking, and the 1D snaking 
via subcritical bifurcations of stripes. We believe that it is more 
honest to also present our results in that order, instead of trying 
to first present as much theory as possible and then use numerics 
only for illustration. More importantly, we believe 
that this makes the paper also more readable.

In \S\ref{dsec} we give with a brief discussion. We believe that the 
results are not specific to the model system \reff{rd} but 
can be expected in any reaction diffusion system (over sufficiently large 
domains) with a bistability 
between different patterns which allow homo--or heteroclinic connections 
in the associated Ginzburg--Landau system. Finally, we close with 
a short list of Open Problems.

\noindent
{\bf Acknowledgements.} We thank Jonathan Dawes for helpful discussions, and 
two anonymous referees for a critical 
reading of the first version of this manuscript and many helpful questions and 
remarks, which in particular motivated us to set up the modified system 
\reff{mod1}, and some of the numerics in \S\ref{addsec}.

\section{Numerical Results}\label{numrel} 
\subsection{Stripes, hexagons, and beans} 
\label{numisec}
We  use the bifurcation and
continuation software {\tt pde2path} \cite{p2p} to  numerically 
calculate patterns 
for \reff{rd}, and their stability. We start with 
domains of type $\Omega=(-l_x,l_x)\times (-l_y,l_y)$, $l_x=2l_1\pi / k_c$,
$l_y=2l_2 \pi/(\sqrt{3}k_c)$, $l_1,l_2\in\N$, with Neumann boundary conditions, 
chosen to accomodate the basic stripe and hexagon 
patterns.\footnote{In \cite[\S4.2]{p2p} we set $l_y=2\del l_2\pi/(\sqrt{3}k_c)$, 
where  the slight ``detuning'' $\delta\approx 1$ is used 
to unfold the multiple bifurcation at $\lam=\lam_c$ 
since {\tt pde2path} currently only deals with simple bifurcations. 
However, by carefully choosing $\lam$ stepsizes for $\lam$ 
near $\lam_c$ the discretization 
error is enough to do this unfolding and thus we drop $\delta$ here, i.e., 
set $\del=1$. 
The bifurcation diagram and solution plots for $\del=1$ and $\del=0.99$ 
are visually indistinguishable. }

\begin{figure}[!ht]
\bce

\begin{minipage}{0.4\textwidth}
\includegraphics[width=\textwidth, height=45mm]{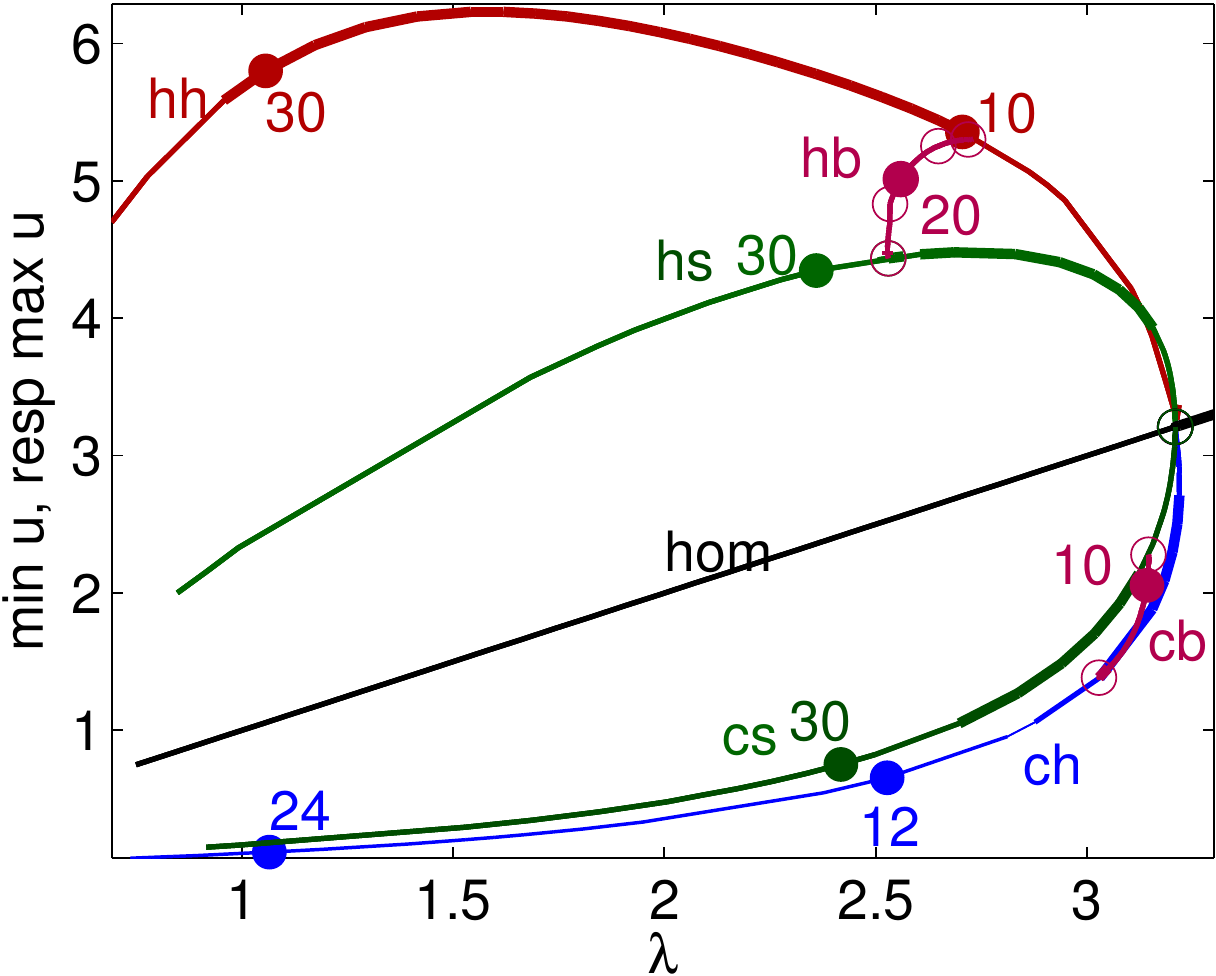}
\end{minipage}\hs{6mm}
\begin{minipage}{0.21\textwidth}
\includegraphics[width=\textwidth]{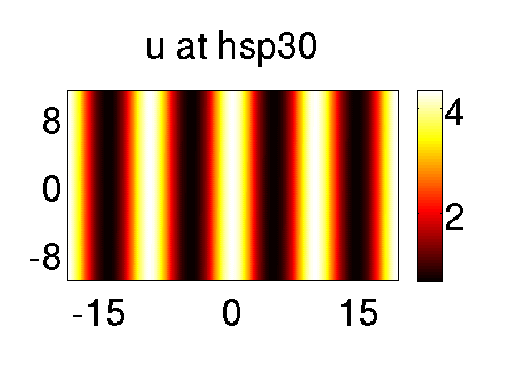}\\
\includegraphics[width=\textwidth]{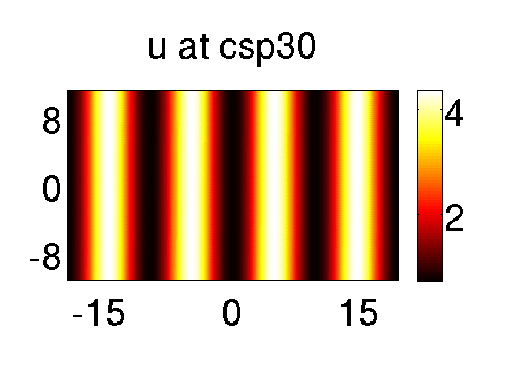}\\
\end{minipage}
\begin{minipage}{0.21\textwidth}
\includegraphics[width=\textwidth]{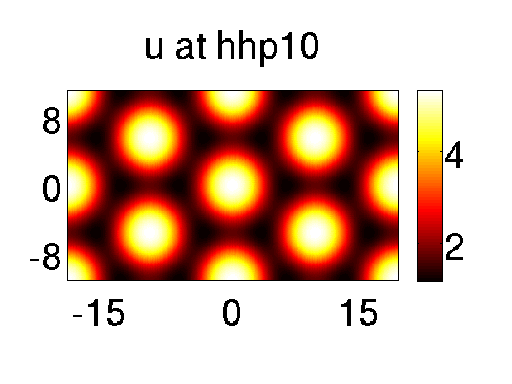}\\
\includegraphics[width=\textwidth]{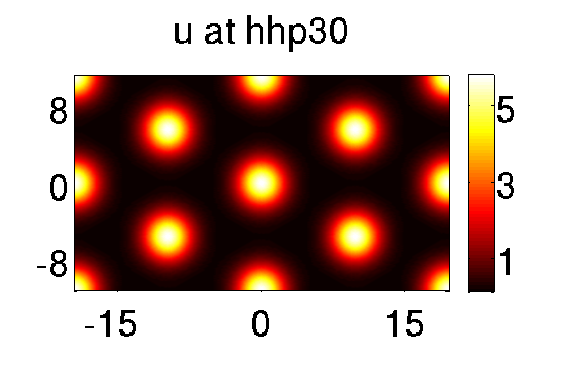}\\
\end{minipage}

\vs{-6mm}
\begin{minipage}{0.21\textwidth}
\includegraphics[width=\textwidth]{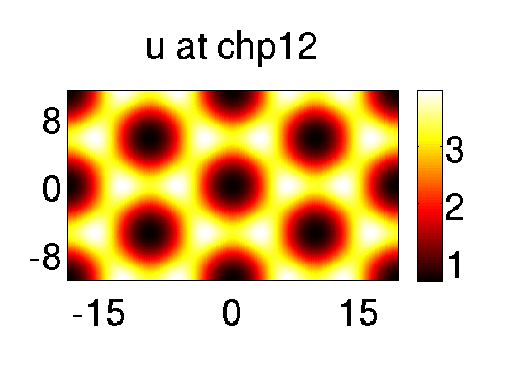}
\end{minipage}
\begin{minipage}{0.21\textwidth}
\includegraphics[width=\textwidth]{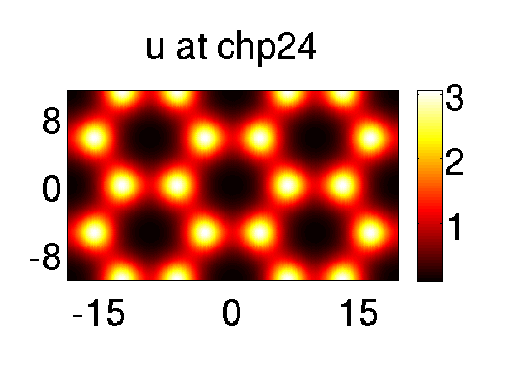}
\end{minipage}
\begin{minipage}{0.21\textwidth}
\includegraphics[width=\textwidth]{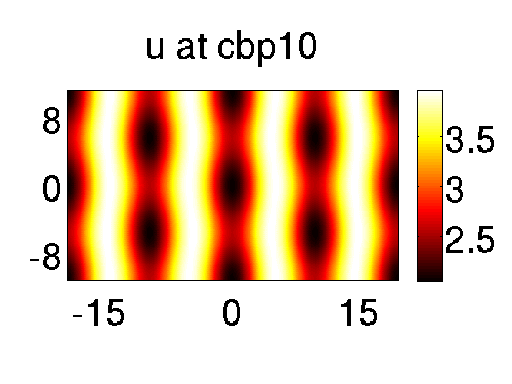}
\end{minipage}
\begin{minipage}{0.21\textwidth}
\includegraphics[width=\textwidth]{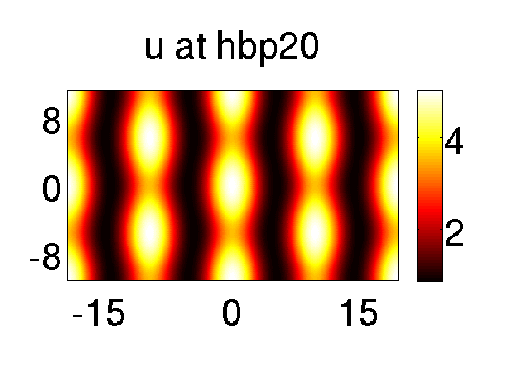}
\end{minipage}
\ece 

\vs{-6mm}
\caption{{\small Bifurcation diagram and example solutions of \eqref{rd}
  obtained with {\tt pde2path} over the ``$2\times 2$''--domain
  $\Omega=(-l_x,l_x){\times}(-l_y,l_y)$ with Neumann boundary
  conditions, $l_x{=}4\pi / k_c$, $l_y{=}4 \pi /(\sqrt{3}k_c)$, 
$k_c=\sqrt{\sqrt{2}{-}1}{\approx}0.6436$. 
The branch {\tt hom} are 
  the homogeneous solutions, {\tt hs} the hot (up) stripes,
  {\tt cs} the cold (down) stripes, {\tt hh} the hot hexagons, and {\tt ch} the
  cold hexagons. For instance {\tt hsp30} stands for the 30$^{\text{th}}$ point 
on the branch {\tt hs}. {\tt hb} and {\tt cb} are mixed mode branches, 
which we call hot
  and cold beans, respectively. For {\tt
    hs}, {\tt hh}, and {\tt hb} we  plot the maximum of $u$, 
and the minimum for {\tt cs}, {\tt ch},
  and {\tt cb}. Stable and
  unstable parts of branches are represented by thick and thin lines,
  respectively. In particular, except for {\tt hhp10, hhp30} 
all displayed patterns are unstable. 
Only a small selection of bifurcation points is indicated by $\circ$. 
See \cite{smov} for a movie.\label{figsch}}}
\end{figure}

In Fig.~\ref{figsch} 
we use the rather small $l_1=l_2=2$, which we call 
a $2\times 2$ domain as in both directions 2 hexagons ``fit''. 
The main panel shows a part of the very rich bifurcation diagram. 
Following biological terminology we classify the stripes ($A\ne 0$, $B=0$ in 
\reff{cos}) 
into hot (also called ``up'') stripes 
({\tt hs}, $A>0$) and
cold (also called ``down'') 
stripes ({\tt cs}, $A<0$), which exhibit a maximum respectively a 
minimum at $x=0$. We always plot $u$. If we have a
hot pattern for $u$, then $v$ is a cold pattern and vice versa. This
follows from the predator-prey structure of the reaction term $N$ of
\eqref{rd}. The stripes are $2 \pi/k_c$-periodic in the horizontal direction, 
bifurcate in a supercritical pitchfork, and are stable from $\lambda_{s}^b
\approx 3.15$ to $\lambda_{s}^e \approx 2.51$, where here and henceforth the 
subscript of the bifurcation parameter $\lambda$ denotes the branch
  and the superscript stands for \textbf{e}nding or \textbf{b}eginning of
  stability.  
The hexagons ($A=B\ne 0$) can be classified
into hot hexagons ({\tt hh}, $A=B>0$) and cold hexagons ({\tt ch}, $A=B<0$), 
which have a maximum resp.~a minimum in the center of every hexagon. 
They are $4 \pi/k_c$-periodic
in the horizontal and $4 \pi/(\sqrt{3}k_c)$-periodic in the vertical
direction. The hexagon branch bifurcates transcritically from the homogeneous 
branch at $\lam=\lam_c$. It ``starts'' in a saddle node or fold bifurcation at 
$\lambda_{ch}^b \approx 3.22$, and 
the stability region of the cold hexagons begins at the fold
and ends in $\lambda_{ch}^e\approx 3.03$.
The hot hexagons are stable from
$\lambda_{hh}^b \approx 2.73$ to $\lambda_{hh}^e \approx 0.98$.

Thus, there is a bistable range of cold stripes and cold hexagons for $\lambda
\in (\lambda_{ch}^e,\lambda_{s}^b)$, of hot stripes and hot hexagons for
$\lambda \in (\lambda_{s}^e,\lambda_{hh}^b)$, and of cold hexagons and the
homogeneous solution for $\lambda \in (\lambda_{c},\lambda_{ch}^b)$. 
A branch of ``skewed hexagon'' or ``mixed mode'' solutions 
of the form \eqref{cos} with $A<B<0$ bifurcates subcritically from 
the cold hexagon branch {\tt ch} in $\lambda_{ch}^e$ and terminates on the
cold stripe branch {\tt cs} in $\lambda_{s}^b$. Following \cite{yang2004} 
we call this type of solutions cold {\em beans} ({\tt cb}). 
There is also a branch {\tt hb} of hot beans with $A>B>0$ 
which bifurcates subcritically from the
hot stripe branch {\tt hs} in $\lambda_{s}^e$ and terminates on the hot
hexagon branch {\tt hh} in $\lambda_{hh}^b$.

\begin{figure}[!ht]
\bce
\begin{minipage}{70mm}(a)\\[-5mm]
\ig[width=65mm]{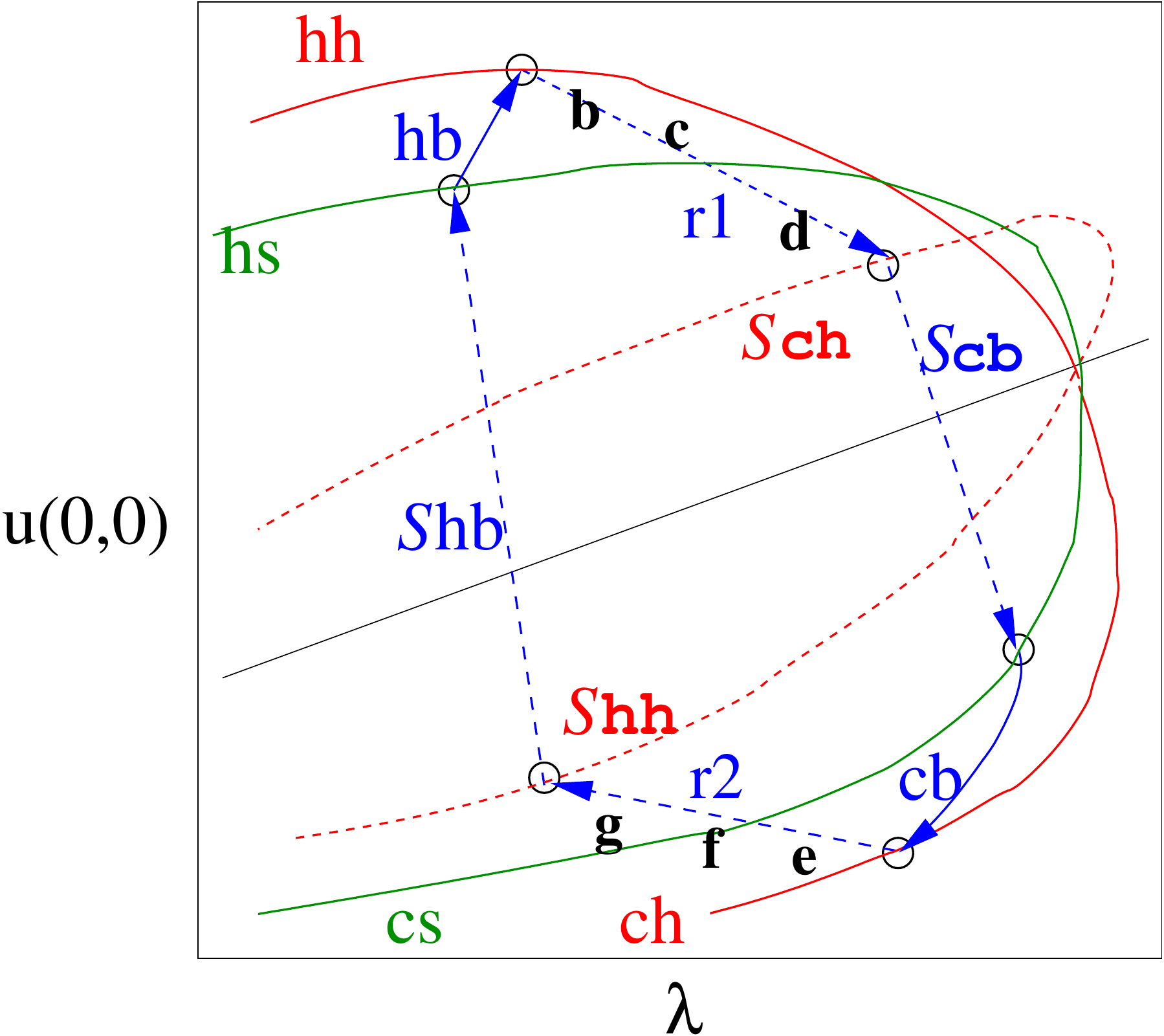}
\end{minipage}
\begin{minipage}{70mm}
\begin{tabular}{p{1mm}>{\vspace{-.7cm}}p{30mm}p{1mm}>{\vspace{-.7cm}}p{30mm}}
(b)&\ig[width=30mm]{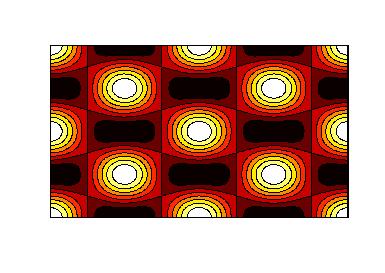}&
(c)&\ig[width=30mm]{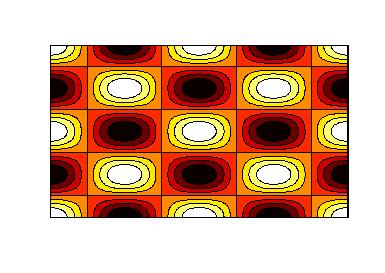}\\
(d)&\ig[width=30mm]{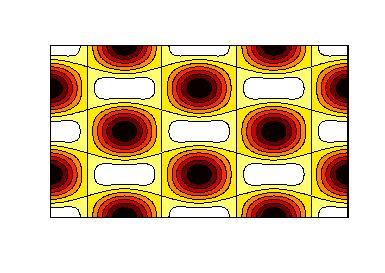}&
(e)&\ig[width=30mm]{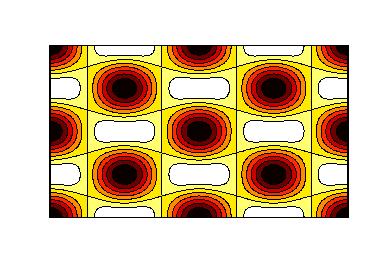}\\
(f)&\ig[width=30mm]{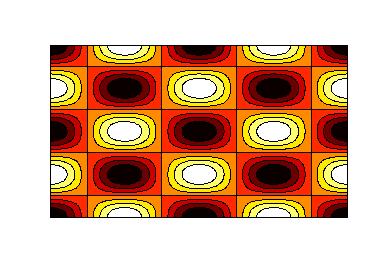}&
(g)&\ig[width=30mm]{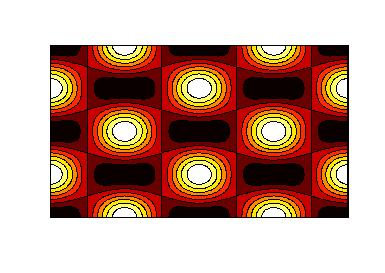}
\end{tabular}
\end{minipage}
\ece 

\vs{-5mm}
\caption{{\small (a) Schematic bifurcation diagram  with $u(0,0)$ as ordinate. Additional branches 
compared to Fig.~\ref{figsch} are dashed, and $\CS{\tt *}$ denotes the 
respective phase shifted pattern, i.e., $\CS{\tt hh}$ are 
the phase shifted hot hexagons.  
{\tt hb} and {\tt cb} both form part of a loop 
${\tt hs}\stackrel{{\tt hb}}{\ra}{\tt hh}\stackrel{{\tt r1}}{\ra}
\CS{\tt ch}\stackrel{\CS{\tt cb}}{\ra}{\tt cs}
\stackrel{{\tt cb}}{\ra}{\tt ch}
\stackrel{{\tt r2}}{\ra}\CS{\tt hh}
\stackrel{\CS{\tt hb}}{\ra}{\tt hs}$, where 
{\tt r1} and {\tt r2} are called rectangles.  (b)--(d) {\tt r1} 
example solutions for which, in the expansion \reff{cos}, 
$B>0$ and $A$ changing from 
$A{>}0$ in (b) via $A{=}0$ in (c) to $A{<}0$ in (d). 
(e)--(f) {\tt r2} with $B<0$ 
and $A$ changing from $A{<}0$ in (e) via $A{=}0$ in (f) to $A{>}0$ in (g). 
(b)--(d) correspond to ``hot'' solutions with a maximum in the middle 
of the domain and hence in the middle of the central hexagon, 
while (e)--(g) illustrate ``cold'' solutions with a minimum in the 
middle of the domain. 
}  \label{huf2}}
\end{figure}

\brem\label{srem} For all these patterns the discrete symmetry 
$\CS=S_{2\pi/k_c}$ given by the shift by $2\pi / k_c$ in the horizontal
direction yields a solution as well. The stripes are of course 
invariant under $\CS$, but for the hexagons the shifts 
generate new branches $ \CS{\tt hh}$ 
and $\CS{\tt ch}$. These make the plotting 
of the bifurcation diagram  a bit complicated graphically. 
Therefore, in Fig.~\ref{huf2} we repeat the bifurcation diagram  from Fig.~\ref{figsch} 
in a schematic way, with ordinate $u(0,0)$, i.e., $u$ in the center of the  
domain, and add the branches of shifted hexagons. 
The bean branches {\tt hb} and {\tt cb} then both take part in a ``loop'' 
involving shifted beans and so called 
rectangles, which again can be expanded as in \reff{cos}, now with 
$|A|\le |B|$, but which are generically unstable. 
This can be worked out on the level of amplitude equations, 
see \S\ref{AnaRes}, and it is also recovered by our numerics for 
the full system. 
The behaviour described below like snaking branches bifurcating from 
hot beans obviously transfers to branches related by $\CS$. 
In the following we mainly focus on the {\tt hb} and {\tt cb} branches, and 
also mostly restrict to following just one direction at bifurcation. 
\eex 
\erem 

\brem The main regime of interest to us is $\lam\in [2.5,3.22]$, 
and below $\lam\approx 2.5$ all branches plotted except 
${\tt hh}$ are unstable. 
We plot a somewhat larger bifurcation picture since, e.g., on the 
{\tt ch} branch some interesting patterns occur. \eex 
\erem 

\brem\label{numrem} (a) {\tt pde2path} uses the Matlab FEM {\tt pdetoolbox} 
to discretize elliptic PDEs like \reff{rd}, including 
some error estimators and adaptive mesh--refinement, 
see \cite{p2p} for details. 
In Fig.~\ref{figsch} we used a regular ``base mesh'' of 20.000 triangles 
which, e.g., on the beans branch is refined to about 60.000 triangles 
on average.  Moreover, for all solutions 
calculated the mesh can be further refined to yield arbitrary small error 
estimates without visible changes to the solutions. 
Calculation time for the full bifurcation diagram  in Fig.~\ref{figsch} 
on a quad core desktop PC is about 40 Minutes. 
The {\tt pde2path} script to generate 
Fig.~\ref{figsch} is also included in the software as demo {\tt schnakenberg},  
and a movie running through the bifurcation diagram of Fig.~\ref{figsch} 
and some more movies are collected at 
{\small {\tt www.staff.uni-oldenburg.de/hannes.uecker/pde2path/schnakmov.html}}. 

(b) The rather large numbers of triangles was needed mainly 
in order to avoid uncontrolled branch switching, and it is the 
large number of different branches which makes the continuation and bifurcation 
numerics of \reff{rd} demanding. For instance, on the $2\times 2$ domain 
the first 10 bifurcations from the homogeneous branch occur 
between $\lam=\lam_c\approx 3.2085$ and $\lam_{10}\approx 
3.1651$. On larger domains, 
these bifurcation points collapse to $\lam_c$. 
For instance, over the $6\times 2$ domain used in Fig.~\ref{huf3} below we have 
$\lam_{10}\approx 3.2042$. 
Similarly, on all branches shown in Fig.~\ref{figsch} 
there are many bifurcation points, and thus we only plot a small 
part of the bifurcation diagram. See \cite{p2p} on details how {\tt pde2path} tries to avoid branch jumping. 

(c) The branches plotted 
in Fig.~\ref{figsch} stay the same over all $l_1\times l_2$ domains 
with $l_1,l_2\in\N$, $l_1,l_2\ge 2$. 
However, stability of solutions of course may change. 
For instance, taking a larger $l_2$, the hot stripes no longer 
loose stability at bifurcation of the {\em regular} beans, i.e., with 
fundamental 
wave-vectors $k_1=k_c(1,0), k_{2,3}=k_c(-\frac 1 2, \pm \frac{\sqrt{3}}{2})$ 
forming an equilateral triangle. Instead, 
 {\em streched} beans with sideband wave vectors 
$\tilde k_{2,3}\approx k_{2,3} $, $\tilde k_{2,3}\neq 
k_{2,3} $, bifurcate more early, i.e., for larger $\lam$. 
See, e.g., \S\ref{flsec}. In the following 
we first keep $l_2=2$ fixed and consider larger $l_1$, which will 
give interfaces between stripes and hexagons parallel to the 
(vertical) stripes. In \S\ref{addsec} we look at some other 
interface angles between, e.g., hexagons and stripes.  
This, however remains only an outlook on the complicated problem of 
interface orientation. 

(d) To use {\tt pde2path} for 1D problems like in Fig.~\ref{ipic1}, see 
also \S\ref{1dsec}, we artificially set up a thin strip 
and use Neumann boundary conditions in transverse (here $y$) 
direction such that solutions are constant in $y$, e.g., 
$\Om=(-12\pi, 12\pi)\times (l_y,l_y)$ with small $l_y=0.1$ in Fig.~\ref{ipic1}. 
To keep computations cheap but accurate, $l_y$ and 
the initial mesh are chosen in such a way that we have only two grid points 
in $y$ direction, and that the triangles are roughly equilateral. 
The number of ``grid points in $y$--direction'' may then 
change under mesh-refinement, but the mesh quality is controlled automatically.  

(e) Finally, the current version of {\tt pde2path} as such 
only detects and deals with simple 
bifurcations in the sense that a simple eigenvalue goes through zero. 
Thus, no Hopf bifurcations are detected. By counting the number of 
negative eigenvalues we can however check for Hopf bifurcations 
a posteriori, see  \S\ref{1dsec} for an example.  
\eex
\erem

 \subsection{Planar fronts, localized patterns and snaking in the hot bistable range}
\label{lpsec}
To reduce computational costs, but also for graphical reasons, in 
Fig.~\ref{figsch} we used a rather small domain for the basic 
bifurcation diagram. We now present simulations on larger domains, 
where however we still restrict to intermediate sizes. Again, this is mainly 
due to graphical reasons, and we remark that all results can be reproduced 
 on larger domains, at least qualitatively.

\begin{figure}[!ht]
\bce 
{\small
\begin{minipage}{60mm}(a) Bifurcations from the bean branch\\[2mm]
\ig[width=60mm,height=65mm]{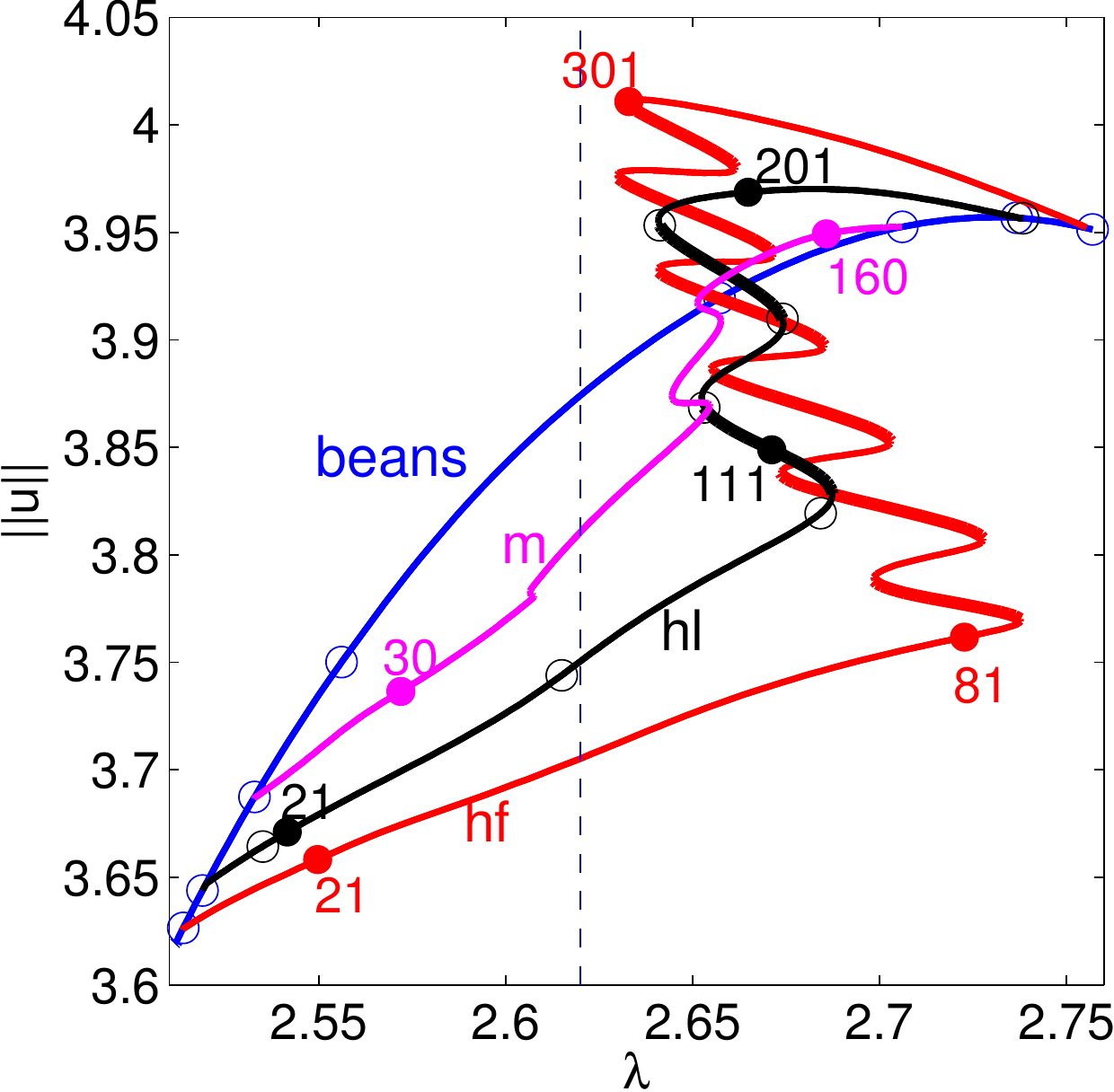}\\
\end{minipage}
\begin{minipage}{50mm}
(b)  hf21, hf81, hf301, hl11\\
\text{}\hs{3mm}\ig[width=50mm]{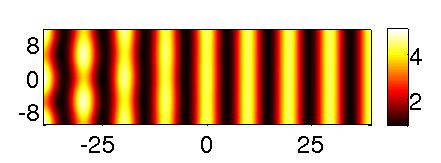}\\
\ig[width=57mm]{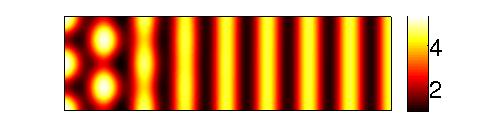}\\
\ig[width=57mm]{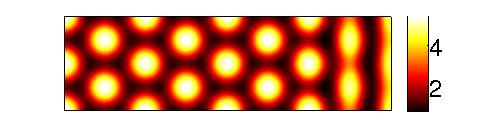}\\
\ig[width=57mm]{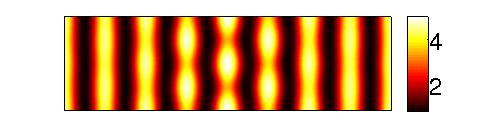}
\end{minipage}\hs{2mm}
\begin{minipage}{50mm}
(c)  hl111, hl201, m30, m160\\
\ig[width=50mm]{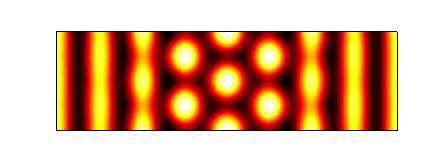}\\[-2mm]
\ig[width=50mm]{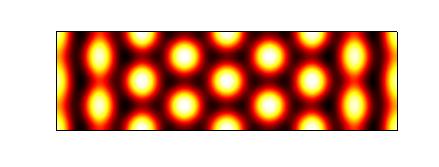}\\[-4mm]
\ig[width=50mm]{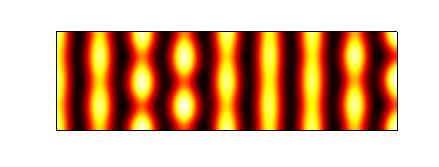}\\[-4mm]
\ig[width=50mm]{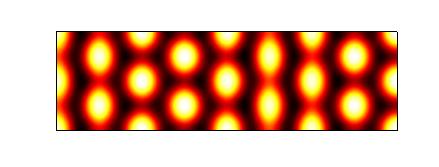}
\end{minipage}
}
\ece 

\vspace{-5mm}
\caption{{\small Bifurcation diagram and example plots of 
fronts and localized patterns of \eqref{rd} 
obtained with {\tt pde2path} over a $4\times 2$ domain; 
$\ld_m\approx 2.62$ is the so called (Ginzburg--Landau) 
Maxwell point, see \S\ref{AnaRes}.  
In (b), for instance {\tt hf21} means the 21st point on the 
``hot front'' branch in (a), and similar for hl* and c*. 
Domain in all plots as in the first plot in (b), and the colormap 
is roughly constant. 
Bifurcation points (in (a) only shown on the bean and the hl branches) 
are indicated by $\circ$, stable and
  unstable parts of branches by thick and thin lines,
  respectively. See \cite{smov} for a movie. 
} \label{snakea}}
\end{figure}

Counting from the {\tt hs} branch, 
{\tt pde2path} yields four bifurcation points {\tt hbbp1}, {\tt hbbp2}, 
{\tt hbbp3}
 and {\tt hbbp4} 
on the hot bean branch {\tt hb}. One bifurcating branch connects 
{\tt hbbp1} and {\tt hbbp4}, and another branch connects 
{\tt hbbp2} and {\tt hbbp3}. 
 By doubling the horizontal length to a $4\times 2$ domain, 
i.e., setting $l_x=8\pi / k_c$, we find eight 
bifurcation points {\tt hbbp1},...,{\tt hbbp8} on {\tt hb}, 
see Fig.~\ref{snakea}, where now as in Fig.\ref{ipic2} we use 
\huga{
\norm{u}=\norm{u}_{L^8}:=\left(\frac 1{|\Om|}\int_\Om 
|u(x,y)|^8{\rm d}(x,y)\right)^{1/8}
}
as the main solution measure. This choice is quite arbitrary, 
but by trial and error we found it 
more suitable to display branches snaking between different 
patterns than the usual $L^2$ norm, which works well for 
snaking of patterns over some homogeneous background. 

The branches bifurcating in {\tt hbbp1} and {\tt
   hbbp2} are called {\tt hf} (hot front) and {\tt hl} (hot localized), 
respectively. Some
 example solutions on {\tt hf} and {\tt hl} are also presented, 
and the so called Ginzburg--Landau Maxell point $\lam_m$. 
The {\tt hf} branch connects  {\tt hbbp1} to 
 {\tt hbbp8}, and contains stationary fronts $\Uhet$ 
 from hot hexagons to hot stripes, while {\tt hl} 
connects {\tt hbbp2} to {\tt hbbp7} and contains homoclinic solutions 
$\Uhom$ in the form of localized hexagons embedded in stripes. 
More precisely, the branches 
take part in closed loops, e.g., the {\tt hf} branch is one half 
of a loop containing 
{\tt hbbp1} and {\tt hbbp8}, the other half containing hexagons 
coming into the domain from the right. 
Again we generally only discuss parts of each 
of these closed loops. 

Both branches, {\tt hf} and {\tt hl}, show a snaking behaviour in the bifurcation diagram, 
and we first discuss the {\tt hl} branch which 
indicates the start of so called homoclinic 
snaking, which becomes more prominent if we further increase the domain size, 
cf.~Fig.\ref{ipic2}. On {\tt hl}, 
looking from {\tt hbbp2} to {\tt hbbp7}, during each 
cycle consisting of two folds, a further hexagon is added on both sides 
of the hexagon patch localized in the middle over a background of 
stripes. The parts of the snake pointing north--west contain stable solutions, 
while on the parts pointing north--east there are 2 unstable eigenvalues. 
The transitions from unstable to stable parts thus proceed via the folds and 
additional bifurcation points, discussed below. 

In contrast to theory and also to 1D problems over very 
large domains, the branch does not snake around 
a vertical line but in a slanted manner. 
This, and the fact that the branch does not directly bifurcate 
from the stripes but from the beans, are finite 
size effects, cf.~\cite{slant2, bbkm2008,hokno2009,  dawes1}. 
During, e.g., the initial traverse from {\tt hbbp2} to the third 
bifurcation point the  
``near stripe bean pattern'' at the bifurcation point is reshaped 
to hexagons in the middle and stripes at the sides. The analogous reshaping  
to ``near hexagon beans'' 
takes place between {\tt hlp201} and {\tt hbbp7}. 

\begin{figure}[!ht]
\bce 
{\small
\begin{minipage}{60mm}(d) Bifurcations from the hl branch.\\
\ig[width=60mm, height=50mm]{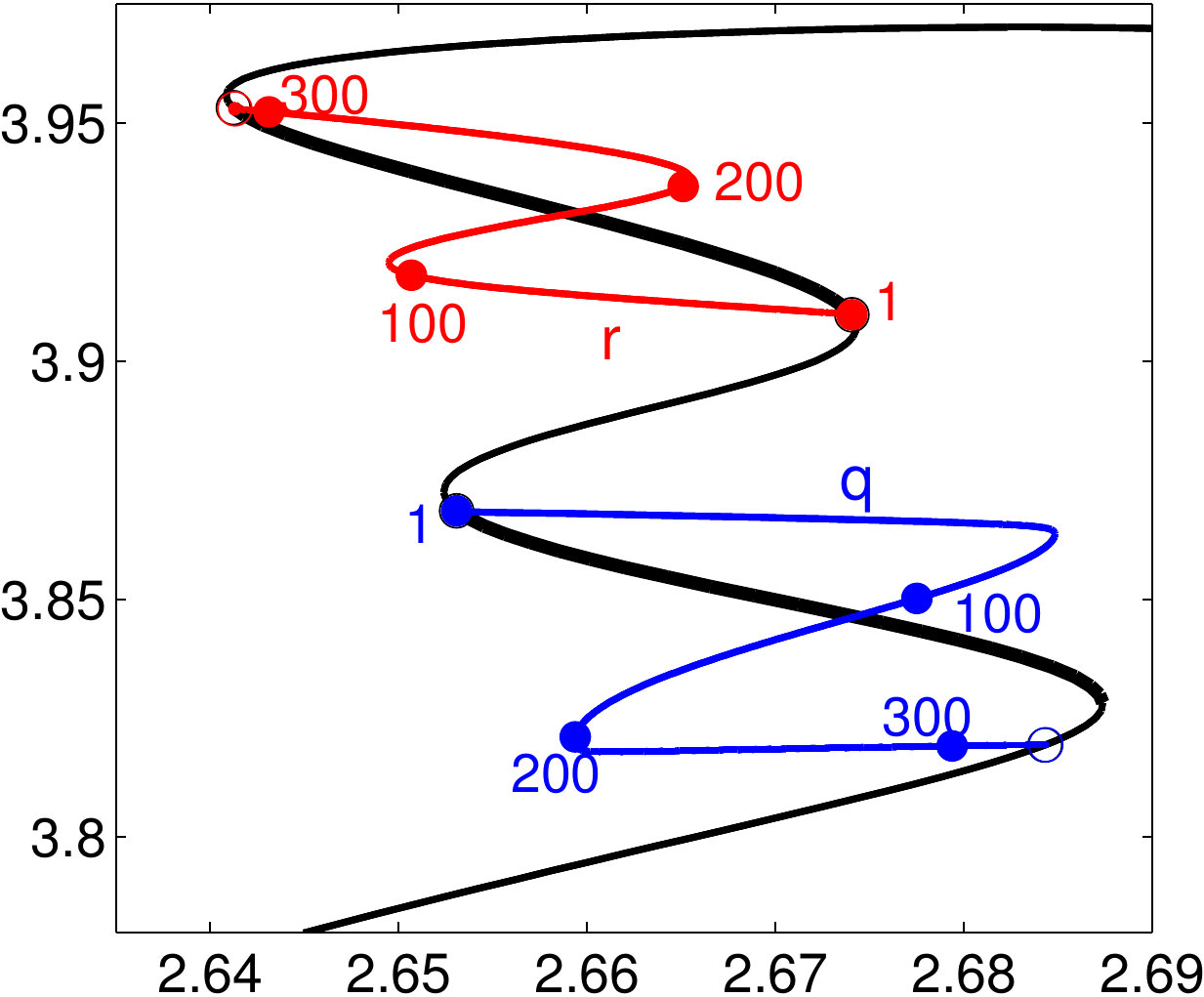}
\end{minipage}\hs{-1mm}
\begin{minipage}{50mm}(e) q1, q100, q200, q300\\
\ig[width=45mm]{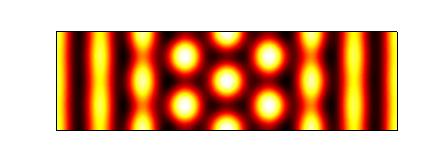}\\[-4mm]
\ig[width=45mm]{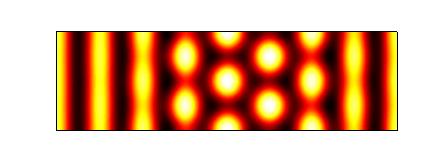}\\[-4mm]
\ig[width=45mm]{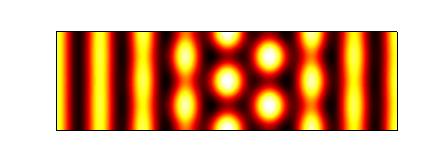}\\[-4mm]
\ig[width=45mm]{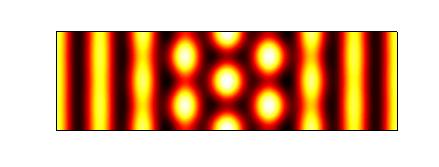}
\end{minipage}\hspace{2mm}
\begin{minipage}{50mm}(f)  r1, r100, r200, r300 \\
\ig[width=45mm]{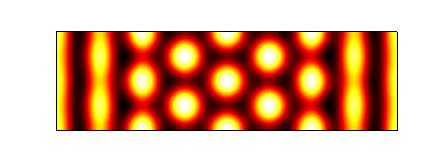}\\[-4mm]
\ig[width=45mm]{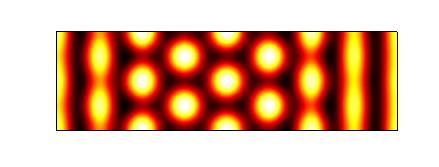}\\[-4mm]
\ig[width=45mm]{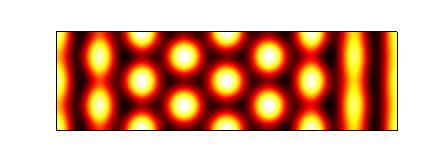}\\[-4mm]
\ig[width=45mm]{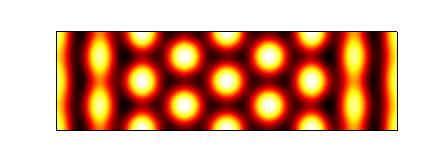}
\end{minipage}
}
\ece 

\vspace{-5mm}
\caption{{\small A magnification of the hl branch from Fig.~\ref{snakea},  
and the ``rung type'' branches connecting bifurcation points on the 
hl branch, and some show example plots.  The diagonal parts 
of the q and r branches can be seen as fragments of a  $\pi/k_cc$ shifted 
odd snake with wave number $k_c$, which does 
not exist over the given domain with homogeneous Neumann boundary conditions. 
See text for details, and \cite{smov} for a movie. 
} \label{snakeb}}
\end{figure}

In 1D models like \reff{qcshe} and \reff{cqshe}, in theory, and 
over sufficiently large domains with periodic boundary conditions, 
bifurcation points near the folds of a homoclinic snake are 
usually associated with ``rungs'' which together with a snake of 
the same wave number patterns but different parity (i.e., 
odd instead of even) form ladders, see in particular \cite{BKLS09}. 
In our case, due to the Neumann boundary conditions there is no snake 
of odd solutions 
with wave number $k_c$ (see however \S\ref{sidesec} for ``odd sideband 
snakes''). 
Thus, it is interesting 
to see the behaviour of branches bifurcating from the bifurcation points 
near the folds, which 
we illustrate in Fig.~\ref{snakeb}.  These points 
are again pairwise connected. On, e.g., the first horizontal 
part of the {\tt q} branch the solutions try to develop an odd 
symmetry around $x=\pi/k_c$ by growing/decreasing hexagons on the right/left. 
However, this cannot connect to the missing odd snake, and 
the branch turns around to a diagonal segment which 
would belong to the odd snake. Thus, this behaviour is again very much a 
finite size and boundary conditions effect, and similar comments apply to 
the {\tt r} branch. 

\brem\label{dawrem1} For both, the 1D cubic quintic SHe 
\reff{cqshe} and the quadratic cubic SHe \reff{qcshe}, 
an analogon of the bean branch is the 
branch $r_-$ of small amplitude unstable stripes that bifurcates 
subcritically from $(u,\lam)=(0,0)$ and ``connects'' to the fold 
where the stripes become stable, cf.~Fig.~\ref{ipic1}.  
In \cite{dawes1} the appearance of bifurcation points on $r_-$ 
is related to modulational instability of the amplitudes $A_-(\lam)$ 
of the unstable stripes in the associated Ginzburg--Landau equation, derived 
by using a scaling of the ``subcriticality parameter'' $\nu$ 
in \reff{cqshe}. 
In particular, this gives a lower bound on the domain size 
necessary for these bifurcations, which is linear in $1/\nu$, 
i.e., inversely proportional to the subcriticality. 
Moreover, it explains why the bifurcations occur 
in connected pairs like (in our case) ({\tt hbbp1},{\tt hbbp8}), 
({\tt hbbp2},{\tt hbbp7}), and so on. 
Similarly, using beyond all order  matched asymptotic expansions, 
a number interesting results 
for finite size effects on the snaking branches in the 1D quadratic cubic 
SHe are obtained in \cite{KAC09}, including that larger 
domains yield smaller pinning ranges, i.e., more narrow snakes. 

At least the relation between subcriticality of the 
equation and necessary domain size for secondary bifurcations 
also holds for our system, see Remark \ref{dawrem2}. We expect that 
calculations similar to those in \cite{dawes1} and in \cite{KAC09}  
can also be done in our case using the system of Ginzburg--Landau 
equations derived in \S\ref{AnaRes} (and extensions as 
in \cite{KAC09}), but naturally they will be 
more complicated.
\eex
\erem

In contrast to the homoclinic snaking of the {\tt hl} branch which at 
least for simpler models on infinite cylinders 
can be also explained and predicted theoretically, see \S\ref{AnaRes}, 
finite size effects should be regarded as essential 
for the ``heteroclinic snaking'' of the hot front branch {\tt hf} in 
Fig.~\ref{snakea}. In fact, at least for the 1D 
Swift--Hohenberg equation, beyond all order asymptotics \cite{chapk09} predict 
that for each $\lam$ near the Maxwell point there are at most 
two pinned front solutions, again see \S\ref{AnaRes}.  Thus, 
the hot fronts should be seen as homoclinics by even extension 
of the solutions over the left boundary (hexagons on a background of stripes) 
or right boundary (stripes on a background of hexagons). See also 
Fig.~\ref{nnf3} below. The {\tt m} (as in multiple) branch, 
consists of solutions with multiple fronts between stripes and hexagons. 

To further assess the finite size effects, in Fig.~\ref{nnf3} we plot 
homoclinic snakes from the second bifurcation points on the 
bean branches over $8\times 2$ and $12\times 2$ domains. 
\begin{figure}[!ht]
\bce
{\small\begin{minipage}{60mm}\ig[width=55mm]{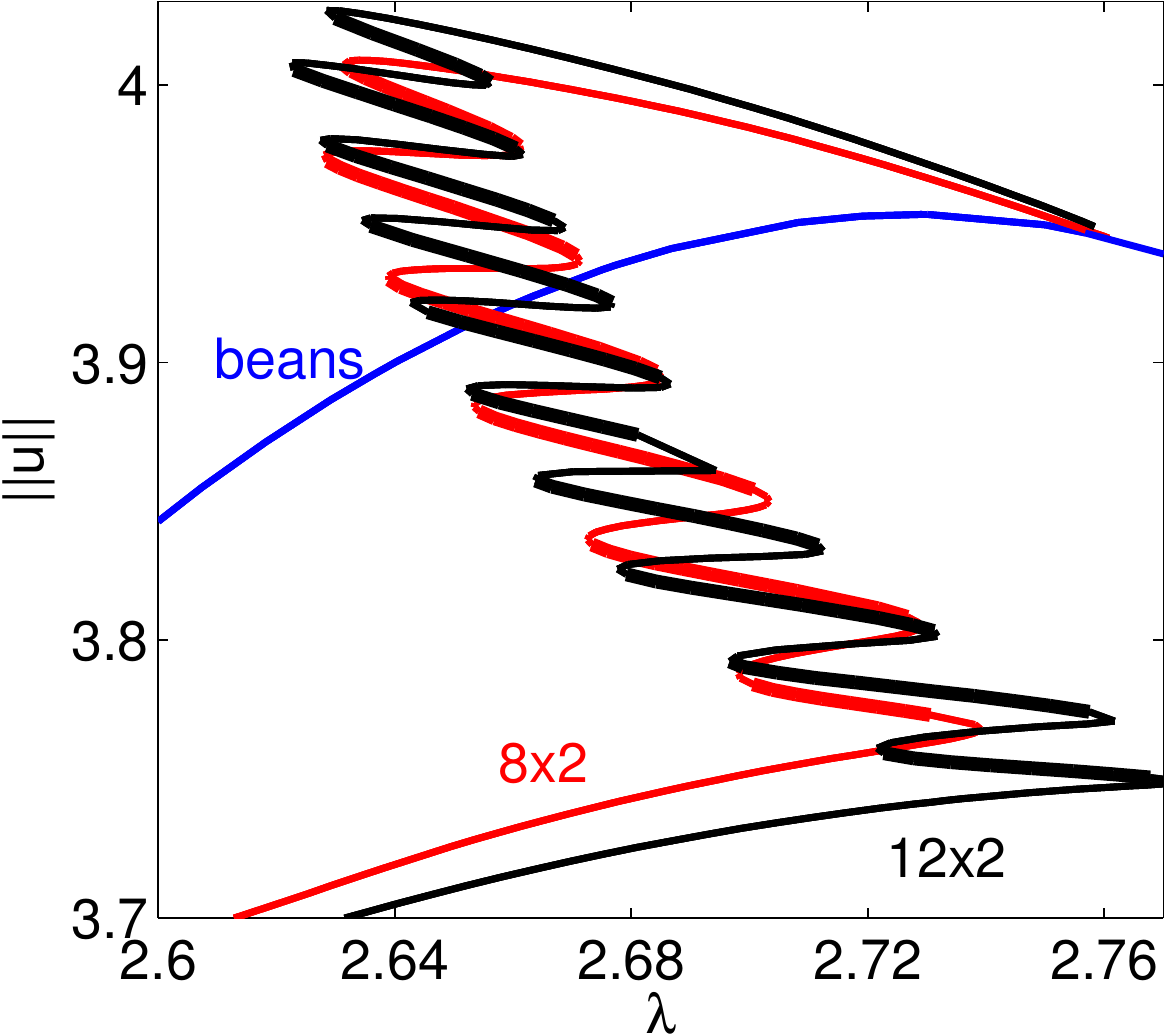}
\end{minipage}}\ece

\vspace{-5mm}
\caption{{\small Hot homoclinic snaking
 over $8\times 2$ and $12\times 2$ domains. \label{nnf3}}}
\end{figure}
The $8\times 2$ branch (see also Fig.~\ref{ipic1}) 
is very close to the {\tt hf} branch of 
the $4\times 2$ domain, thus justifying the claim that the $4\times 2$ {\tt hf} 
branch should be seen as a $8\times 2$ {\tt hl} branch by mirroring 
it over the left boundary. Remarkably, the 
$12\times 2$ {\tt hl} branch is {\em more} slanted than the $8\times 2$ branch. 
This agrees with the trend already seen in Fig.~\ref{snakea} by comparing 
the $4\times 2$ {\tt hf} and {\tt hl} branches, and can be further 
confirmed by considering even larger domains, but is contrary to numerical 
and analytical results for 1D snaking, see in particular \cite{KAC09} 
for the 1D quadratic cubic SHe \reff{qcshe}. There, as the domain size 
increases, snaking 
becomes {\em more} confined around the Maxwell point of the variational 
equation \reff{qcshe} and less slanted. At this point we cannot 
offer an explanation for Fig.~\ref{nnf3}, 
but can only point out that our system is 
not variational. Therefore, Maxwell points and energy arguments 
can only be given approximately 
by asymptotic expansions, see \S\ref{AnaRes}, and the $\lambda$ regime 
discussed so far is quite far from the primary bifurcation. Moreover, 
we also have to deal with a possibly 
rather complicated interaction between the length in $x$, which we increased, 
and the length in $y$, which we fixed, on the one hand 
to keep calculational costs at base, but also to first vary only one 
size parameter. In fact, 
below we shall find snakes which are much less slanted than those 
in Fig.~\ref{nnf3}. These are typically found much closer to $\lam_c$, and 
moreover involve different ratios between $x$ and $y$ wavenumbers, see 
for instance Fig.~\ref{nnf4}.

\subsection{Time integration}\label{tintsec}
The stability information in Figures \ref{figsch}, \ref{snakea}, \ref{snakeb} 
and \ref{nnf3} 
is based on the spectrum of the linearization around a (FEM) solution. Given 
the rather large and complicated space of solutions it is interesting 
to assess the nonlinear stability and basin of attraction of solutions. 
For this we use time integration of the spatial FEM discretization of 
\reff{rd} by some standard semi implicit method. This is also quite 
useful to quickly obtain nontrivial 
starting points for continuation and bifurcation. Generally speaking, 
we obtain desired stable target solutions from only rough initial guesses, 
where typically we integrate for some rather short time and then use 
a Newton loop for the stationary problem to get to the stationary solution. 
\begin{figure}[!ht]
\bce
{\small
\begin{minipage}{50mm}
\ig[width=50mm, height=45mm]{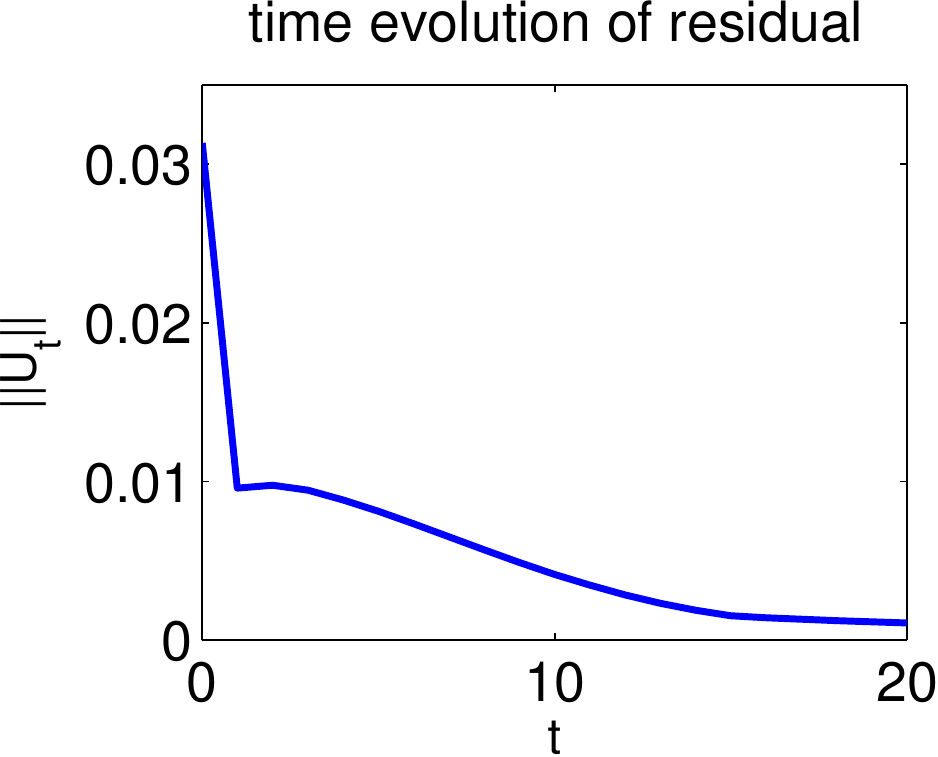}
\end{minipage}\hs{0mm}
\begin{minipage}{90mm}
\ig[width=80mm]{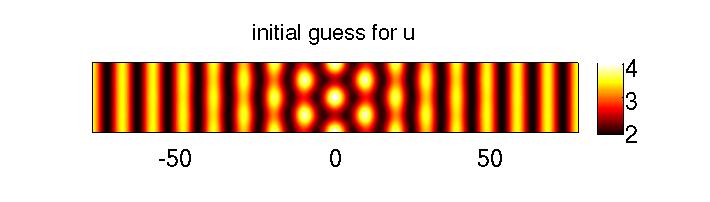}\\[-4mm]
\ig[width=80mm]{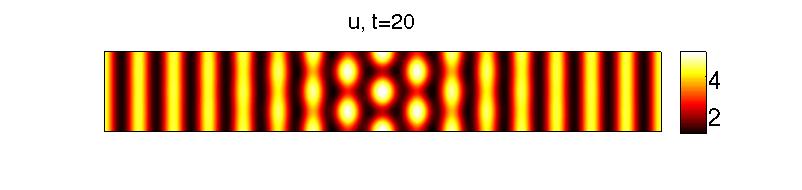}\\[-4mm]
\ig[width=80mm]{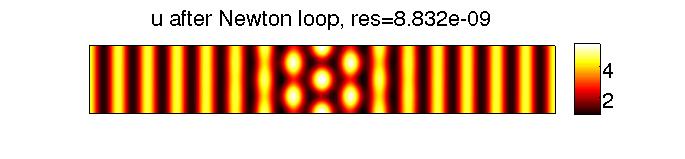}
\end{minipage}
}\ece

\vspace{-5mm}
\caption{{\small Finding a pattern from an initial guess and 
time integration, on an $8{\times}2$ domain. Here $\lam{=}2.7$, 
$U_0$ from \reff{u0g}. For this $U_0$, a Newton loop for the 
stationary problem does not converge, but it does after decreasing the 
residual $\|\pa_t{\tt U}\|_\infty$ by a number of time steps, 
where $\tt U$ denotes the spatial discretization of $U$. 
\label{nnf5}}}
\end{figure}
See Fig.~\ref{nnf5} for an example, where, with $A{=}0.3$, $B{=}0.15$ and 
$L{=}12$, 
\huga{\label{u0g}
U_0=\bpm \lam\left(1+A\cos(k_cx) 
+{\rm sech}(\frac x L)[2B\cos(\frac{k_c} 2x)\cos(k_c\frac{\sqrt{3}}2y)
-0.1\cos(k_cx)]
\right)\\
\frac 1 \lam\left(1-\frac 1 2\left(A\cos(k_cx) 
+{\rm sech}(\frac x L)[2B\cos(\frac{k_c} 2x)\cos(k_c\frac{\sqrt{3}}2y)
-0.1\cos(k_cx)]\right)\right)
\epm.
}
Allthough motivated by the desired mode structure, the precise form of 
\reff{u0g} is rather arbitrary. We obtain convergence to ``hot localized'' 
solution, where of course the width of the hexagon patch over the stripe 
background, i.e., the position in the snake, depends on the initial 
width $L$. 

\begin{figure}[!ht]
\bce
{\small
\begin{tabular}{p{55mm}p{50mm}p{50mm}}
(a) $\lam{=}2.65$, conv.~to localized hexagons further up the snake&
(b) $\lam=2.6$, conv.~to space filling hexagons&
(c)  Norm evolution\\[-0mm]
\ig[width=50mm]{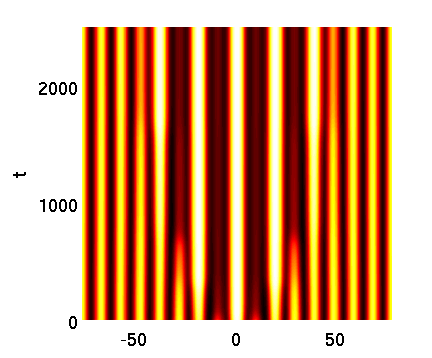}&
\ig[width=50mm]{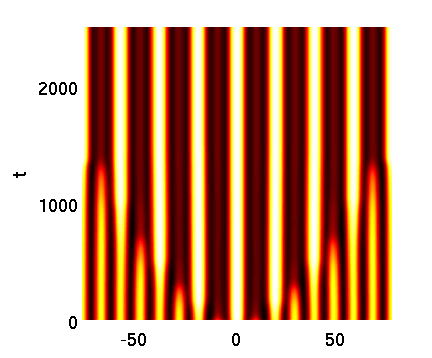}&
\ig[width=40mm,height=40mm]{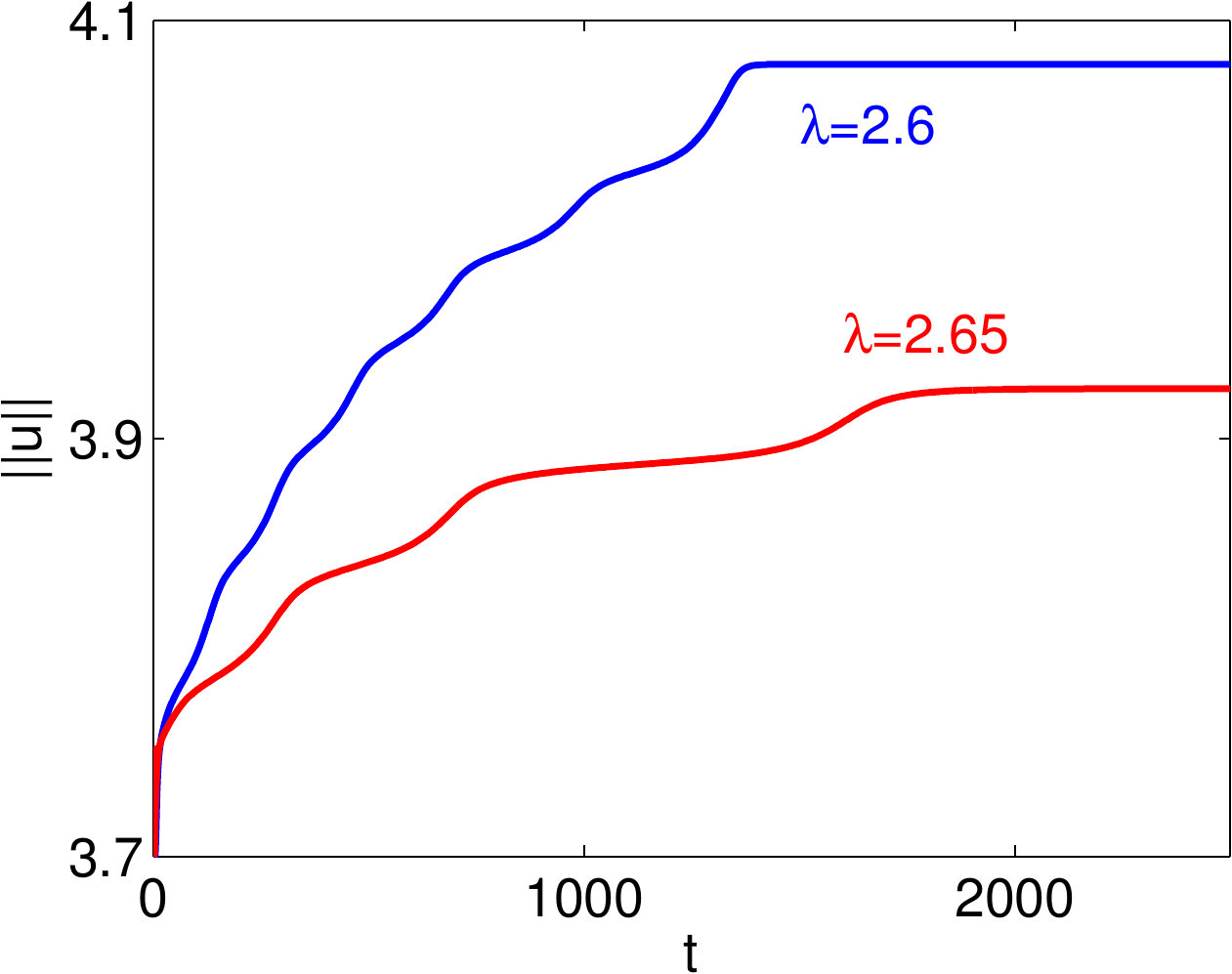}
\end{tabular}
}\ece

\vspace{-5mm}
\caption{{\small Stick-slip motion after ``stepping out of the snake''; 
(a), (b) $u(t,x,-l_y)$ in space--time diagram, initial 
conditions from Fig.~\ref{nnf5}. 
\label{ssfig}}}
\end{figure}

Another interesting numerical experiment is to choose an initial point 
on a snaking branch and externally modifying $\lam$ to some other value. 
For instance, in Fig.~\ref{ssfig} we use the stationary solution $U$ for 
$\lam=2.7$ from Fig.~\ref{nnf5} and modify $\lam$ to $\lam=2.65$, which 
is still inside the snaking region, and to $\lam=2.6$, which is just outside. 
To illustrate the time evolution of $u(t,x,y)$ we plot a space time diagram 
of $u(\cdot,\cdot,-l_y)$, i.e., of $u$ on the lower boundary 
of $\Om$. In these the hexagons appear 
as ``rolls'' with double  the wavelength of the stripes due to 
$|k_x|=\frac 1 2 k_c$ for hexagons. For both,   $\lam=2.65$ and $\lam=2.6$, 
we observe two fronts moving outwards with a ``stick--slip'' motion, 
where $n(t):=\|u(t,\cdot,\cdot)\|$ 
becomes flat as $u(t,\cdot)$ comes close to the (here left) folds in the 
$8\times 2$ 
snaking branch from Fig.~\ref{nnf3}. For $\lam=2.6$ we then end up 
with space filling hexagons, and for $\lam=2.65$ with some localized hexagons
``further up'' on the snake, which again shows stability of these solutions. 

\subsection{Planar fronts, localized patterns and pinning in the cold range}
\label{crsec} 
\begin{figure}[!ht]
\bce
{\small
\begin{minipage}{59mm}(a)\\
\ig[width=60mm, height=55mm]{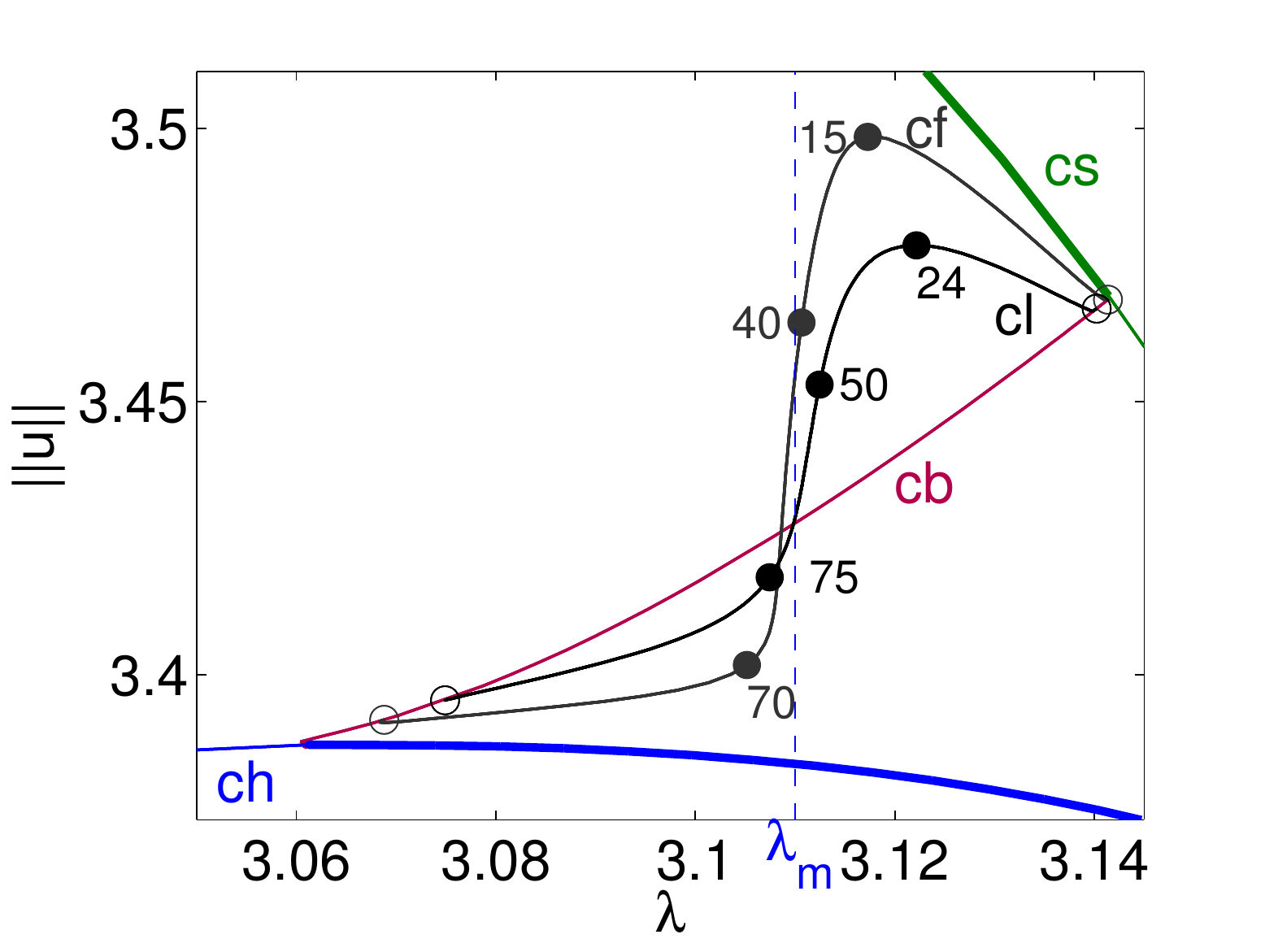}
\end{minipage}\hs{-5mm}
\begin{minipage}{95mm}(b) {\tt cf15, cf40} and {\tt cl24, cl50}\\
\ig[width=105mm]{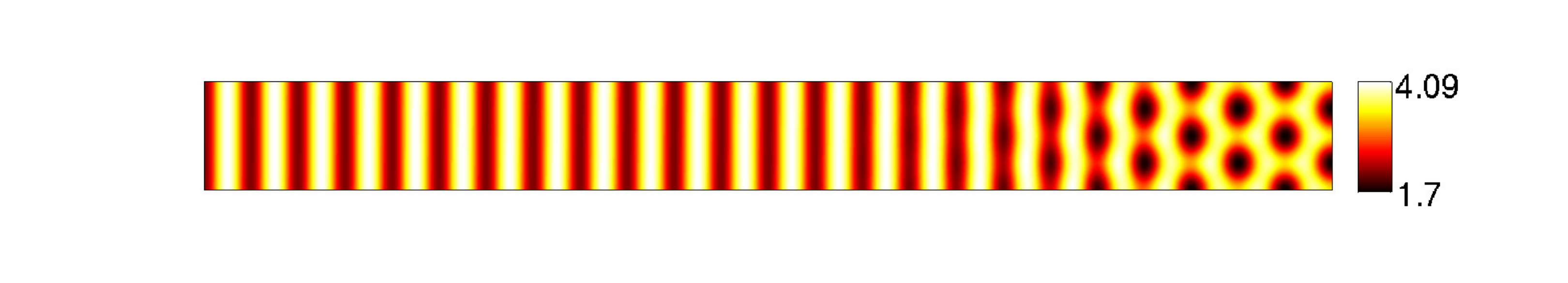}\\[-4mm]
\ig[width=105mm]{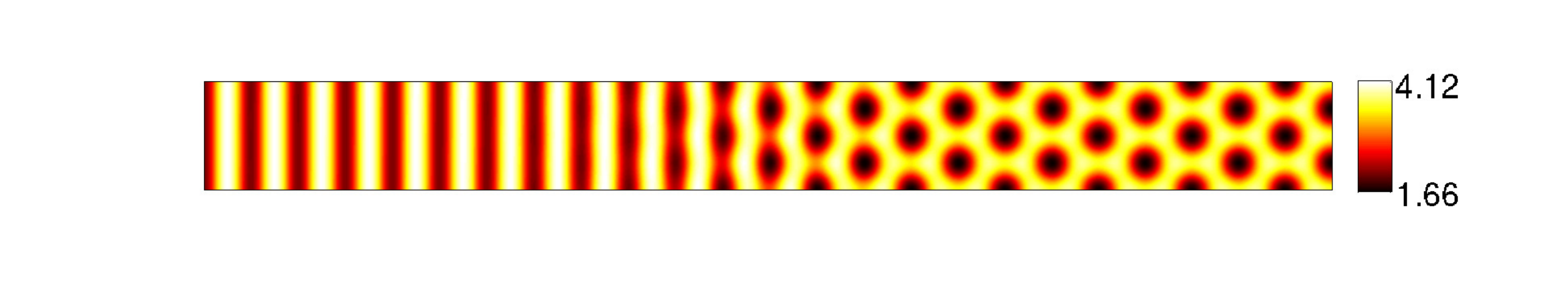}\\[-4mm]
\ig[width=105mm]{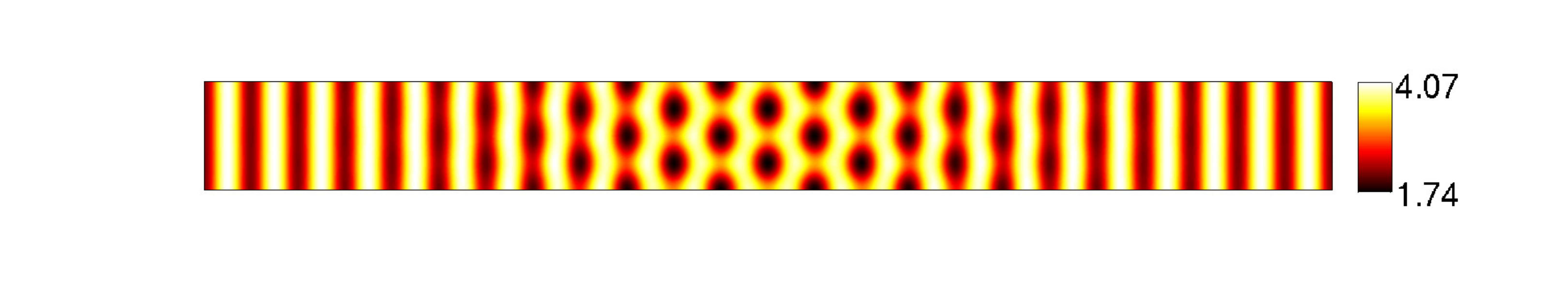}\\[-4mm]
\ig[width=105mm]{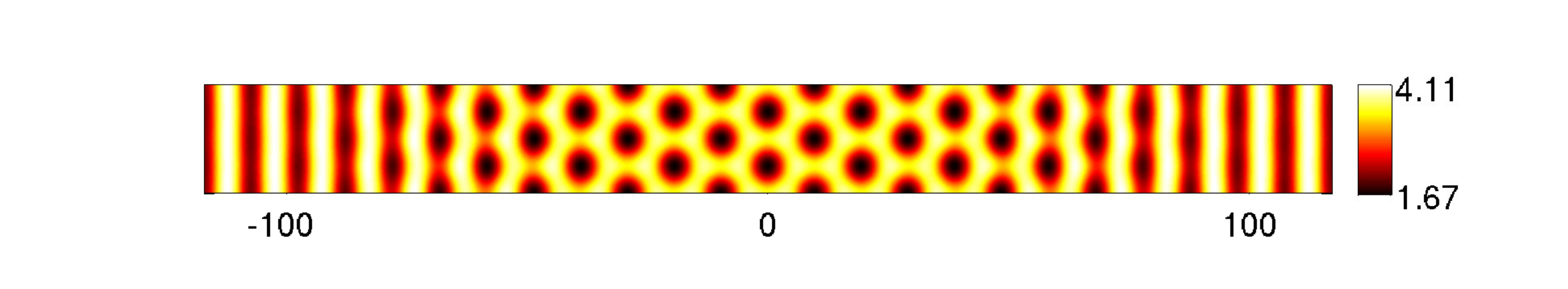}
\end{minipage}
}\ece

\vspace{-5mm}
\caption{{\small Bifurcation from the cold beans on a $12\times 2$ domain.  
(a) {\tt cf} branch connecting {\tt cbbp1} with {\tt cbbp8} in a monotonous 
way (no snaking), cold localized branch {\tt cl}, 
and the cold Maxwell point (see \S\ref{AnaRes}). (b) Some example solutions. 
See \cite{smov} for a movie.}
\label{huf3}}
\end{figure}

On the cold bean branch {\tt cb} there are no bifurcation points over 
a $2\times 2$ or $4\times 2$ domain.  
However, Fig.\ref{huf3}(a) shows 4 of the 8 bifurcation points obtained 
on the cold bean branch over a $12\times 2$ domain, 
namely {\tt cbbp1, cbbp2, cbbp7, cbbp8}, the branch {\tt cf} (cold front) 
connecting {\tt cbbp1} and {\tt cbbp8}, and the branch {\tt cl} connecting 
{\tt cbbp2} and {\tt cbbp7}, and consisting of hexagons on 
a stripe background.  Panel (b) shows 4 example 
solutions. In, e.g., {\tt cl50} it can be seen that the localized 
hexagons are shifted by $2\pi/k_c$ in $x$, such that here we really 
consider $\CS{\tt cb}$ which connects {\tt cs} to $\CS{\tt ch}$, 
cf.~Remark \ref{srem}, but for simplicity we omit the $\CS$. 

Similar to Fig.~\ref{snakea}, after bifurcation from, e.g., {\tt cbbp1}, 
in the traverse to {\tt cf15} the solutions reshape 
into stripes on the left and hexagons on the right. The growth  
of the hexagon part then happens in a narrow $\lam$  regime around $\lam_m$ 
between {\tt cf15} and {\tt cf70}, while between {\tt cf70} and 
{\tt cbbp8} the solutions are reshaped to near hexagon beans.  
Similar remarks hold for the {\tt cl} branch. 

\begin{figure}[!ht]
\bce
{\small
\begin{tabular}[t]{lp{70mm}lp{70mm}}
(a)&(b)\\
\ig[width=60mm,height=55mm]{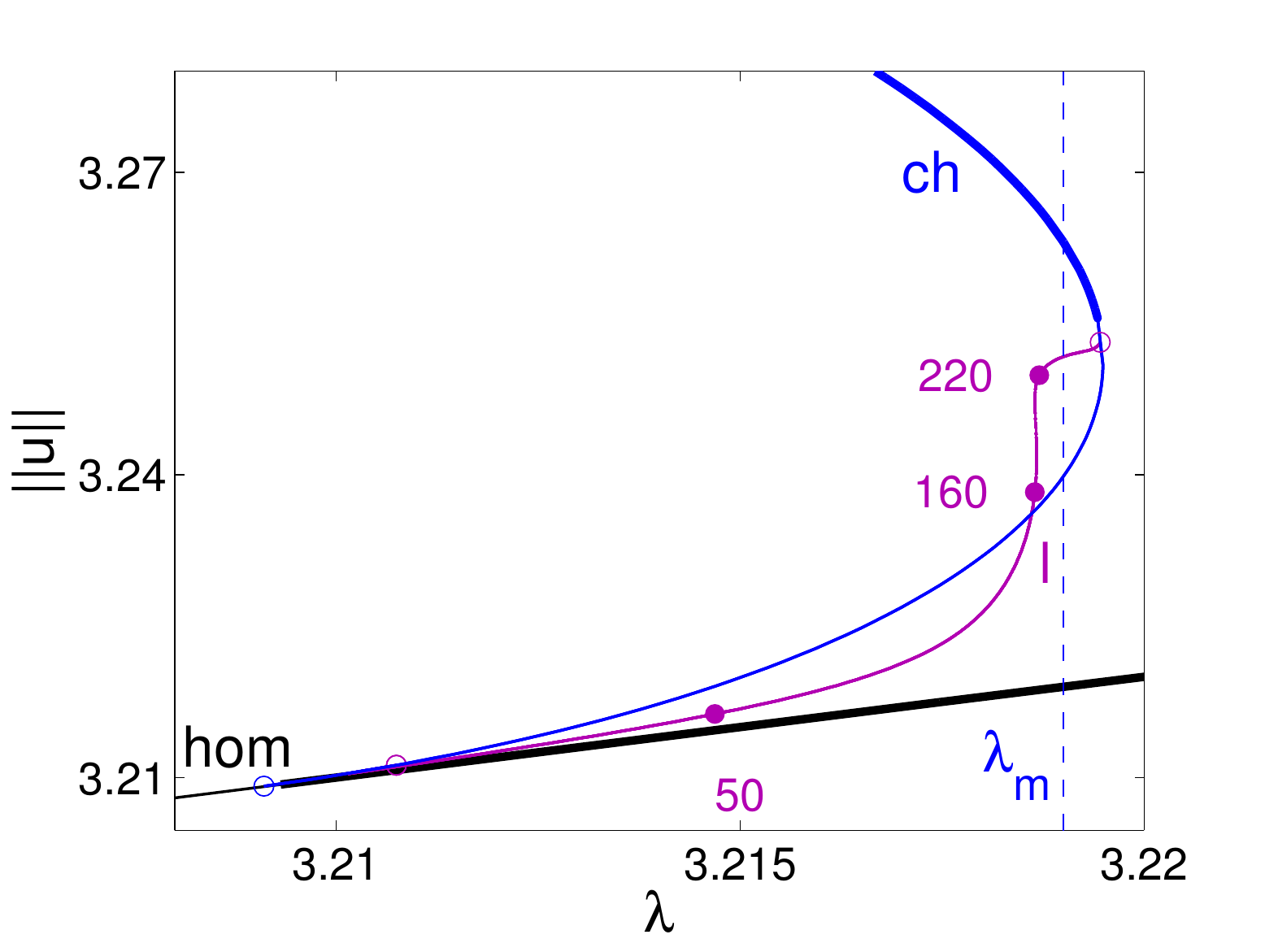}&
\ig[width=60mm,height=50mm]{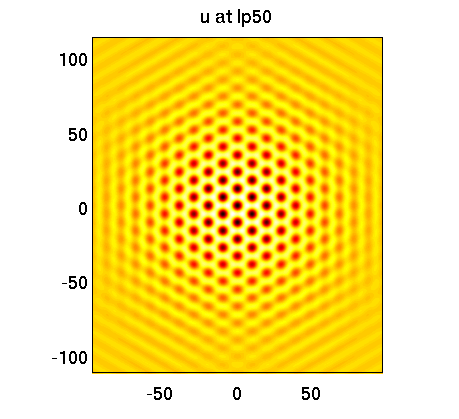}\\
(c)&(d)\\
\ig[width=60mm,height=50mm]{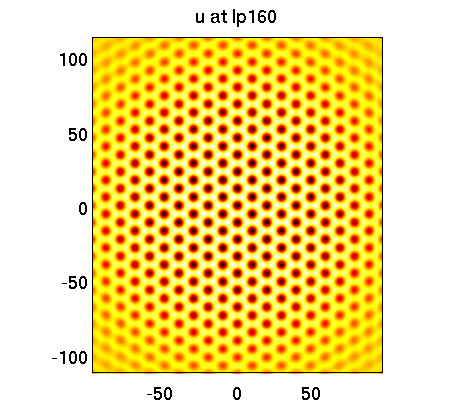}&
\ig[width=70mm]{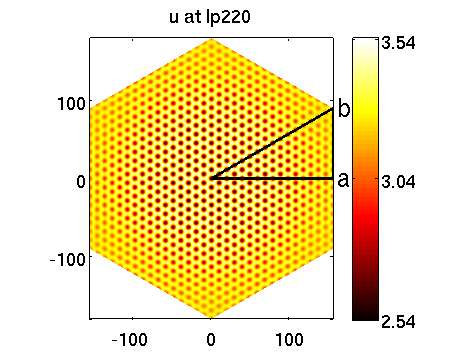}\\
(e)
\end{tabular}\\
\ig[width=150mm]{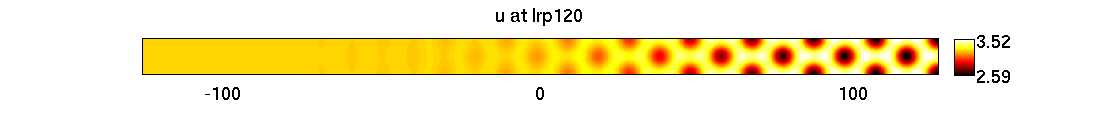}
}\ece

\vspace{-5mm}
\caption{{\small Bistable range between homogeneous solution and cold 
hexagons. 
(a) partial bifurcation diagram and the homogeneous Maxwell point (see \S\ref{AnaRes}). 
Additional to quasi 1D solutions we find 
fully 2D localized hexagon patches. The triangle in (d) indicates 
the computational domain with $a=(32\pi/k_c,0)$ and 
$b=(32\pi/k_c,32\pi/(\sqrt{3}k_c)$, and about 90.000 triangles. 
We use Neumann BC on all three sides and for plotting first make an 
even extension over $y=0$ and then five rotations by 60 degrees. 
(e) Quasi 1D front between $w^*=(\lam,1/\lam)$ and spots 
on a $13\times 1$ domain, $\lam\approx 3.219$. See \cite{smov} for a movie.}
\label{huf4}}
\end{figure}

The third bistable range we discuss is between the homogeneous solution 
$w^*=(\lam,1/\lam)$ and the cold hexagons 
for $3.21\approx\lam_c<\lam<\lam_{ch}^b\approx 3.22$. 
Here, over sufficiently large domains, we have bifurcation points 
on the small amplitude (unstable) cold hexagon branch, which now 
corresponds to the $r_-$ branch for \reff{cqshe} in Fig.~\ref{ipic1}. 
As above we can find quasi 1D fronts between $w^*$ and cold 
hexagons, see Fig.~\ref{huf4}(e), and the associated 1D localized patterns. 
However, as in the qcSHe \reff{qcshe} and the cqSHe \reff{cqshe}, 
see \cite{strsnake} and the references therein, 
we can also calculate fully localized patches of hexagons, 
see Fig.~\ref{huf4}(a)--(d). 

\brem\label{dawrem2}
An essential difference between Fig.~\ref{snakea} and Figures \ref{huf3} 
and \ref{huf4} is that in the latter there is no snaking. Starting 
with Fig.~\ref{huf3} and following Remark \ref{dawrem1} we believe 
that the need for a large 
domain to obtain bifurcations from the {\tt cb} branch can be explained 
from the ``weak subcriticality'' of the system in this range, such 
that the cold bistable range 
is rather narrow, compared to the hot bistable range. 
Related to this, the difference in 
(Ginzburg--Landau--)energy between {\tt ch} solutions and 
{\tt cs} solutions is much 
smaller than between {\tt hh} solutions and {\tt hs} solutions in their 
bistable range and this results in flatter fronts 
for the associated Ginzburg--Landau system, which we discuss in more 
detail in \S\ref{glsec}. 

However, here increasing the domain gives bifurcation points on {\tt cb} 
(e.g., 4 on a $12\times 2$ domain), 
and branches connecting the first and the last, the second and 
the second to last, and so on, as in \cite{dawes1}, 
but no snaking. The reason is most likely an exponentially small width 
of the snaking region, which cannot be resolved by our numerics, analogous 
to a number of results for the qcSHe \reff{qcshe} and the cqSHe \reff{cqshe}, 
which can directly be related to Fig.~\ref{huf4}. 
For instance, \cite[Fig.33,34]{hexsnake}, see also \cite[Fig.22]{strsnake}, 
show the 
(numerically calculated) positions of the left $\lam_-(\nu)$ 
and right folds $\lam_+(\nu)$ in the snake of hexagons over $u=0$ in 
\reff{qcshe}. Then, as expected, $\lam_-(\nu)<\lam_M(\nu)<\lam_+(\nu)$, 
where $\lam_M(\nu)$ is the relevant Maxwell point, but the width 
$s(\nu)=\lam_+(\nu)-\lam_-(\nu)$ goes to $0$ quickly as $\nu$ becomes small, 
and the snaking cannot be continued down to $\nu=0$. 
Analytically, \cite{chapk09, suma11, dean11, KC13} explain that the 
pertinent small parameter in \reff{qcshe} and \reff{cqshe} 
is the subcriticality  
$\eps(\nu)=\lam_c-\lam_f(\nu)$, where $\lam_c=0$ and $\lam_f(\nu)<0$ 
denotes the position of the fold where the subcritical rolls become 
stable. Formulas for 
$s(\eps)=\lam_+(\eps)-\lam_-(\eps)$ 
can then derived using beyond all order asymptotics, 
and show the exponential smallness of $s(\eps)$ in $\eps$, i.e., 
\huga{\label{swform}
s(\eps)\sim C\eps^{-\al}\exp(-\beta/\eps)\text{ as }\eps\ra 0, 
}
where $C,\al,\beta>0$ are numerical constants 
depending on the parameters of the problem and, e.g., front orientation, 
but not on $\eps$. These formulas 
show good agreement with numerical simulations, 
see, e.g., \cite[Fig.2]{KC06},\cite[Fig.5]{dean11}, \cite[Fig.1(c)]{suma11}, \cite[Fig.2,3]{KC13}. 

If we also assume a dependence like \reff{swform} here, with 
$\eps=\lam_{s}^b-\lam_{ch}^e\approx 0.12$ the ``bistability width'' of 
stripes and cold hexagons in Fig.~\ref{huf3}, and 
$\eps=\lam_{ch}^b-\lam_c\approx 0.02$ the subcriticality in Fig.~\ref{huf4}, 
then the snaking width might be too small 
to resolve numerically. In fact, increasing the domain size and 
keeping the discretization reasonably fine, the steep 
part near $\lam_m$ becomes steeper, 
but in our numerics this process converges 
to a monotonous in $\lam$ branch. To elaborate this we set 
up the modification \reff{mod1} which 
for  $\sigma<0$ increases the subcriticality of the cold hexagons. 
See \S\ref{sigsec} where we find snaking of the analoga of the {\tt cl}, 
and {\tt l} branches of Fig.~\ref{huf3} and \ref{huf4}, respectively. 
\eex
\erem

\subsection{Sideband patterns}
\label{sidesec}
By our choice of  $l_1\times l_2$ domains 
$\ds\Om=\left(-\frac{2l_1\pi}{k_c},\frac{2l_1\pi}{k_c}\right)
\times \left(-\frac{2l_2\pi}{\sqrt{3}k_c},\frac{2l_2\pi}{\sqrt{3}k_c}\right)$ 
with homogeneous Neumann boundary conditions, 
in a neighborhood of $(w,\lam)=(0,\lam_c)$ {\em only} patterns with 
basic wave vectors $k_c(1,0)$ and 
$k_c(-\frac 1 2,\pm\frac{\sqrt{3}} 2)$ exist, as these are 
the only unstable modes that fit into the domain. 
However, for larger $\lam_c-\lam$ more modes become unstable 
and hence patterns with sideband wave vectors like 
$k=k_c(\frac{4l_1+m_1}{4l_1},0)$, $m_1=\pm1,\pm2,\ldots$ (sideband stripes), 
or $k=k_c(-\frac 1 2\frac{4l_1+m_1}{4l_1}, 
\pm\frac {\sqrt{3}} 2\frac{4l_2\pm m_2}{4l_2})$, 
$m_{1,2}=\pm1,\pm2,\ldots$ (with 
$\|(\frac 1 2 \frac{4l_1+m_1}{4l_1}, \frac {\sqrt{3}} 
2\frac{4l_2\pm m_2}{4l_2})\|$ close to unity to build sideband hexagons)
come into play: branches of solutions with 
such modes bifurcate from the homogeneous branch 
below $\lam_c$ and should be expected to behave roughly similar as 
the basic stripe and hexagon branches, i.e., for $\lam<\lam_c$ 
we get an existence balloon of patterns. Moreover, allthough unstable 
at bifurcation they may become stable away from bifurcation, forming 
the so called Busse balloon as a subset of the existence balloon. 
Of course, all this heavily depends on the size of the domain: for instance, 
on the $2\times l_2$ domain, the sideband (vertical) stripes have 
wave vectors $(k_1,0)$ 
with $k_1/k_c=\frac 1 4,\frac 1 2,\frac 3 4,\frac 5 4,\ldots$, while on 
the $8\times l_2$ domain we have 
$k_1/k_c=\ldots,\frac {15}{16},\frac {31}{32},\frac {33}{32},\frac {17}{16}, 
\ldots$, where for instance we call the stripe with expansion 
$A\sin(\frac{31}{32} k_cx)+\hot$ an odd stripe, in contrast to the 
even stripes considered so far. 

\begin{figure}[!ht]
\bce
{\small\begin{minipage}{80mm}(a)\\
\ig[width=65mm, height=60mm]{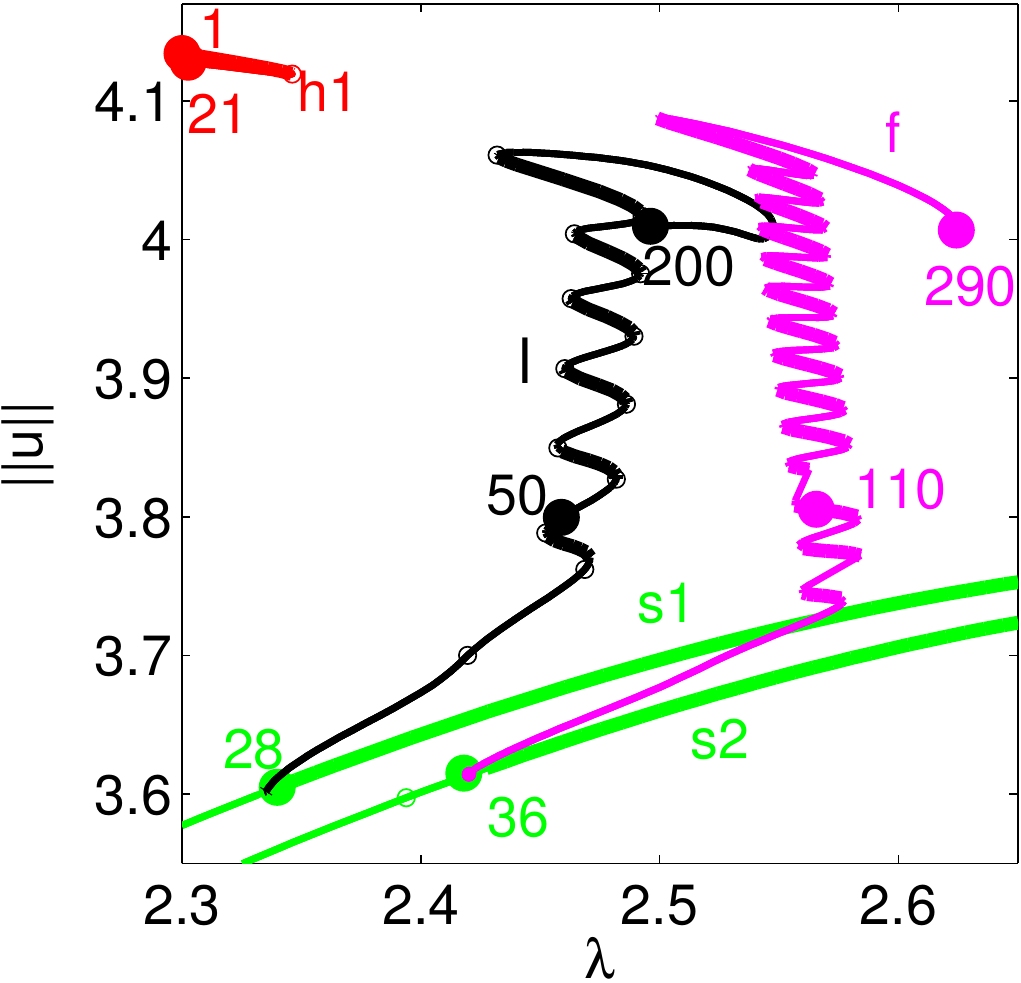}
\end{minipage}\hs{-5mm}
\begin{minipage}{80mm}(b) {\tt s1p28, s2p36}, and tangents at bif.~for 
{\tt l} and {\tt f}\\[-0mm]
\ig[width=85mm]{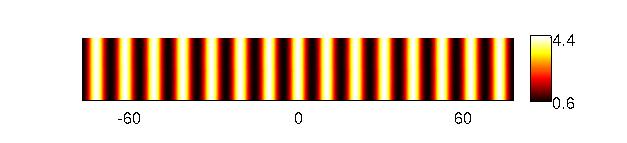}\\[-2mm]\ig[width=85mm]{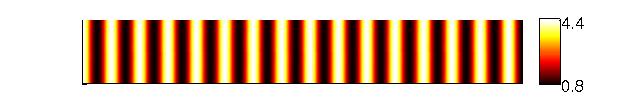}\\
\text{}\hs{2mm}\ig[width=85mm]{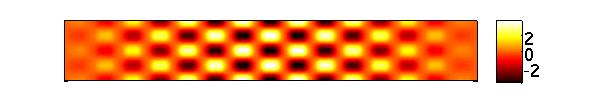}\\[-4mm]
\text{}\hs{-0.5mm}\ig[width=88.5mm]{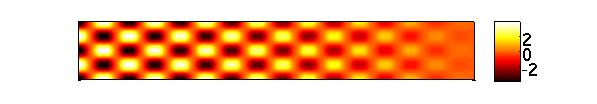}\\[-2mm]
\end{minipage}\\
\begin{minipage}{80mm}(c) l50 and l200\\[-0mm]
\ig[width=81mm]{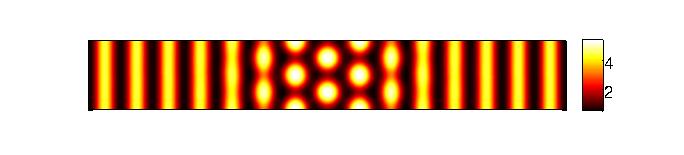}\\[-4mm]\ig[width=81mm]{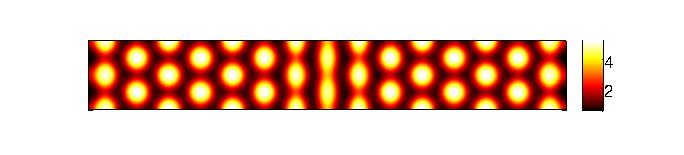}
\end{minipage}
\begin{minipage}{80mm}(d) f110 and f290\\
\text{}\hs{-1mm}\ig[width=82mm]{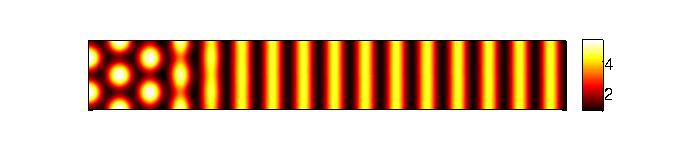}\\[-4mm]
\text{}\hs{-1mm}\ig[width=82mm]{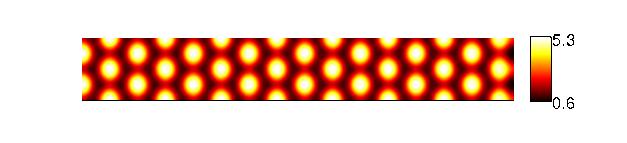}
\end{minipage}\\
\begin{minipage}{160mm}(e) {\tt h1p1 and h1p21}\\
\text{}\hs{-1mm}\ig[width=82mm]{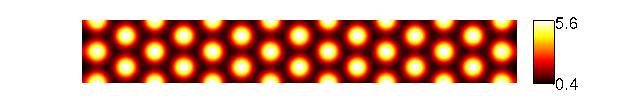}
\hs{-2mm}\ig[width=82mm]{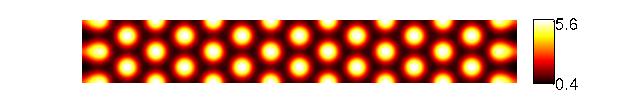}
\end{minipage}
}\ece

\vspace{-5mm}
\caption{{\small (a) Bifurcation diagram for patterns with sideband 
wavenumbers on the $8\times 2$ domain. 
(b) Stripes with $k_1=\frac {15}{16}k_c$ ({\tt s1}) 
and $k_1= \frac {31}{32}k_c$ ({\tt s2}), and 
bifurcation directions from these. (c),(d) Some example plots on the bifurcating 
branches. (e) Slightly distorted hexagons with  
$k_1=\frac {15}{16}k_c$, obtained from an initial guess and 
time integration followed by a Newton loop at $\lam=2.3$. 
See \cite{smov} for a movie.
\label{nnf4}}}
\end{figure}

Just as an indication of how such sideband patterns can further complicate 
the bifurcation diagram of \reff{rd}, in Fig.~\ref{nnf4} we plot a number 
of ``sideband'' branches and sample solutions on the $8\times 2$ domain, 
thus allowing many more wave numbers $k_1=0,\frac 1 {32} k_c, \frac {2}{32} k_c, 
\ldots$ in $x$ direction than wave numbers $k_2=0,\frac 1 {4} k_c, \ldots$ 
in $y$ direction. The branches {\tt s1, s2} belong to (hot) ``subharmonic'' 
stripes with $k_1=\frac {30} {32} k_c, \frac {31}{32} k_c$, respectively. 
Interestingly they stay stable much longer than the basic hot stripe 
branch {\tt hs} 
which looses stability at $\lam=\lam_{s}^e\approx 2.51$. However, this is not a 
contradiction  
as  $k_c$ need not be in the 
Busse balloon of stable wavenumbers for $|\lam-\lam_c|=\CO(1)$, see, e.g., 
\cite{DRS12} for an interesting example.

In contrast to the {\tt hs} branch, 
when loosing stability, the bifurcation is no longer to beans, but 
directly to a localized branch {\tt l} in case of {\tt s1} and 
to a front branch {\tt f} in case of {\tt s2}. For {\tt s1}, an explanation is  
that hot hexagons with $k_1$ wave numbers $\frac {15} {32}k_c$ 
do not exist for $\lam>\lam_h\approx 2.35$: for $\lam=2.3$ we can generate 
such hexagons from a suitable initial guess and time--integration, 
see \S\ref{tintsec}, but as we try to continue the {\tt h1} branch to 
larger $\lam$ we get a fold at $\lam_h$, and on the lower branch the hexagons 
loose their shape as $\lam$ increases again. This is a finite 
size effect, in particular of the rather small size in $y$. 
The wave vector $k_{h1}=k_c(-\frac {15}{32},\pm\frac {\sqrt{3}} 2)$ of 
the {\tt h1} branch which yields the slightly distorted hexagons 
shape is the one with $|k|$ closest to $k_c$ in the family 
$k=k_c(\frac {15}{32}, \frac {\sqrt{3}} 2\frac{8\pm m_2}{8})$,  
and apparently this is not yet in the existence balloon for 
$\lam>\lam_h$. Consequently, we should not expect beans involving 
$k_{h1}$ for $\lam>\lam_h$, and this explains the direct bifurcation from 
the {\tt s1} branch to the {\tt l} branch. 

Similar remarks apply to the bifurcation from the {\tt s2} branch to 
the {\tt f} branch, where moreover no hexagons with 
$k=k_c(-\frac{31}{64},\pm\sqrt{3}/2)$ exist 
on the $8\times 2$ domain at all, for any $\lambda$ in the depicted 
range. Thus, bifurcation 
to the {\tt f} branch is the natural candidate, and by mirroring 
solutions over the left or right boundary this corresponds to an {\tt l} branch 
over the $16\times 2$ domain. On both snakes the segments pointing north--west 
are stable, and there are bifurcation points near the folds. 
Also note that point 200 
on the {\tt l} branch has stripes on a hexagon background and thus 
is already on the ``way back'' to the bifurcation 
point of {\tt l} from {\tt s1}. Here we close our brief discussion of 
2D sideband patterns, but refer to \S\ref{1dsec} for (1D) branches 
connecting different wave numbers, and to \S\ref{flsec} for fully localized 
sideband hexagons over stripes. 

\section{Ginzburg-Landau reduction}\label{AnaRes}
We now approximate the hexagons, stripes and mixed modes by the 
Landau formalism, and the fronts and 
localized patterns with planar interfaces by the Ginzburg--Landau formalism. 
In particular, using an energy argument 
and the so called Maxwell point for the Ginzburg--Landau system we 
find an approximate prediction where to find the 
respective branches in the bifurcation diagram for \eqref{rd}. 
Similar ideas are worked out much more deeply for the 1D 
quadratic--cubic SHe \reff{qcshe} in \cite{chapk09} (see also 
\cite{KAC09}) 
and the cubic--quintic SHe \reff{cqshe} in \cite{dean11}, where, 
by augmenting the Ginzburg--Landau ansatz with beyond all order 
asymptotics, accurate bifurcation diagrams for homoclinic snaking 
were rigorously derived. 

In our case we need a system of Ginzburg--Landau equations, and 
our analysis is more formal 
since even a consistent derivation of the Ginzburg--Landau system 
is difficult as we refrain from 
scaling assumptions for quadratic interactions but rather work with the 
numerical coefficients in \reff{rd}. However, from the Ginzburg--Landau 
system we calculate the Maxwell point as a necessary condition for 
Ginzburg--Landau fronts, and the pinning argument from \cite{pomeau} then 
suggest the existence of stationary fronts for \reff{rd}. 
This approximation turns out to be qualitatively and at least in 
the cold regime also quantitatively correct, and thus it gives a 
lowest order approximation for the numerical solutions, although 
it cannot explain the snaking (or the non--snaking). 

In a second step we use the Ginzburg--Landau formalism as a predictive 
tool to choose the $\sigma$ direction in the modified model \reff{mod1} 
in which we can increase the subcriticality of the cold regime 
to obtain cold snaking branches (\S\ref{sigsec}), and in which we have 
subcritical bifurcation of stripes and hence can 
expect 1D snaking (\S\ref{1dsec}).

\subsection{Landau description of hexagons, stripes, and 
mixed modes}
To formally describe planar fronts for \eqref{rd} by a 
Ginzburg-Landau system the idea 
is to treat $x$ as an unbounded variable, while 
$y\in [-\frac{m\pi}{\sqrt{3}k_c}, \frac{m\pi}{\sqrt{3}k_c}]$ with 
Neumann boundary conditions, as in the 
numerics. At least close to $\lam_c$ the most unstable modes 
of the linearization $L(\Delta)$ of \reff{dglredu0} around $0$ 
are then $e_1\Phi$, $e_2\Phi$, and $e_3\Phi$, where 
$e_1=\er^{\ri k_cx}$, $e_2=\er^{\ri k_c(-x+\sqrt{3}y)/2}$,
 $e_3=\er^{\ri k_c(-x-\sqrt{3}y)/2}$, and 
$\Phi=\Phi(\ld)\in\R^2$ is the eigenvector 
of $\hat{L}(k_c,\ld)$ to the eigenvalue $\mu_+(k_c,\ld)$. 
First we consider slowly varying complex amplitudes 
$A_j=A_j(t)$ of these modes, $j=1,2,3$, i.e., our 
ansatz reads 
 \begin{align}\label{newansatz}
   w= \sss A_i e_i \Phi+\frac12 \sss |A_i|^2 \phi_0 + \sss A_i^2
   e_i^2\phi_{1} + \sum_{1\le i < j \le 3} A_i \ov{A_j} e_i \overline{e_j}
   \phi_{2}+\text{c.c.}+\text{h.o.t.}, 
 \end{align} 
where $\ov A_j$ means the complex conjugate of $A_j$, c.c.~stand for 
the complex conjugate of all preceeding terms, and \text{h.o.t.} stands 
for higher order terms. This is taken as a weakly nonlinear 
expansion, i.e., the goal is to successively 
remove terms of order $A_j^m$, $m=1,2,3$ from 
the residual $L(\Delta)w{+}G(w)$, where the  
amplitudes $A_j$ are assumed to be small, though later we will use 
the expansion also for $\CO(1)$ amplitudes. 
The vectors  
$\phi_{0}$, $\phi_1$, and $\phi_{2}$ are introduced 
to remove quadratic terms at wave vectors $\fettk$ with 
 $k=|\fettk|=0$, $k=2k_c$ and $k=\sqrt{3}k_c$  from the residual.  

The calculations are best organized by writing \reff{dglredu0} in the form 
\begin{align}
   \partial_t w = L(\Delta) w+B(w,w)+C(w,w,w),\label{dglredu}
 \end{align}
 where $B$ and $C$
 are symmetric bilinear and trilinear forms. To remove terms 
of order $A_i A_j$ from \eqref{rd} we need
\huga{\label{a3} 
\begin{split}
&\phi_0(\ld)
   =-2\Lh(0,\lam)^{-1}B(\Phi,\overline{\Phi}), \qquad 
\phi_{1}(\ld)=-\Lh(2k_c,\lam)^{-1}B(\Phi,\Phi), \\
&\phi_{2}(\ld)=-2\Lh(\sqrt{3}k_c,\lam)^{-1} B(\Phi,\overline{\Phi}). 
\end{split}
} 
These terms arise due to quadratic interactions 
of the forms, e.g., $B(A_1e_1\Phi,\ov{A_1e_1\Phi})=|A_1|^2B(\Phi,\ov{\Phi})$, 
$B(A_1e_1\Phi,A_1e_{1}\Phi)=A_1^2\er^{2\ri k_c x}B(\Phi,\Phi)$, and 
$B(A_1e_1\Phi,\ov{A_2e_2\Phi})=A_1\ov{A_2}
\er^{\ri \frac{k_c} 2(3 x-\sqrt{3}y)}B(\Phi,\ov{\Phi})$. 
Allthough $\Phi\in\R^2$ in our case we keep the notation 
$\ov{\Phi}$ as this makes it easier to see where the respective terms come from. 
The matrices $\Lh(0,\lam)$, $\Lh(2k_c,\lam)$ and $\Lh(\sqrt{3}k_c,\lam)$ 
are invertible at least for $\lam$ not too far from $\lam_c$. 
From the Fredholm alternative we obtain the Landau ODE system as the 
solvability conditions for removing terms up to cubic order at the 
critical modes, namely 
\ali{
\text{at }e_1:\quad  &\partial_t A_1=f_1(A_1,A_2,A_3):=
c_1 A_1+c_2 \ov{A_2}\ov{A_3} +c_3 |A_1|^2A_1+c_4 A_1(|A_2|^2 +|A_3|^2), 
\nonumber \\
\text{at }e_2:\quad   &\partial_t A_2=f_2(A_1,A_2,A_3):=
c_1 A_2+c_2 \ov{A_1}\ov{A_3} +c_3 |A_2|^2A_2+c_4 A_2(|A_1|^2 +|A_3|^2), 
\label{lan}\\
\text{at }e_3:\quad    &\partial_t A_3=f_3(A_1,A_2,A_3):=
c_1 A_3+c_2 \ov{A_1}\ov{A_2} +c_3 |A_3|^2A_3+c_4 A_3(|A_1|^2 +|A_2|^2), 
\nonumber 
} 
with $c_1(\ld)=\mu_+(k_c,0,\ld)$,  
$c_2(\ld)= 2 \langle B(\overline{\Phi},\overline{\Phi}),\Phi^*\rangle$, 
$c_3(\ld)=\langle
 3C(\Phi,\Phi,\overline{\Phi})+2B(\overline{\Phi},\phi_{1})+2B(\Phi,\phi_{0}),\Phi^*\rangle$, and $c_4(\ld)= \langle
 6C(\Phi,\Phi,\overline{\Phi})+2B(\Phi,\phi_{2})
+2B(\Phi,\phi_{0}),\Phi^*\rangle$. 
Here $\Phi^*(\ld)$ is the adjoint eigenvector of $\hat{L}(k_c,\ld)$
 to the eigenvalue $\mu_+(k_c,\ld)$, normalized such that
 $\langle\Phi,\Phi^*\rangle=1$.   
At $\ov{e_j}$, $j=1,2,3$, we obtain the 
complex conjugate equations. 
See, e.g., \cite{gsk84, pomeau}, \cite[\S IV,A,1,a(iii)]{ch93}, \cite[\S2]{dsss03}, 
\cite[\S5]{hoyle},\cite{nepo06, pismen06}, which also explain 
that \reff{lan} is the generic form of the amplitude equations for 
hexagonal symmetry. The next step is the actual calculation of the 
coefficients. 

\brem\label{arem} In \reff{newansatz}, $\Phi=\Phi(k_c,\lam)$ varies 
with $\lam$ and is not fixed 
at $\lam=\lam_c$, which would be the more classical ansatz. 
Similarly, terms like \reff{a3} and the coefficients 
$c_2,\ldots,c_4$ are evaluated at $\lam$, not as usual at $\lam=\lam_c$. 
The formalisms are equivalent for $\lam\ra\lam_c$ 
with the differences hidden in 
the h.o.t.~in \reff{newansatz}. The reason why we always evaluate 
at $\lam$ is that we want to use the formalism also for 
$\lam_c-\lam=\CO(1)$ and then expect (and find) better approximations 
with $\Phi=\Phi(\lam)$. On the other hand, for the comparison 
with the numerics we want to keep the wave vectors fixed and 
thus evaluate at $k_c$, which approximately stays the most unstable wave 
number also for $\lam_c{-}\lam{=}\CO(1)$. 
\eex 
\erem 

\begin{figure}[!ht]
\bce{\small\begin{tabular}[t]{lp{50mm}lp{55mm}lp{45mm}}
(a)&(b)&(c)\\\ig[width=43mm, height=47mm]{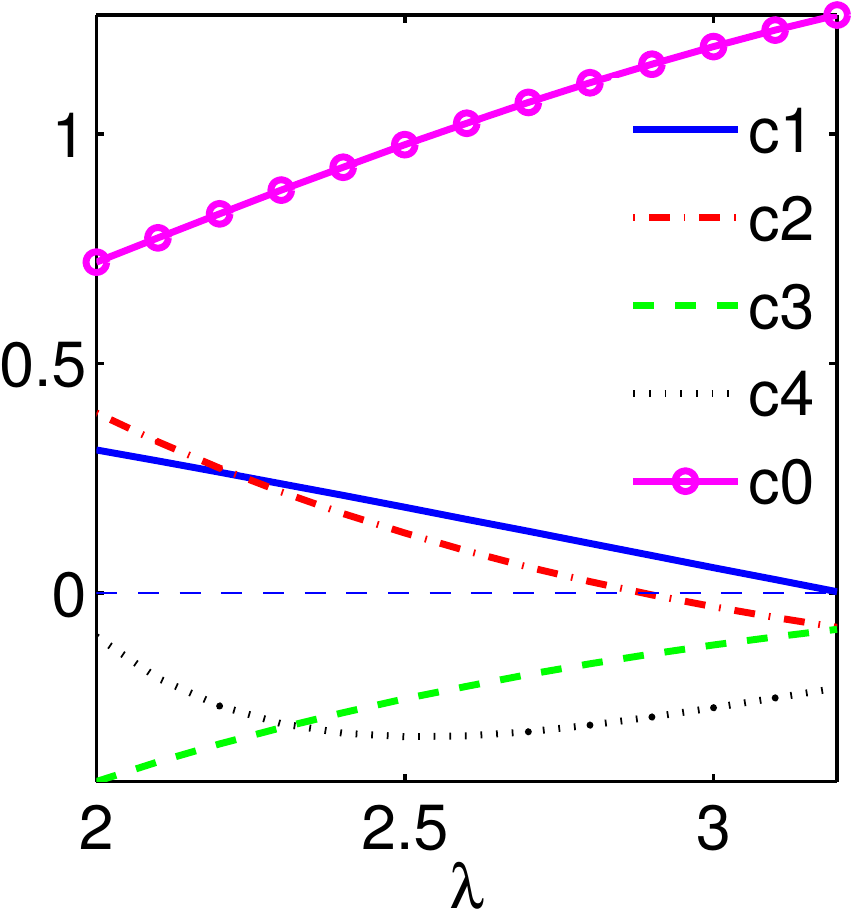}
&\ig[width=50mm,height=52mm]{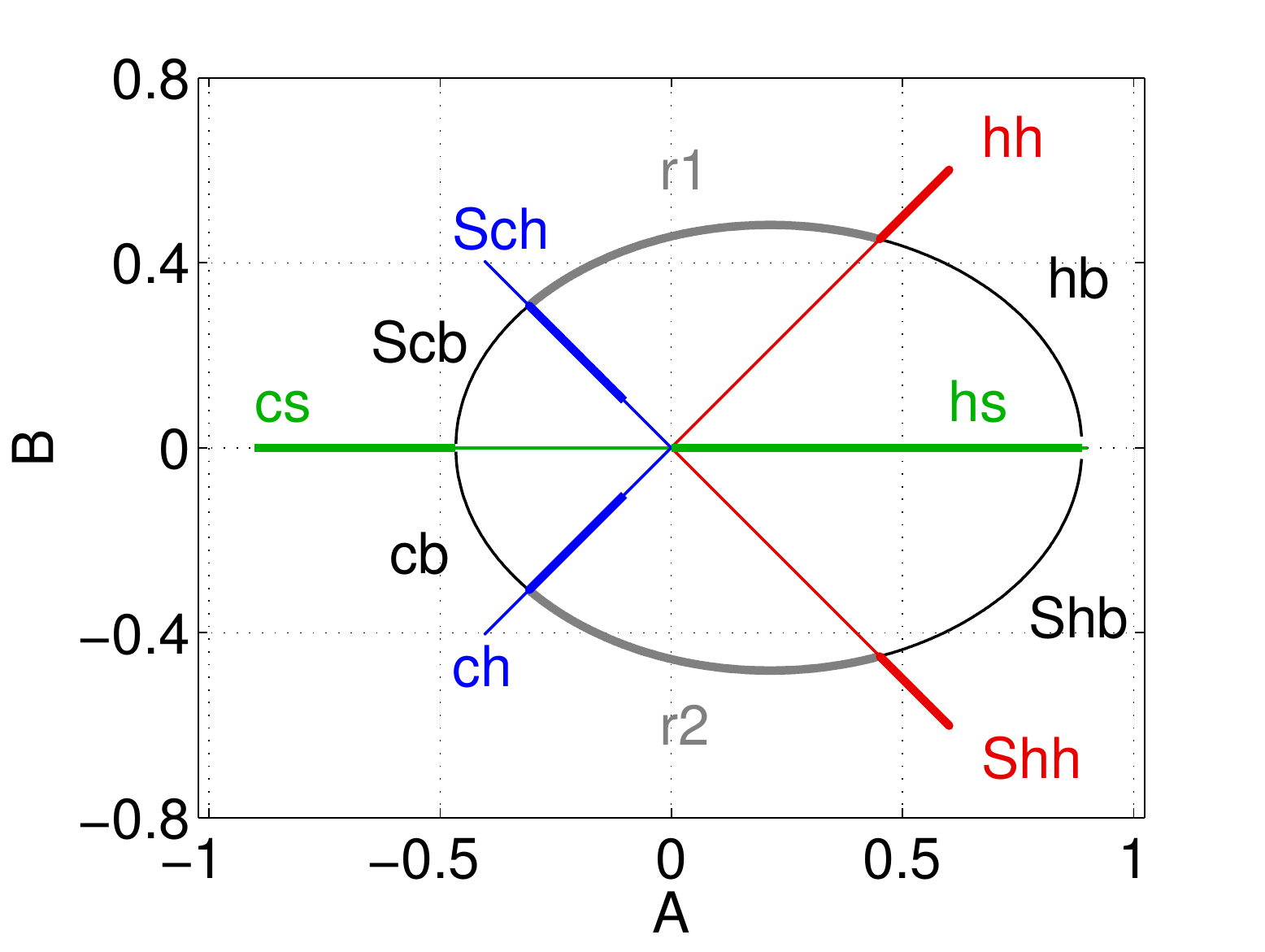}
&\ig[width=45mm,height=47mm]{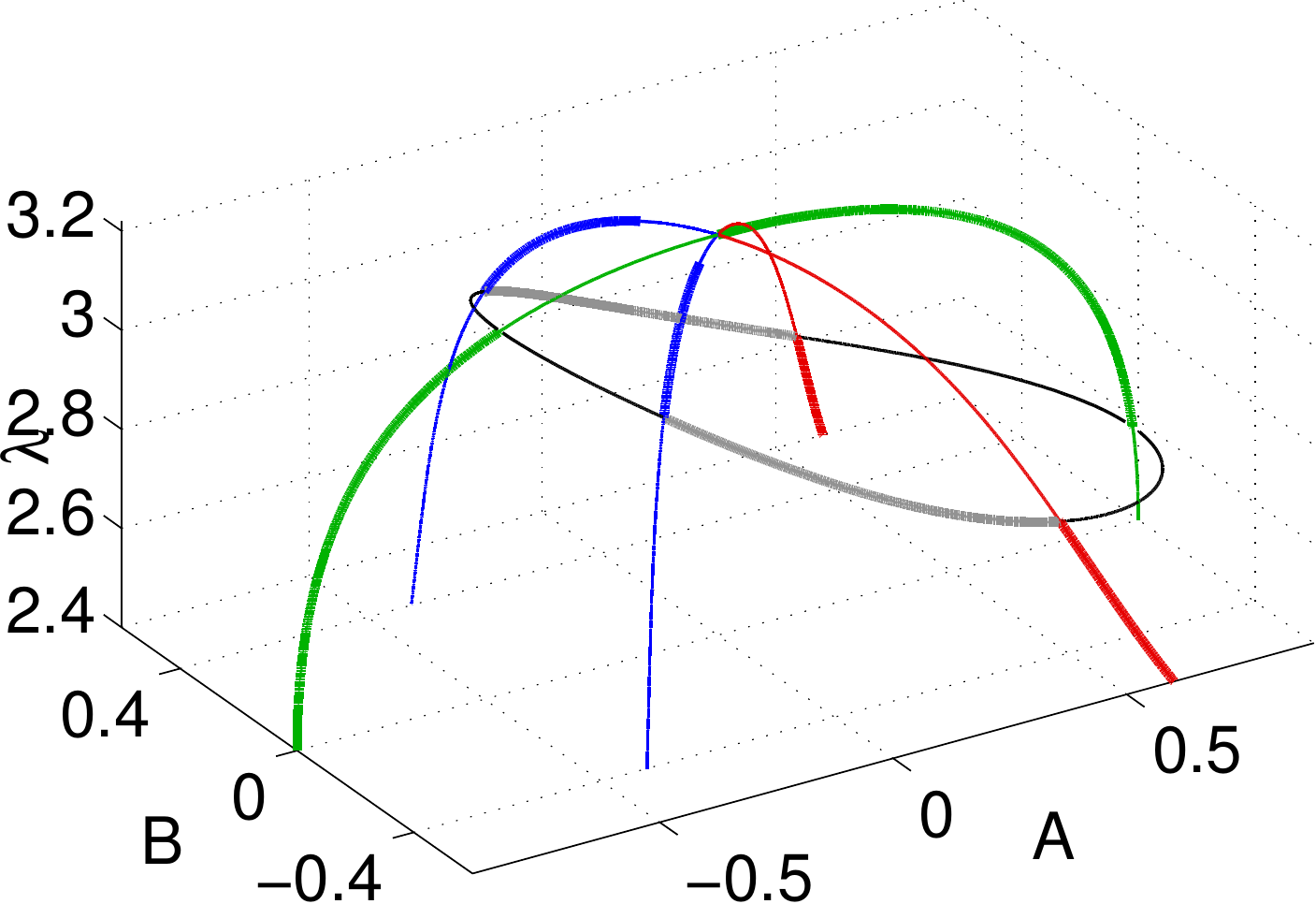}\\
(d)&(e)&(f)\\\ig[width=50mm,height=40mm]{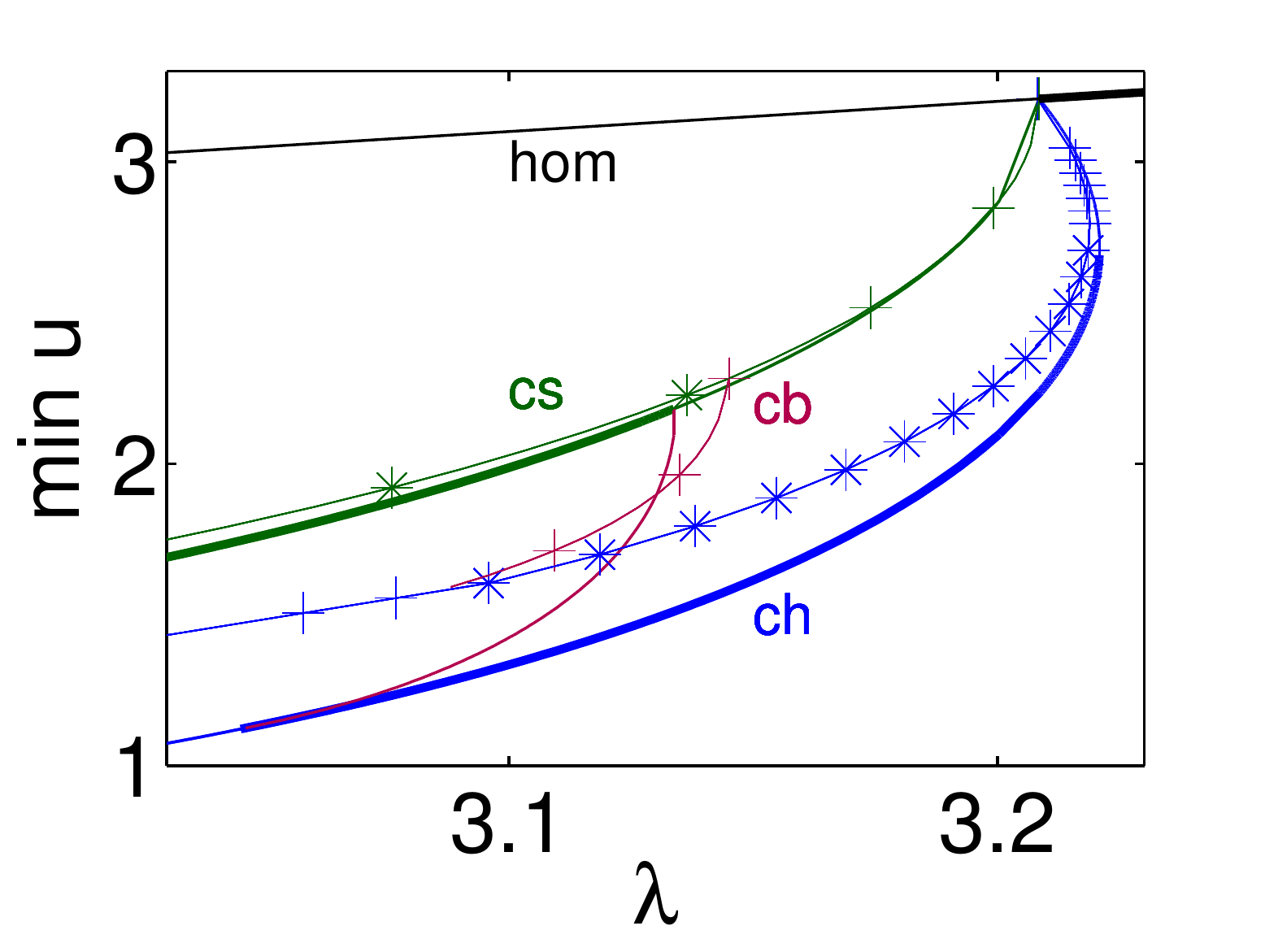}&
\ig[width=50mm,height=40mm]{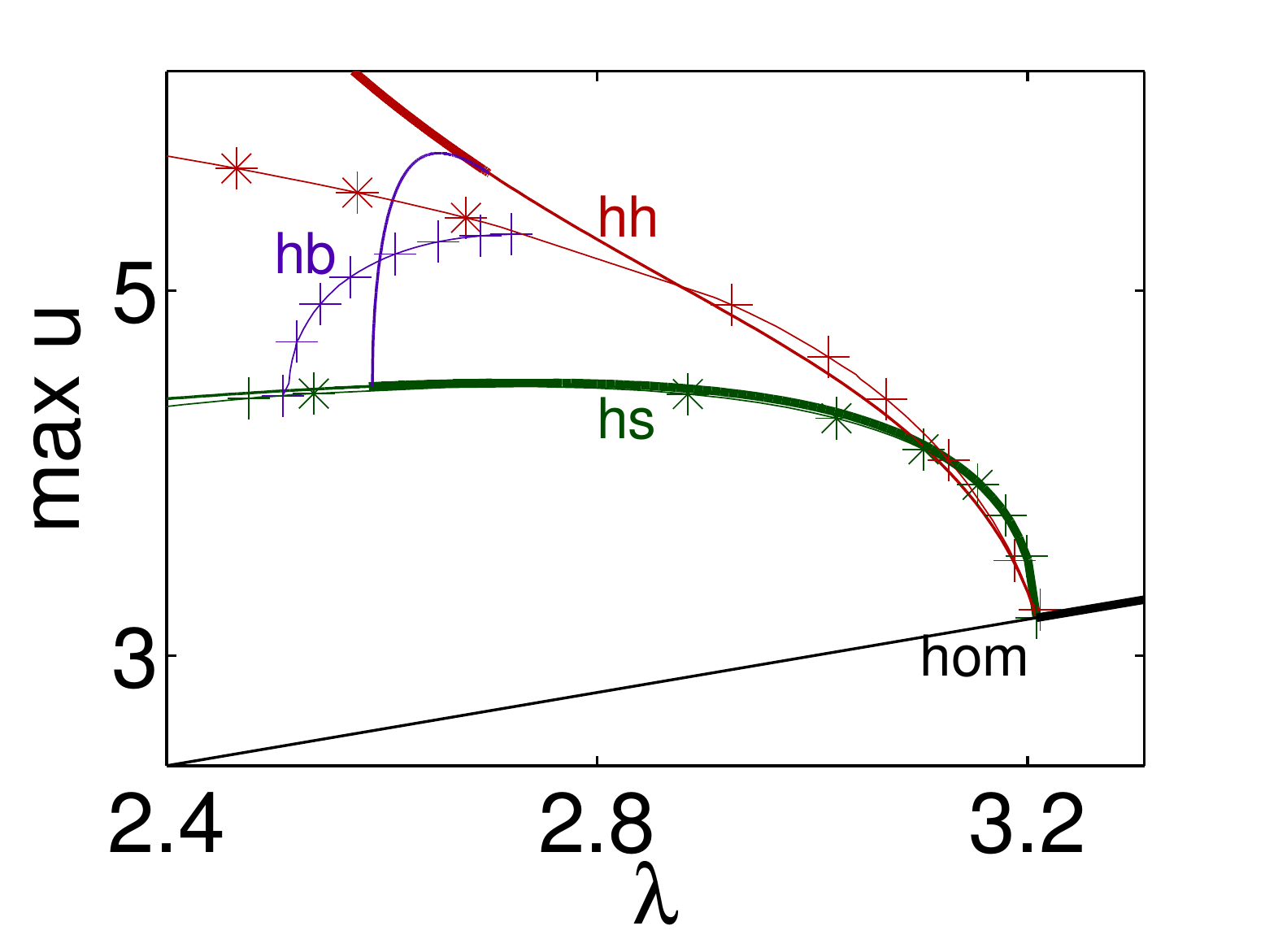}&
\quad\ig[width=40mm,height=40mm]{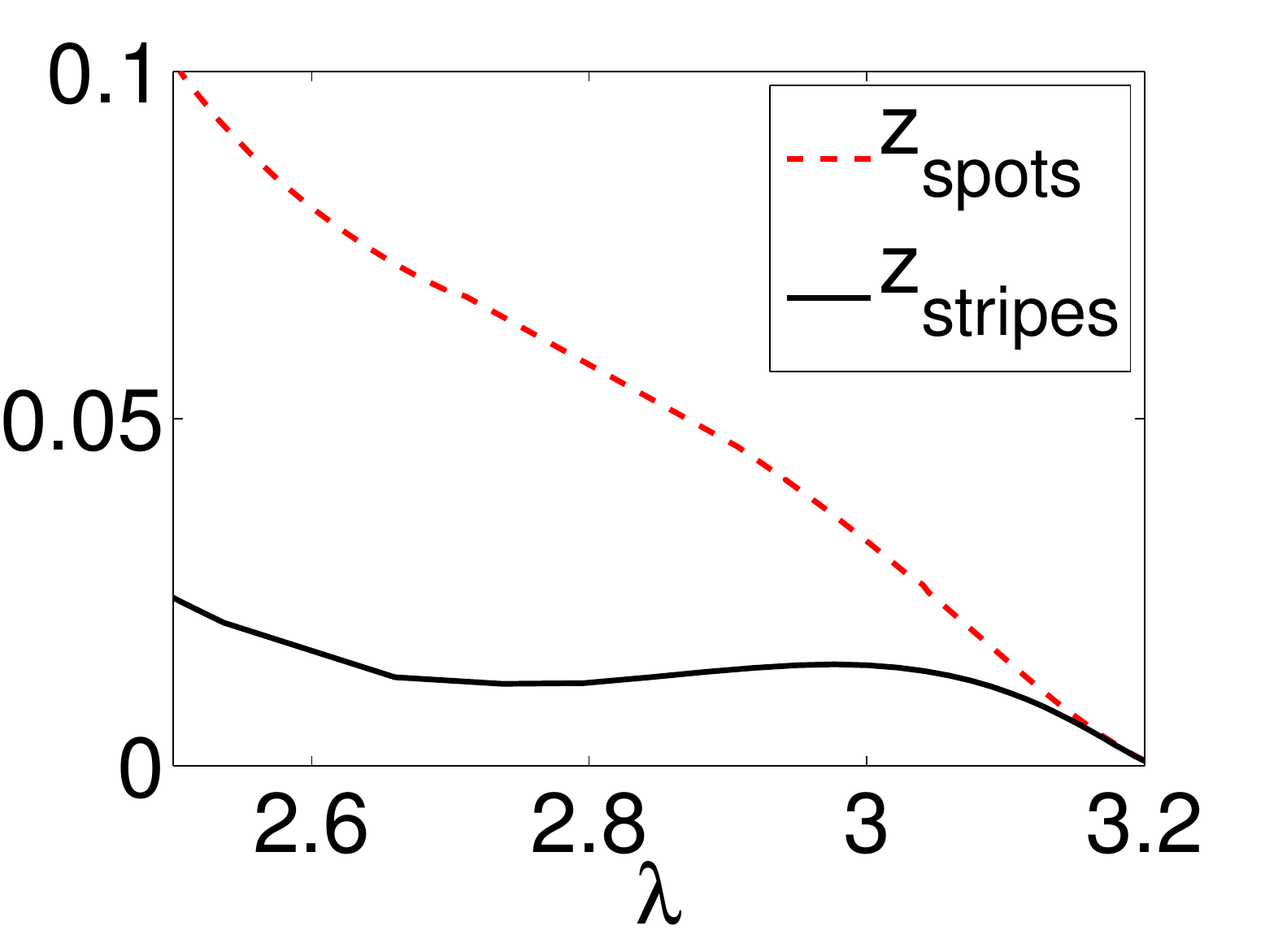}
\end{tabular}}\ece

\vs{-5mm}
\caption{{\small (a) Landau coefficients, and Ginzburg--Landau coefficient 
$c_0$, see \S\ref{glsec}. (b),(c) Landau bifurcation diagram in the invariant subspace 
$A_1=:A\in\R$ and $A_2=A_3=:B\in\R$, stability/instability 
indicated by thick/thin lines. 
(d),(e) 
Comparisons of numerical bifurcation diagram from Fig.~\ref{figsch} (indicated by * for stable 
points and + for unstable points)  with Landau bifurcation diagram 
near the cold and hot bistable regimes. 
(f) relative $L^1$ errors between hot numerical and Landau solutions defined in 
\reff{errdef}.}
\label{huf1}}
\end{figure}

From the coefficients $c_j(\lam)$, $j=1,\ldots,4$, 
see Fig.~\ref{huf1}(a), it follows that 
\ali{ 
T_\pm=\pm\sqrt{-\frac{c_1}{c_3}} \quad \text{and} \quad
   P_{\pm}=-\frac{c_2}{2(c_3+2c_4)}\pm\sqrt{\frac{c_2^2}{4(c_3+2c_4)^2}-\frac{c_1}{c_3+2c_4}}, \label{lanpoints} 
 } 
are real in the Turing-unstable range $\lam<\lam_c$, 
respectively for $\lam<\lam_{{\rm GLfold}}\approx \lam_{cb}\approx 3.22$. 
The triples $(T_+,0,0)$, $(T_-,0,0)$,
 $(P_+,P_+,P_+)$, and $(P_-,P_-,P_-)$ solve \eqref{lan} and via 
\reff{newansatz} generate hot stripes, cold stripes, hot hexagons, and cold hexagons,
 respectively. Mixed modes are obtained 
from setting $A_3=A_2$ and solving \reff{lan} for $(A_1,A_2)$. 

The formulas, together with the signs of $c_j(\lam)$, 
also explain the supercritical pitchfork bifurcation of stripes, and the 
transcritical bifurcation of the hexagons. Figures \ref{huf1}(b),(c) show 
the complete bifurcation diagram for stationary solutions of \reff{lan} in the invariant 
subspace $A_1=:A\in\R$ and $A_2=A_3=:B\in\R$, while (d),(e)  
compares the bifurcation diagram of stripes, hexagons and beans for \reff{lan} with the 
numerical bifurcation diagram from Fig.~\ref{figsch} for the full system 
\reff{rd}. Though qualitatively correct down to $\lam=2.4$, 
the approximation errors grow with $|\lam-\lam_c|$, as expected. 
Additionally to the $L^\infty$ like 
error already shown, for $i\in\{\text{stripes,hexagons}\}$ we define the 
relative $L^1$ errors 
\begin{align}\label{errdef}
  z_{i}=\sum\limits_{(x,y)\in \Omega_d}
    \norm{U_{\text{num},i}(x,y)-U_i(x,y)}_1\bigg/\sum\limits_{(x,y)\in
      \Omega_d} \norm{U_{\text{num},i}(x,y)}_1 
\end{align}
between the numerical solutions $U_{\text{num}}$ and 
approximate solutions $U=w^*+w$ from \reff{newansatz}, where 
$\Omega_d$ is the discretization used to calculate the numerical 
solution from Fig.~\ref{figsch}. 
It turns out that 
$z_{\text{stripes}}$ stays rather small also for $\lam_c-\lam=\CO(1)$, 
while $z_{\text{hex}}$ behaves worse, see Fig.~\ref{huf1}(f), and also 
Remark \ref{a2rem} for further comments. 

Another notable discrepancy between the Landau bifurcation diagram and the (numerical) bifurcation diagram 
for \reff{rd} 
is that in the former the (hot) stripes are stable all the way to the 
bifurcation point, while in the PDE this depends on the domain size. 
On a $1\times 1$ domain we obtain exactly the same stability as displayed 
in Fig.~\ref{huf1}(b),(c), and the instability near 
bifurcation of the stripes already on a $2\times 2$ domain 
as in Fig.~\ref{figsch} is due to a long zig-zag 
instability, which is not captured in our Landau ansatz. 
Similar remarks apply for instance to the 
rectangle branches, which are mostly unstable 
for \reff{rd} over a $2\times 2$ domain.

\subsection{Ginzburg--Landau formalism and fronts}
\label{glsec}
Motivated by the acceptable approximation of the hexagons, stripes and beans 
via \reff{newansatz} and \reff{lan} we proceed to use the 
Ginzburg--Landau (GL) reduction to predict the stationary fronts 
and localized patterns. 
If instead of $A_j=A_j(t)$ we assume that $A_j=A_j(t,x)$ are slowly varying 
functions also of $x\in\R$, then instead of the Landau system \reff{lan} 
we obtain the Ginzburg--Landau system 
\ali{
  &\partial_t A_1= c_{0}\partial_{x}^2 A_1+f_1(A_1,A_2,A_3), \quad 
\partial_t A_j=\frac{c_{0}}{4}\partial_{x}^2 A_j+f_j(A_1,A_2,A_3),\ j=2,3, 
\label{gl}
} 
where $c_{0}(\ld)=-\frac12\partial_{k_1}^2\mu_+((k_c,0),\ld)>0$. 
See \cite{gs94, bsvh95, gs99b, alex} for background on this formal procedure, 
and for so called attractivity and approximation theorems which 
 estimate the 
difference between a true solution of \reff{rdt} and 
an approximation described by \reff{newansatz} and \reff{gl}, 
 close to bifurcation. 
These theorems involve some small amplitude assumption 
for $(u,v)$, related slow scales for $t$ and $x$ in $A_j(t,x)$, and, 
for the present case of three resonant modes, a suitable scaling for 
the quadratic interactions, i.e., small $c_2$. 

\begin{figure}[!ht]
\bce
{\small
\begin{tabular}[t]{lp{35mm}lp{35mm}lp{35mm}lp{35mm}}
(a)&(b)&(c)&(d)\\
\ig[width=35mm,height=37mm]{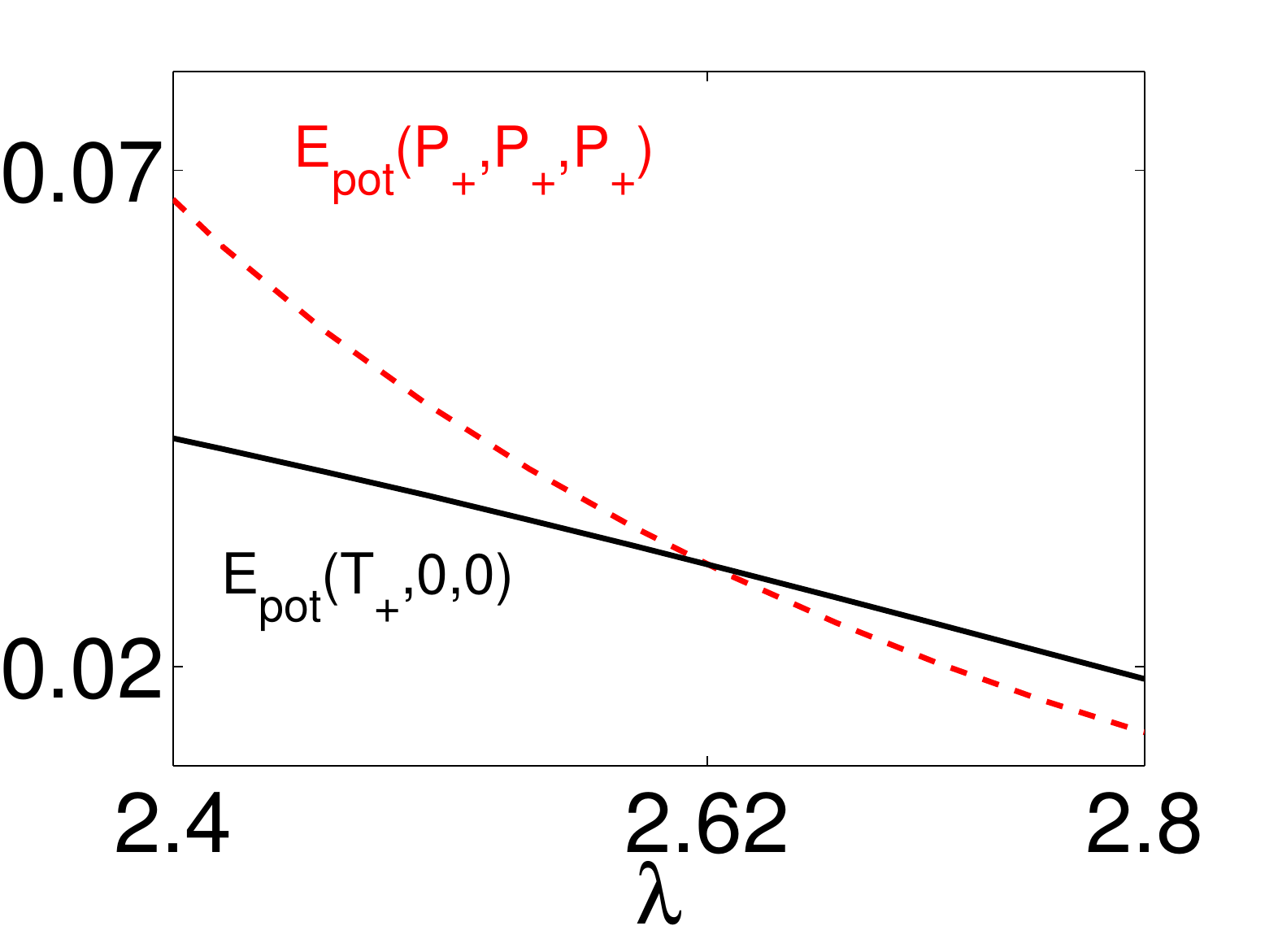}&
\ig[width=35mm,height=35mm]{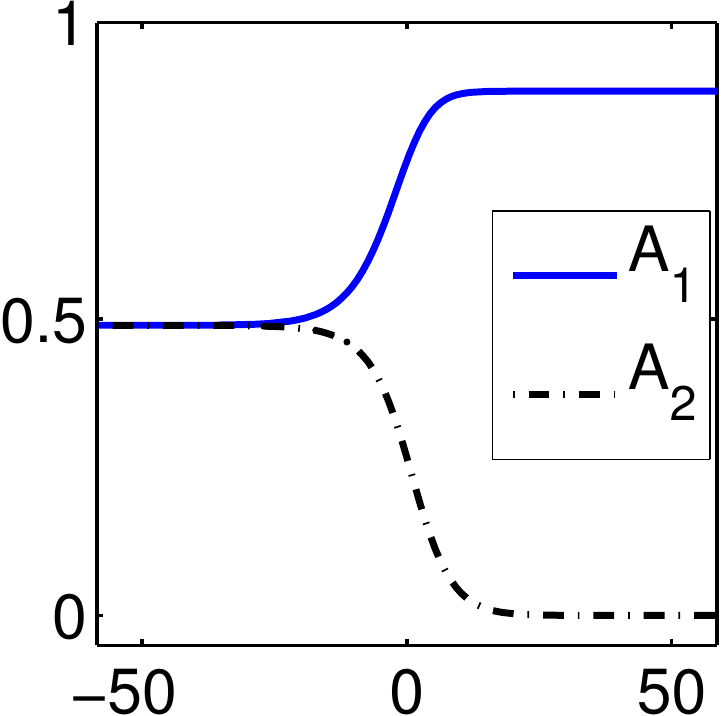}&
\ig[width=35mm,height=37mm]{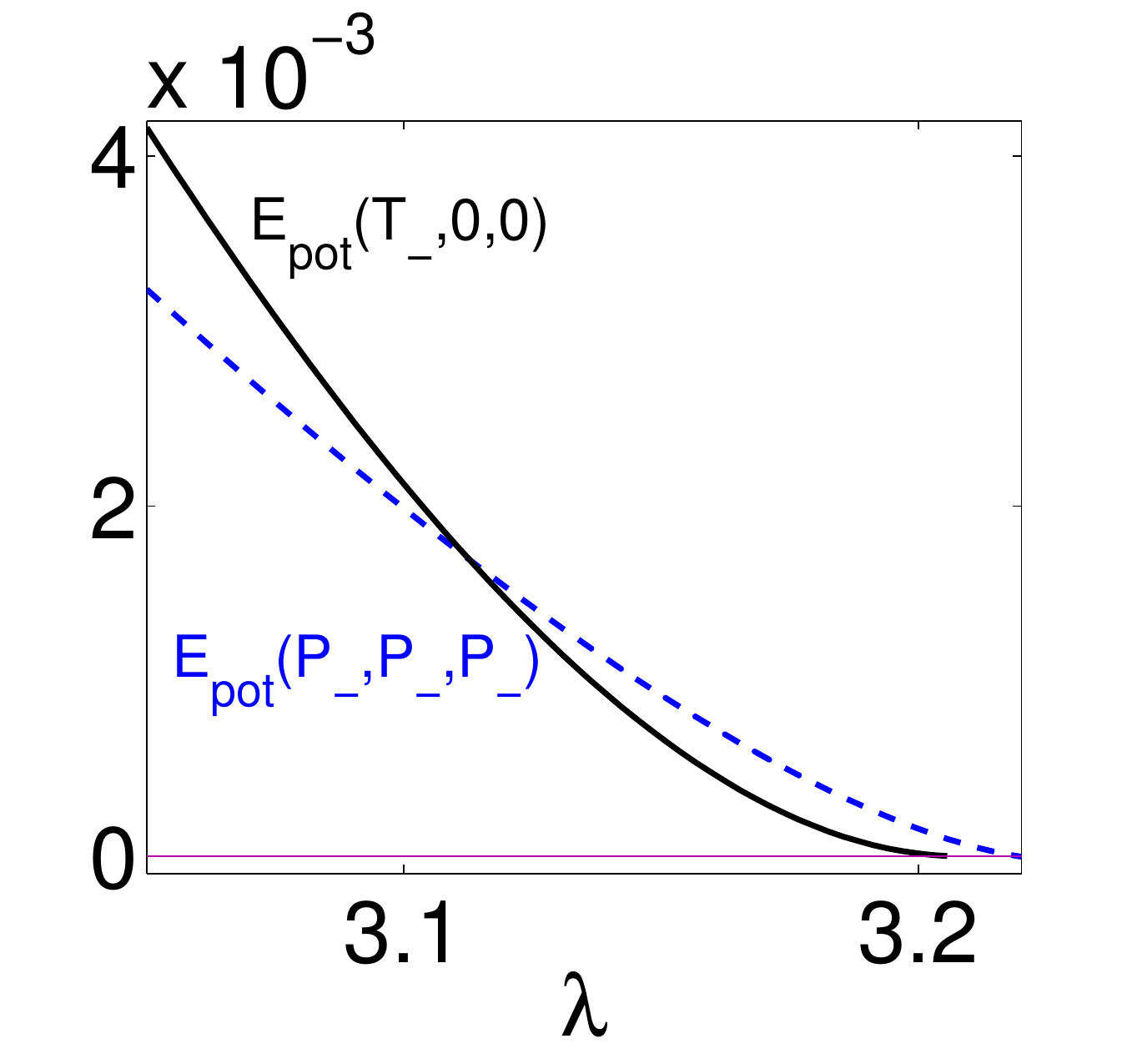}&
\ig[width=40mm,height=35mm]{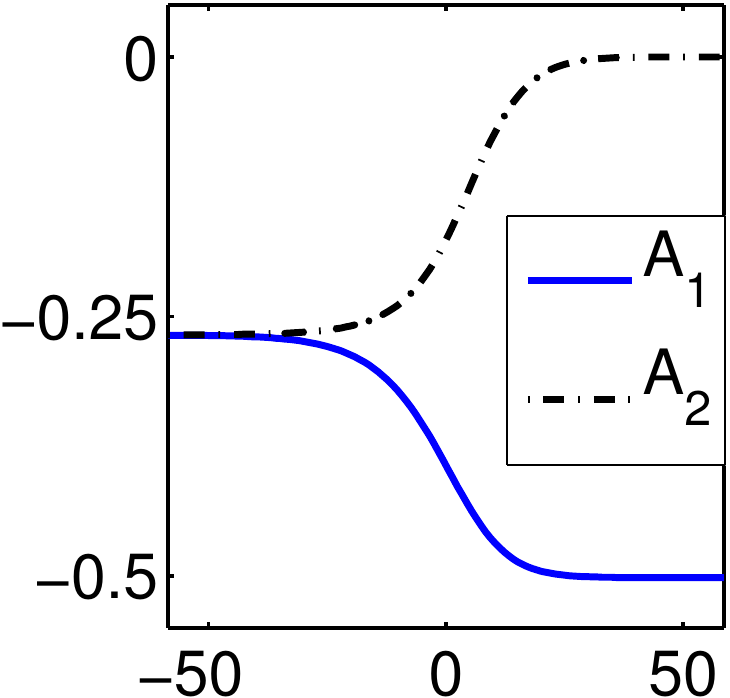}
\end{tabular}
}
\ece

\vs{-5mm}
\caption{{\small (a), (c) Energies of hexagons and stripes near 
the hot resp.~cold bistable regimes; their intersections 
define the respective Maxwell points $\lam=2.62$, $\lam=3.11$, and 
the zero of $E_{{\rm pot}}(P_-,P_-,P_-)$ defines 
the third Maxwell point $\lam=3.219$.  
Hot (b) and cold (d) stationary Ginzburg--Landau 
fronts between stripes and hexagons 
at the respective Maxwell points, for comparison 
both over domain $[-12\pi/k_c, 12\pi/k_c]$, $A_3=A_2$. }  \label{huf5}}
\end{figure}

Here we again want to use the GL system \reff{gl} at an $\CO(1)$ 
distance from $\lam_c$ and find a necessary condition 
for stationary fronts between hexagons like $A\equiv (P_+,P_+,P_+)$ 
and stripes like $A\equiv (T_+,0,0)$ in \reff{gl}. 
Thus we consider the stationary Ginzburg--Landau system 
\ali{
  c_{0}\partial_{x}^2 A_1{+}f_1(A_1,A_2,A_3)=0, \quad 
  \frac{c_{0}}{4}\partial_{x}^2 A_j{+}f_j(A_1,A_2,A_3)=0,\ j=2,3, 
  \label{glsd}
} 
as a dynamical system in the spatial variable $x$. 
Now restricting to real amplitudes $A_j$, 
the total energy of \eqref{glsd} is given by $\etotal=\ekin+\epot$,
where 
\alinon{ &\ekin=\frac{c_0}{2} \lkl(\partial_{x}
  A_1)^2+\frac14(\partial_{x} A_2)^2+\frac14(\partial_{x} A_3)^2 \rkl, 
\quad \text{and} \\
&\epot=\sum_{i=1}^3 \lkl\frac{c_1}2 A_i^2+\frac {c_3} 4 A_i^4\rkl+c_2A_1A_2A_3 
  +\frac{c_4}{2}(A_1^2A_2^2+A_1^2A_3^2+A_2^2A_3^2) 
}
are the kinetic and potential energy, respectively. 
Then $\frac{d}{dx}\etotal=0$, i.e., 
$\etotal$ is conserved. Thus,  a 
necessary condition for, e.g., a heteroclinic orbit $\Afront$ between 
$(P_+,P_+,P_+)$ and $(T_+,0,0)$ to exist in \reff{glsd} 
is $\epot(T_+,0,0)=\epot(P_+,P_+,P_+)$. Again we first focus on the 
hot bistable range and in Fig.~\ref{huf5}(a) 
plot $\epot(T_+,0,0)$ and $\epot(P_+,P_+,P_+)$.  
Their intersection defines the so called (hot) Maxwell point $\ld_m$. 
Though we refrain from discussing the 
general energy landscape and dynamics of \reff{glsd}, it turns out 
that at least numerically 
the necessary condition $\lam=\lam_m$ is also sufficient. 
Figure \ref{huf5}(b) shows a stationary front for \reff{gl} 
(with Neumann boundary conditions), which can either 
be obtained from time evolution of \reff{gl} with a suitable initial guess, 
or from solving the stationary boundary value problem. 
See also \cite{MNT90} for some more analysis for front solutions of \reff{glsd}, 
including some implicit solution formulas. 
On the other hand, for $\lam<\lam_m$ with $\epot(T_+,0,0)<\etotal(P_+,P_+,P_+)$ 
we obtain a front for \reff{gl} {\em travelling} towards higher energy, i.e.,
$(P_+,P_+,P_+)$ invades $(T_+,0,0)$, and vice versa for $\lam>\lam_m$. 

Thus, from the GL approximation  {\em standing} 
fronts $\Uhet$ between different patterns are only predicted at precisely 
$\lam=\lam_m$. 
In the physics literature, e.g., \cite{pomeau}, the basic argument 
for the existence of $\Uhet$ in an interval around $\lam_m$ is 
that the patterns create an effective periodic potential 
which yields a pinning of fronts. 

It is not obvious whether 
the stationary GL system \reff{glsd} has homoclinic solutions $\Ahom$ with, 
say, $\Ahom(x)\ra (P_+,P_+,P_+)$ as $x\ra\pm\infty$, and that 
pass near $(T_+,0,0)$ near $x=0$. 
However, (suitably shifted) 
fronts $\Afront(x_0+\cdot)$ and ``backs'' $\Aback(x_1+\cdot)$ 
with $\Aback(x)=\Afront(-x)$ can be glued together to 
give approximate homoclinics with long plateaus but also with 
some dynamics in time which can be expected to be 
exponentially slow in the separation distance between 
$\Afront$ and $\Aback$, see, e.g., \cite{CP90}. In fact we can 
generate almost stationary pulses numerically, but eventually 
the solution decays to a homogeneous rest state $(T_+,0,0)$ 
or $(P_+,P_+,P_+)$. Similarly, \reff{glsd} may have 
periodic orbits which stay close to $(P_+,P_+,P_+)$ 
resp.~$(T_+,0,0)$ over very long $x$-intervals. 
Clearly, for these ``approximate homoclinics'' 
a similar pinning argument as for the heteroclinics should apply 
and predict the existence of localized patterns as in Fig.~\ref{snakea} 
near the Maxwell point $\lam_m$. 

In \cite{chapk09, dean11} it is worked out mathematically  
that the pinning and hence also the snaking are exponentially 
small effects and thus cannot be predicted at any order by 
Ginzburg--Landau type asymptotic expansions alone. However, as already 
indicated in Remark \ref{dawrem2}, using beyond all order asymptotics 
snaking in the model problems 
\reff{qcshe} and \reff{cqshe} can be 
described very accurately. See also, e.g., \cite{kno2008,strsnake} 
 and the references therein for alternative arguments explaining the 
snaking via so called heteroclinic tangles in the 
spatial dynamics formulation of again \reff{qcshe} and \reff{cqshe}.

The same pinning arguments apply to the cold bistable range 
(Fig.~\ref{huf3}), but the difference to the 
hot range is that the bistability resp.~subcriticality as measured 
by $\eps$ in Remark \ref{dawrem2} in 
the cold range is much smaller than in the hot range. Here we repeat, 
from a slightly different point of view, 
the heuristic argument, why this makes the snaking region very small  
and thus snaking impossible to detect in our numerics. 
Small $\eps$ in the cold range relates to the fact that the dependence of  
$|E_{{\rm pot}}(T_-,0,0)-E_{{\rm pot}}(P_-,P_-,P_-)|$ on $\lam$ in the cold 
range is much flatter than 
that of $|E_{{\rm pot}}(T_+,0,0)-E_{{\rm pot}}(P_+,P_+,P_+)|$ on $\lam$ 
in the hot range. 
A short calculation then yields that relatively to the diffusion 
constant $c_0$ the vector fields 
$f_1,f_2$ in \reff{glsd} are smaller in the cold range than in the hot range. 
Therefore, cold Ginzburg--Landau fronts are 
flatter and wider than hot ones, cf.~Fig.~\ref{huf5}(b) and (d). Thus, if 
for convenience we write $u(x,y;\lam)=\sum_{j=1}^3 A_j(\eps x;\lam)
e_i\Phi(k_c)+\cc$ 
for a localized pattern $u$, assuming that 
the $A_j$ vary on an $\CO(1)$ scale, then 
there is a stronger separation of 
scales in the cold range than in the hot range. 
Consequently, in the cold range the Fourier transform of $u$ in $x$ 
is more localized around the pertinent $k_j$ due to  
\huga{\label{pfour} 
\hat{u}(k,y;\lam)=\frac 1 \eps\left(\sum_{j=1}^3 
\hat{A}\left(\frac{k-k_j}{\eps};\lam\right)\er^{\ri l_j y}+\ccf\right)\Phi(k_c), 
}
where $k_1=k_c, k_{2,3}=-\frac 1 2 k_c$, $l_1=0, 
l_{2,3}=\pm \frac{\sqrt{3}}2$. 
 
On the other hand, repeating the Ginzburg--Landau analysis leading 
to \reff{gl} with the $x$--wave number $k$ as a 
parameter, $k\approx k_c$, we obtain a wave number dependent 
Ginzburg--Landau energy $\etotal(\lam,k)$ for patterns. Then 
$\frac{d}{dx}\etotal=0$ along a localized pattern can 
be seen as a selection principle for the average wave number $k$, or conversely 
defines a function $\lam=\lam(k)$ with  generically
$\lam'(k)\ne 0$, see, e.g., \cite{BuDa12}. 
Therefore, if $\hat{u}(k;\lam)$ is only weakly localized in the different 
$k_j$ (i.e., $\eps$ not very small in \reff{pfour}), then we  
expect a strong $\lam$ dependence in the branch of localized patterns, 
and hence that $\lambda$ varies over a significant range during growth of 
a localized pattern, and we may expect $\lam$ to snake back and forth 
during that growth. 
Conversely, for small $\eps$ in \reff{pfour} we expect a narrow 
$\lam$ range for some growing pattern, as in Fig.~\ref{huf3} and \ref{huf4}, 
and snaking is difficult to detect numerically.  See also \S\ref{sigsec} for more 
quantitative arguments. 

\brem\label{a2rem} 
An ad hoc way to derive a Landau system like \reff{lan} from 
a system like \reff{rd} is to {\em only} consider the solvability 
conditions \reff{lan} at the critical modes without first 
removing the $A^2$ residual at second harmonics. This amounts to  
the {\em first order ansatz} 
$w=\sum_{i=1}^3 A_i e_i \Phi+\text{c.c.}$, or in other 
words setting $\phi_0,\phi_1,\phi_2=0$ in \reff{a3}. We denote 
the new coefficients in \reff{lan} by $c_{31}$ and $c_{41}$.
This is obviously 
simpler than \reff{newansatz}, but not formally consistent. 
However, at order $\CO(1)$ distance from the bifurcation 
the first order ansatz {\em may} give a better approximation of solutions than 
the ansatz \reff{newansatz}, 
since \reff{newansatz} represents only an asymptotic 
expansion, and not the first terms in some convergent series.

It turns out that in some sense this is indeed the case for \reff{rd}, 
and this can be used to improve the prediction of the Maxwell point. 
While in the hot bistable range the first order ansatz gives a much larger error
for the stripes, its error for the hot hexagons is in fact
smaller than the one in Fig.~\ref{huf1}(c). The idea is to use the
coefficients $c_3(\lam)$ obtained from \reff{newansatz} for the stripes 
and the coefficients $c_{31}(\lam), c_{41}(\lam)$ 
from the first order ansatz 
for the hexagons in a ``mixed'' Ginzburg--Landau system that 
retains the variational structure. 
Thus, let $S$ be
$T_+$ and $H$ be $P_+$ determined by using \reff{newansatz} and
the first order ansatz, respectively, and consider the system
\ali{
  c_{0}\partial_{x}^2 A_1 +c_1 A_1+c_2 A_2 A_3 +\lkl c_{3} \lkl 
\frac{A_1-H}{S-H}\rkl+c_{31} \lkl \frac{A_1-S}{H-S}\rkl\rkl A_1^3+c_{41} A_1(A_2^2 +A_3^2)&=0, \nonumber \\
  \frac{c_{0}}{4}\partial_{x}^2 A_2 +c_1 A_2+c_2 A_1 A_3 +\lkl 
c_{3}\lkl \frac{A_2-H}{S-H}\rkl+c_{31} \lkl \frac{A_2-S}{H-S}\rkl \rkl A_2^3+c_{41} A_2(A_1^2 +A_3^2)&=0, \label{glsd2}\\
  \frac{c_{0}}{4}\partial_{x}^2 A_3 +c_1 A_3+c_2 A_1 A_2 +\lkl c_{3}
  \lkl \frac{A_3-H}{S-H}\rkl+c_{31} \lkl \frac{A_3-S}{H-S}\rkl \rkl
  A_3^3+c_{41} A_3(A_1^2 +A_2^2)&=0. \nonumber }  In the hot bistable
range we have  $S > H$ such that 
\eqref{glsd2} is well-defined in this
range. Moreover, \reff{glsd2} has again a conserved energy 
$\etotal=\ekin+\epot$ where now 
$$
\epot{=}\sum_{i=1}^3 \frac{c_1}2 A_i^2{+}c_2A_1A_2A_3
+\sum_{i=1}^3 A_i^4\left(c_{3}\frac{\frac15 A_i -\frac14 H}{S-H}
+c_{31}\frac{\frac15 A_i -\frac14 S}{H-S}\right)
{+}\frac{c_4}2(A_1^2A_2^2+A_1^2A_3^2+A_2^2A_3^2).
$$
Using this energy to calculate the hot Maxwell--point for \reff{glsd2} 
we find $\lam_m=2.67$ which is more in the center of the snaking region 
than $\lam_m=2.62$. Thus, using the additional information that 
the hexagons are better approximated by the first order ansatz than by 
\reff{newansatz} we can obtain a better prediction for the Maxwell point. 
\eex 
\erem

\subsection{A modification with cold snakes and 1D snaking}\label{sigsec}
In \S\ref{crsec} and before Remark \ref{a2rem} we conjecture that the 
observed non snaking of cold connecting 
 branches is due to the weak difference in 
(Ginzburg--Landau) energy of the associated patterns, 
which are also related to the narrow range of bistability between 
the cold hexagons and the cold stripes or the homogeneous solution, 
respectively, and which we expect to yield exponentially narrow 
snaking regions, cf.~\reff{swform}. 
To study this quantitatively we now consider a modification of \reff{rd}. 
Clearly, many such modifications are possible, for instance using 
the Selkov--Schnakenberg parameters $a,b$, or varying the diffusion 
constant $d$, which we fixed to $d=60$. 

\begin{figure}[!ht]
\bce 
  \begin{minipage}{0.29\textwidth}
    (a)\\
    \ig[width=1\textwidth, height=45mm]{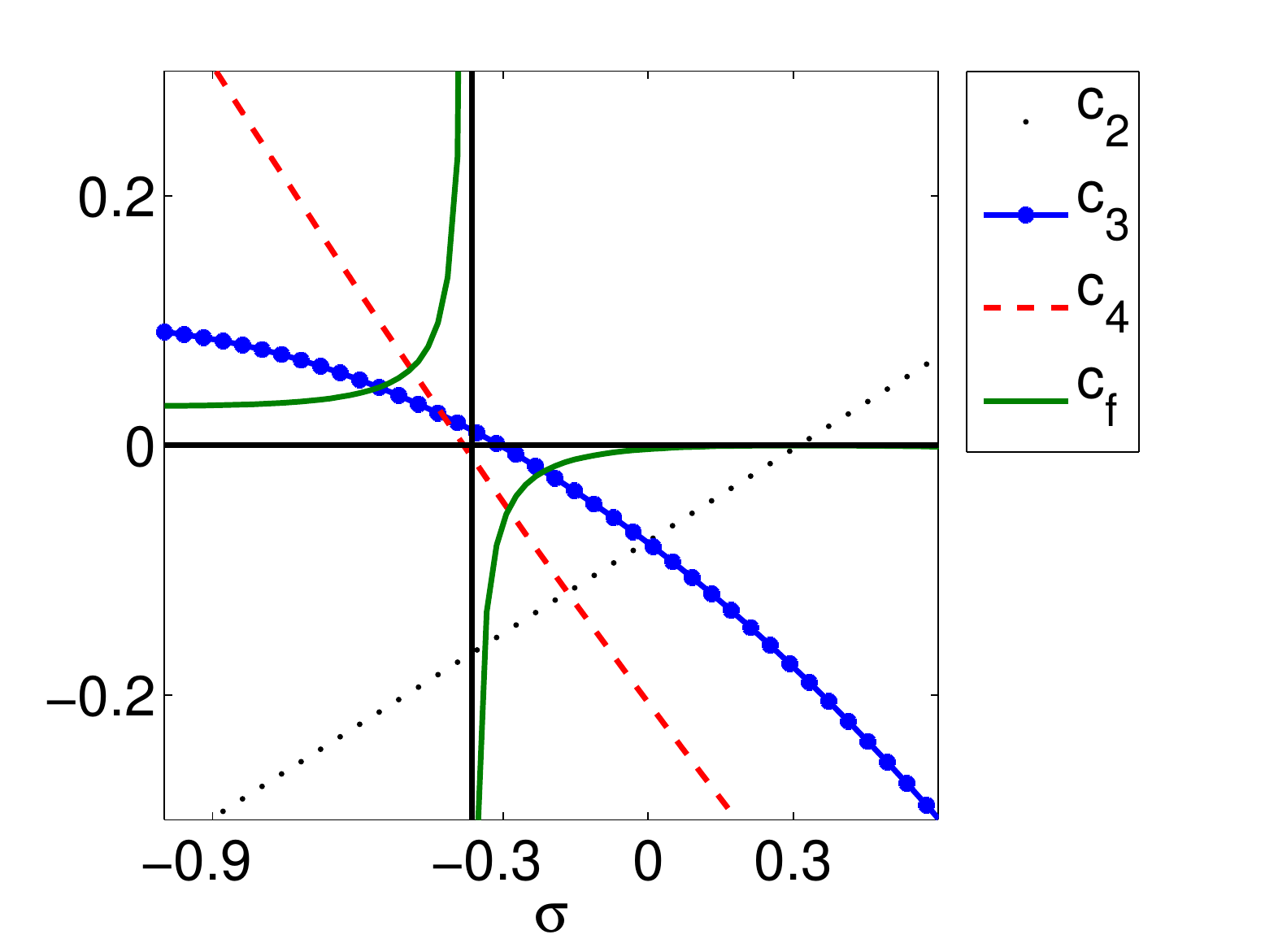}
  \end{minipage}
  \begin{minipage}{0.29\textwidth}
    (b)\\\ig[width=1\textwidth, height=45mm]{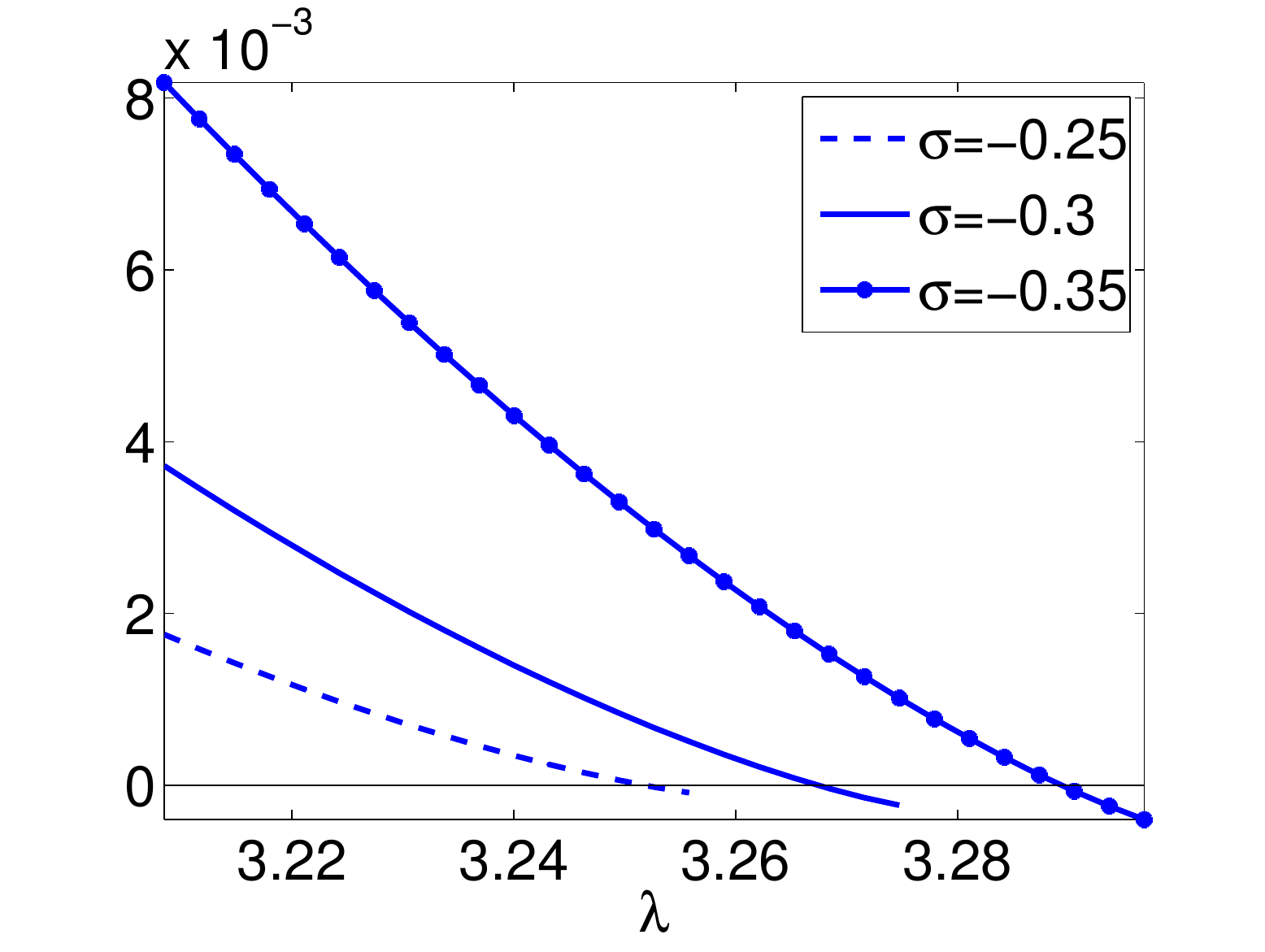}
  \end{minipage}
 \begin{minipage}{0.4\textwidth}
    (c)\\\ig[width=1\textwidth]{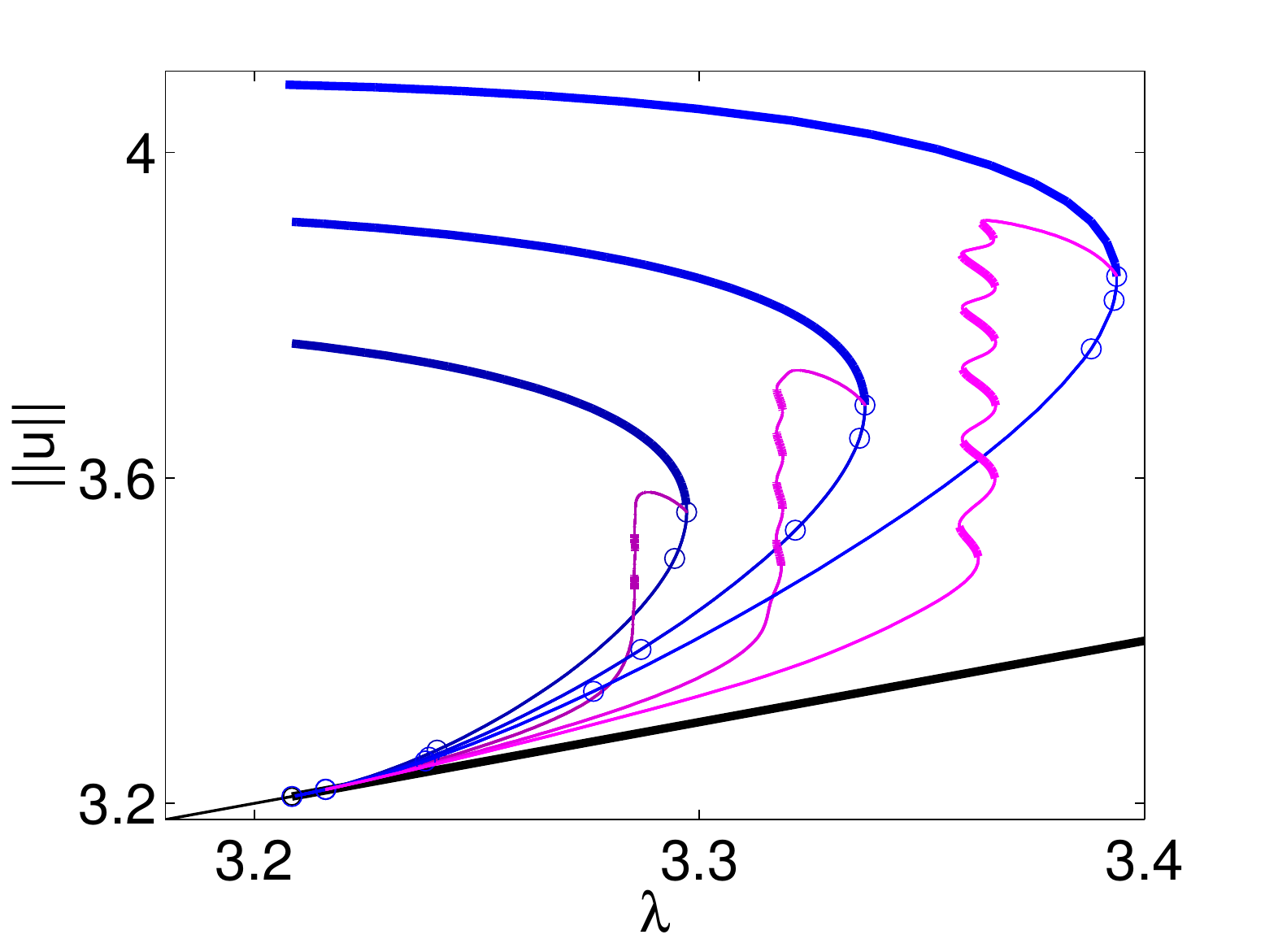}
  \end{minipage}
\\
  \begin{minipage}{0.29\textwidth}
    (d)\\\ig[width=1\textwidth,height=45mm]{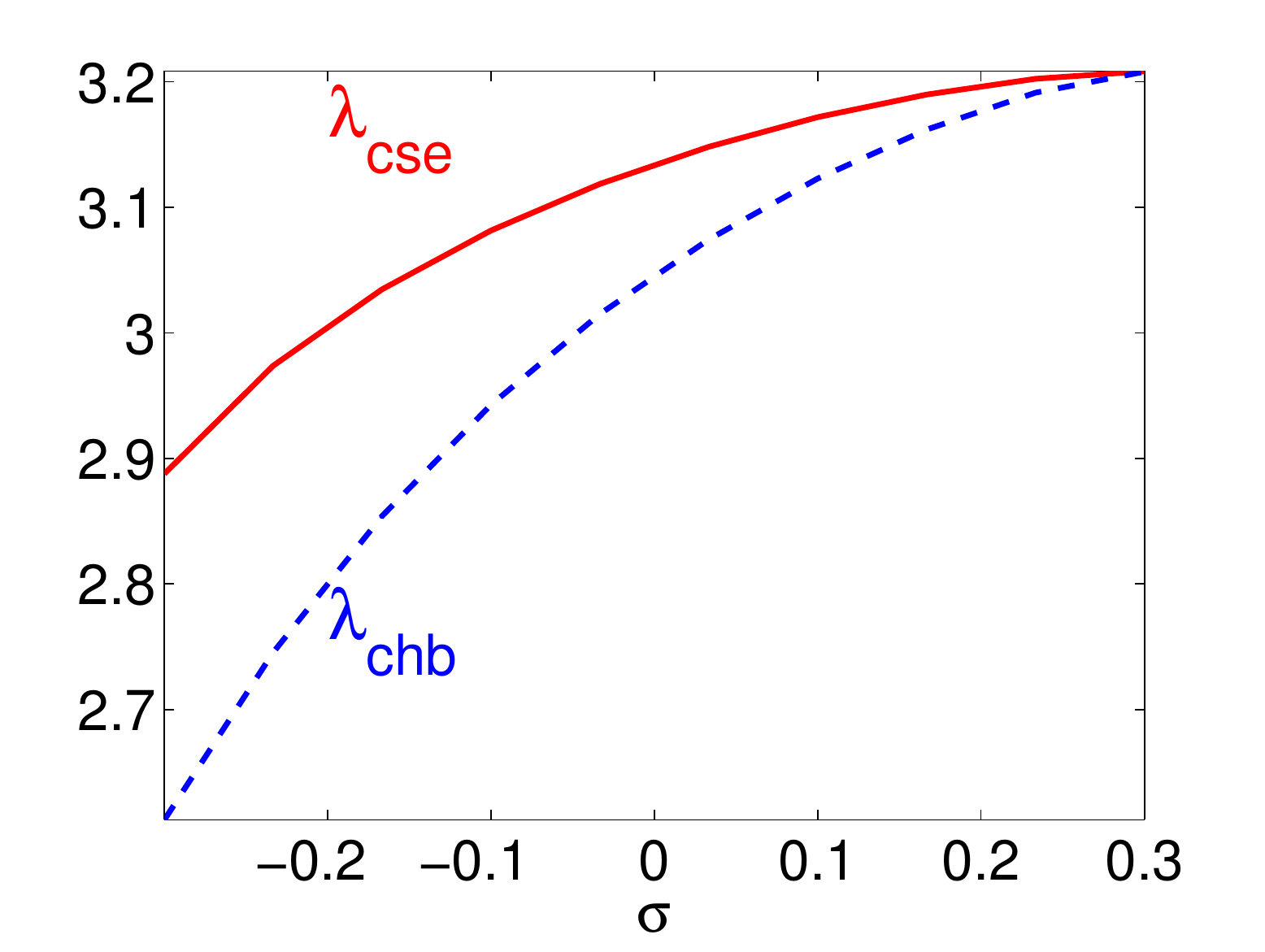}
  \end{minipage}
  \begin{minipage}{0.29\textwidth}
    (e)\\\ig[width=1\textwidth,height=45mm]{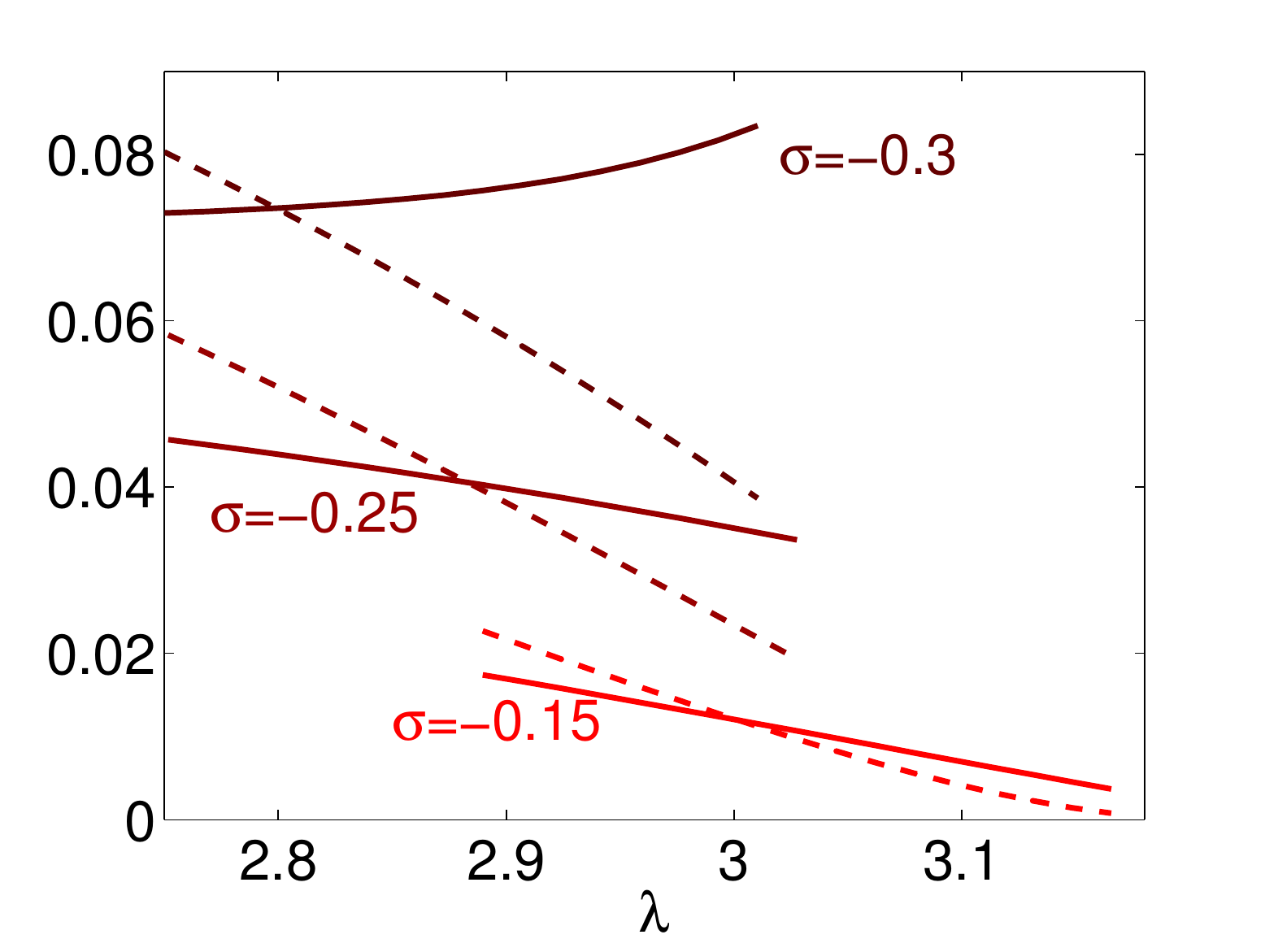}
  \end{minipage}
\begin{minipage}{0.4\textwidth}
    (f)\\\ig[width=1\textwidth]{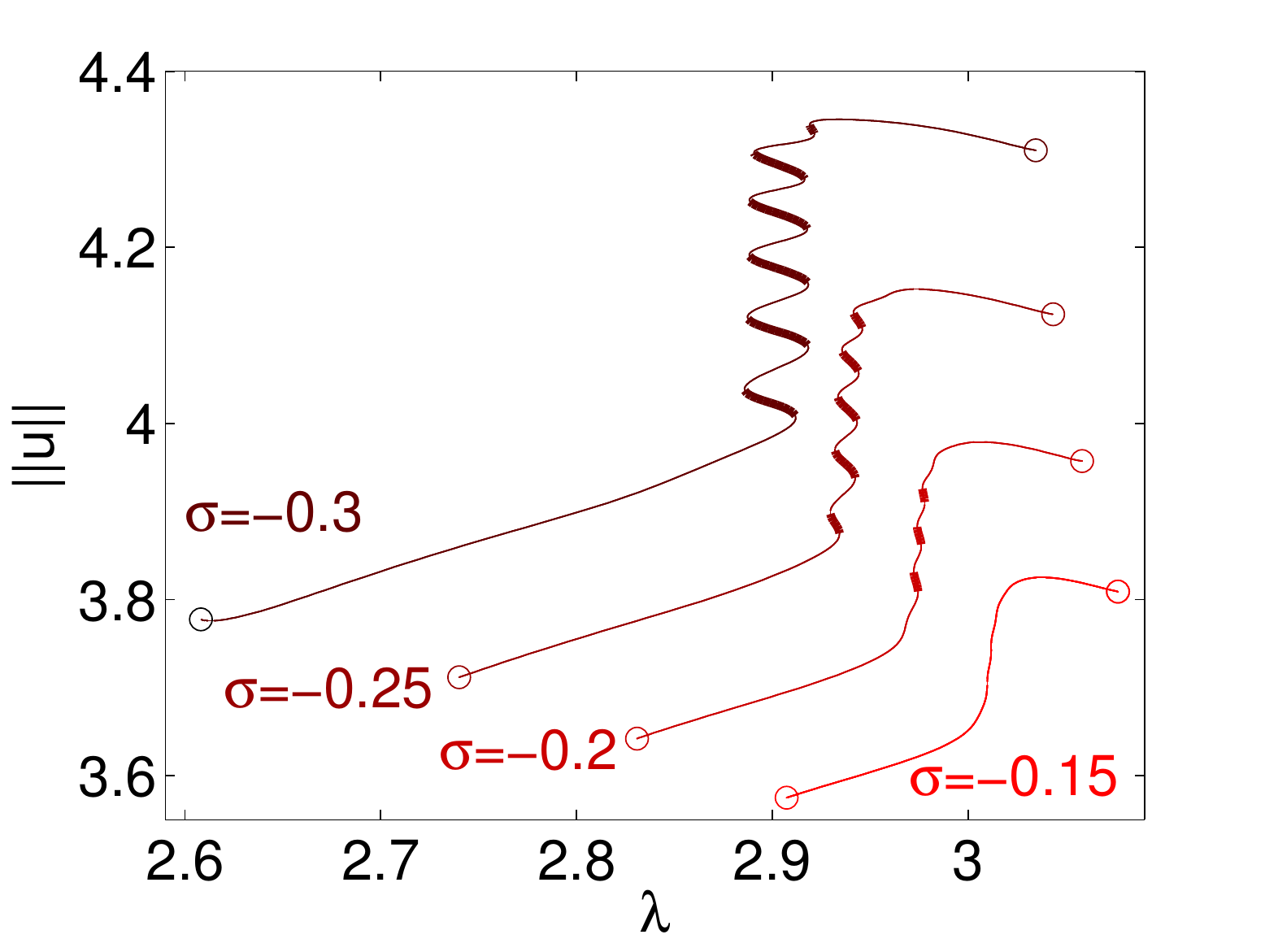}
  \end{minipage}
\ece

\vs{-5mm}
  \caption{{\small Switching on cold snaking for the modified model \reff{mod2} 
by decreasing $\sigma$. (a) Landau coefficients evaluated in $\lam_c$, 
and $c_f$, 
    as functions of $\sigma$. (b) Potential energies of cold hexagons for
    different $\sigma$ values. (c) Branches of localized hexagons over 
the homogeneous state on the $4\times 2$ domain 
for $\sigma=-0.25$, $-0.3$, and $-0.35$.
(d) $\lam_s^e(\sigma)$ and $\lam_{ch}^b(\sigma)$ such that $\lam_s^e-\lam_{ch}^b$ gives 
the width of the bistable range of 
    cold stripes and hexagons as a function of $\sigma$. (e) 
 Potential energies (solid lines) of cold hexagons and stripes (dashed).
(f) branches of localized patterns in bistable
    ranges between cold hexagons and cold stripes on the $4\times
    2$ domain.} \label{sigplot}}
\end{figure}

Here consider the modification \reff{mod1}, i.e., 
for the stationary case, 
\huga{\label{mod2}
D\Delta U+N(U,\lam)+\sigma \left(u-\frac 1 v\right)^2\bpm 1\\-1\epm=0, 
}
which has the advantage that the homogeneous solution $w=(\lam,1/\lam)$, 
the linearization 
around $w$ and thus $\lambda_c$ do not depend on $\sigma$. Consequently, 
proceeding as in \S\ref{glsec}, $\sigma$ only changes 
the Landau coefficients $c_2$,
$c_3$, and $c_4$.  For solutions of \reff{lan} we calculate that
there is a fold on the hexagon branch if the discriminant of
$P_{\pm}$ of \eqref{lanpoints} vanishes, i.e., at 
\ali{ c_1=c_f:=\frac{c_2^2}{4(c_3+2c_4)}.  
} 
Let $\lam_f$ be the corresponding $\lambda$ value, such that $\lam_f-\lam_c$ 
is a measure of the strength of the subcriticality on the Landau level. 
Near $\lam_c$ this is proportional to $-c_1$ since $c_1$ is essentially 
a linear function of $\lam-\lam_c$ for $\lam$ near $\lam_c$, see 
Fig.~\ref{huf1}(a). 
In Fig. \ref{sigplot}(a) we plot $c_f$ and the Landau
coefficients $c_2, \ c_3, \ c_4$ in $\lam_c$ as functions of 
$\sigma$. 
The value $c_f(\sig)$ and thus the subcriticality 
predicted by the Landau formalism increases
monotonously by reducing $\sigma$ from 0, and this also 
increases the potential energy difference 
for $\lam>\lam_c$, see Fig.~\ref{sigplot}(b). 
Following Remark \ref{dawrem2} and the end of \S\ref{glsec} 
this suggests that for negative $\sigma$ we may find 
snaking branches for localized cold hexagons on the homogeneous 
background, and this Ginzburg--Landau prediction 
is confirmed numerically for the system 
\reff{mod2} in (c), where snaking starts around $\sigma=-0.25$ 
and becomes stronger when further decreasing $\sig$.

Similarly, the width of the bistable range between cold hexagons and cold 
stripes and the energy differences between these patterns 
increase with decreasing $\sig$ 
(Fig.~\ref{sigplot}(d),(e)). This suggest to also find snaking branches 
between cold hexagons and cold stripes for negative $\sig$, and 
this confirmed by the numerics in (f), where 
again the snaking increases with decreasing $\sig$. 

For $\sigma_c\approx -0.369$ we have $c_3+2c_4= 0$ and
$c_f\rightarrow \mp \infty$ for $\sigma \rightarrow \sigma_c \pm 0$. 
Clearly, this cannot be used as a prediction of the fold position 
as for large $\lam-\lam_c$ higher order terms must be taken into account. 
Nevertheless, in Fig.~\ref{sigplot}(a) the plot continues 
down to $\sig=-0.9$ as the 
sign change of $c_3$ near $\sig=\sig_0\approx -0.3$ has another 
interesting consequence: for $\sig<\sig_0$ we have a subcritical 
bifurcation of stripes, and hence the possibility of 1D homoclinic 
snaking between stripes and the homogeneous background, which we illustrate 
in the next section. 

\brem\label{fcrem} (a) From $c_2(\sig_0)=0$ in Fig.~\ref{sigplot}(a), with 
$\sig_0\approx 0.3$, we expect a codimension 2 point for \reff{mod2} at  
$(\lam,\sig)=(\lam_c,\sig_c)$ with $\sig_c\approx \sig_0$. 
Then a formally consistent Ginzburg--Landau expansion 
could be performed in $(\lam_c,\sig_c)$, which together 
with beyond all order asymptotics should yield a result like  \reff{swform}, 
 i.e., $s(\eps)=C\eps^{-\al}\exp(-\beta/\eps)$. \\
(b) From the subcriticality $\eps(\sig):=\lam_f(\sig)-\lam_c$ 
and the width $s(\eps)$ of the snakes in, e.g., Fig.~\ref{sigplot}(c) 
we could also fit 
the coefficients $C,\al,\beta$ in expansions like 
\reff{swform}, and thus estimate the snaking width at $\sig=0$, 
and similar for the snaking 
between stripes and cold hexagons in Fig.~\ref{sigplot}(f). 
Here, however, we refrain from this, since we should first establish 
\reff{swform}, and since even if we assume \reff{swform}, then 
$\eps(\sig)$ and $s(\eps)$ should be calculated numerically from 
fold continuation for $\lam_f(\sig)$ and $\lam_{\mp}(\sig)$, 
where $\lam_{\mp}(\sig)$ are, e.g., the first 
left and right folds in the snakes in Fig.~\ref{sigplot}(c). 
Our current version of \pdep\ does not have this fold continuation, 
which will however be implemented in the next version. 
\eex\erem 

\subsection{1D snaking}\label{1dsec}
For \eqref{rdt} in 1D there is no bistable
range of patterns
because the hexagons do not exist and
stripes do not bifurcate subcritically, 
and thus we do not expect any localized patterns. 
For the modification \eqref{mod2} of
\eqref{rdt} we have a bistable range: the stripes bifurcate 
subcritically if $c_3$ is positive, i.e., 
for $\sig<\sig_0\approx -0.3$. We choose $\sigma=-0.6$ 
and use {\tt pde2path} to numerically investigate \reff{mod2} over the one 
dimensional domain $\Omega=(-8\pi/k_c, 8\pi/k_c)$, cf.~Remark \ref{numrem}(d). 
Then the stripe patterns branch 
$S8$ with wave number $k_c$ and hence eight periods bifurcates in 
$\lam_c$. 
From the second, third, and
fourth bifurcation point on the homogeneous branch 
stripe patterns with 8.5, 7.5, and 7 (S7 branch) 
periods bifurcate, respectively. We chose the rather small domain as 
in this section we want to give a somewhat complete picture of 
the secondary bifurcations on the S8 and S7 branches. 

\begin{figure}[h]
{\small 
  \begin{minipage}{0.5\textwidth}

(a)\\
    \ig[width=0.9\textwidth]{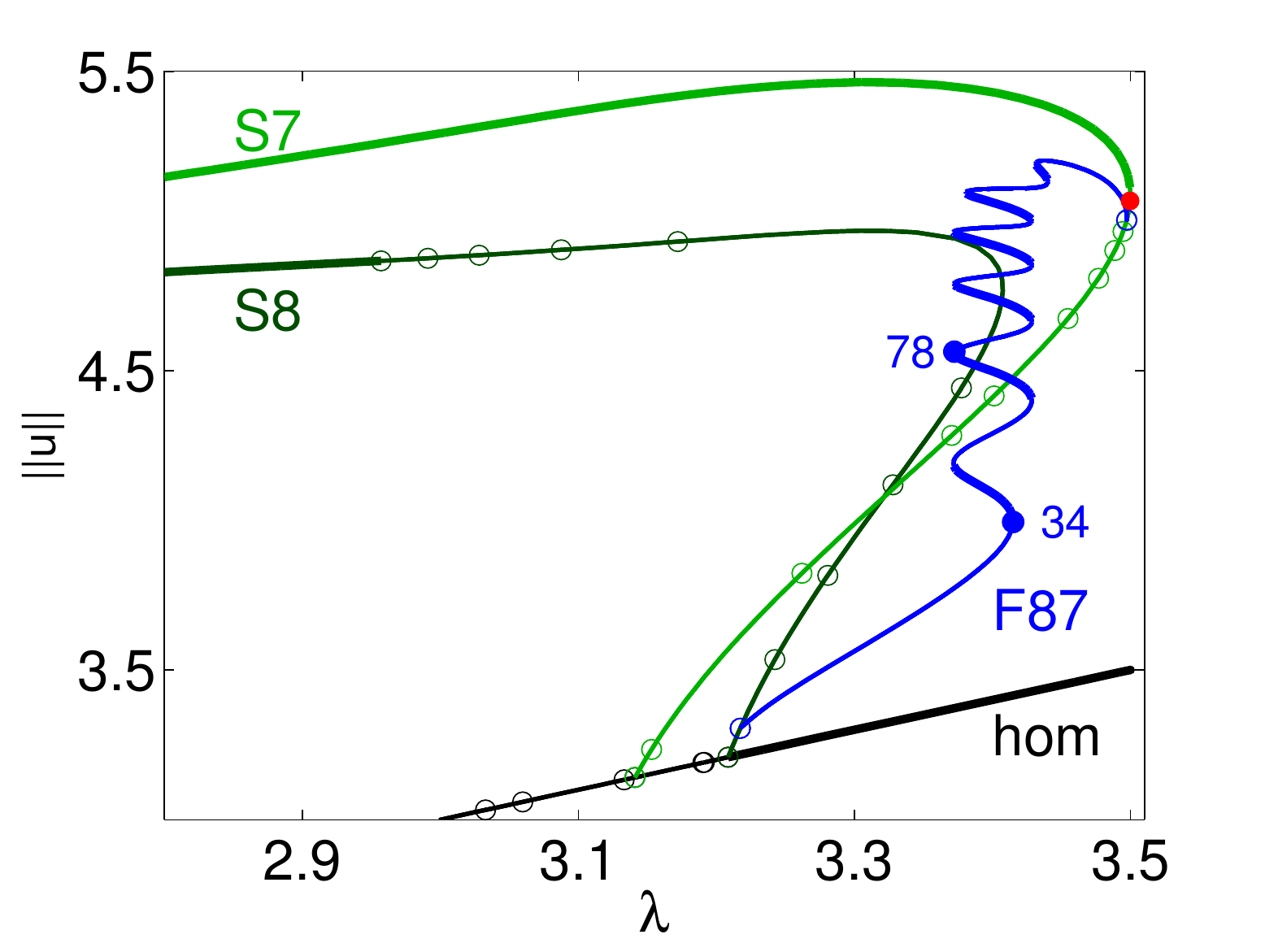}
  \end{minipage}
  \begin{minipage}{0.5\textwidth}
(b)\\
\ig[width=1\textwidth]{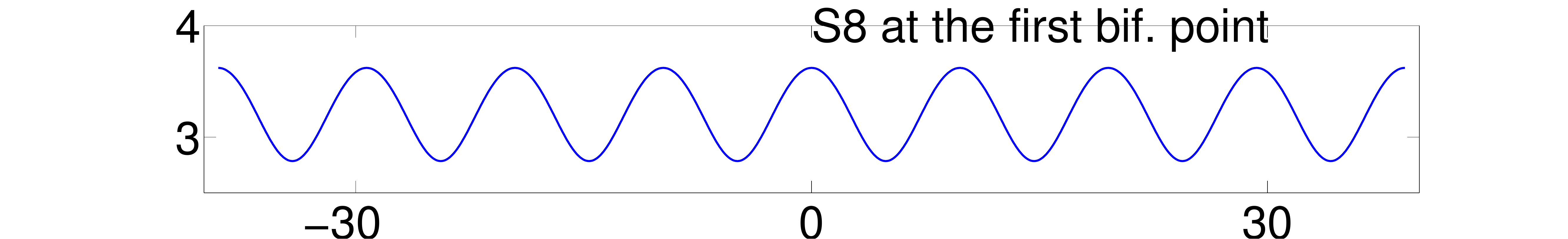}\\
 \ig[width=1\textwidth]{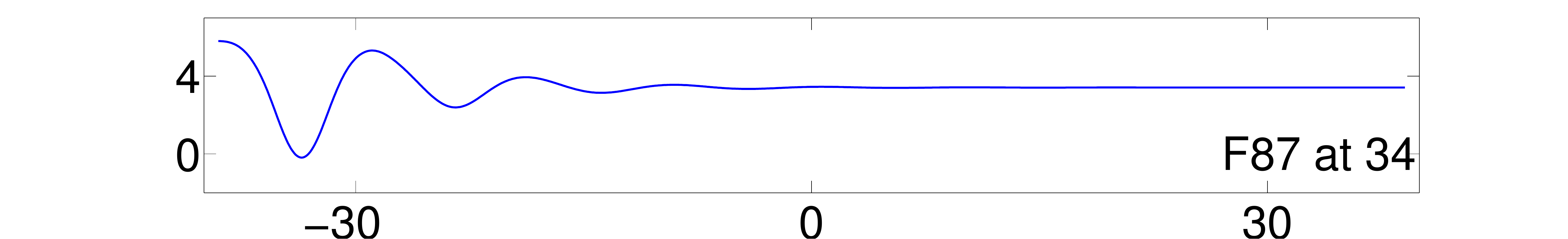}\\
    \ig[width=1\textwidth]{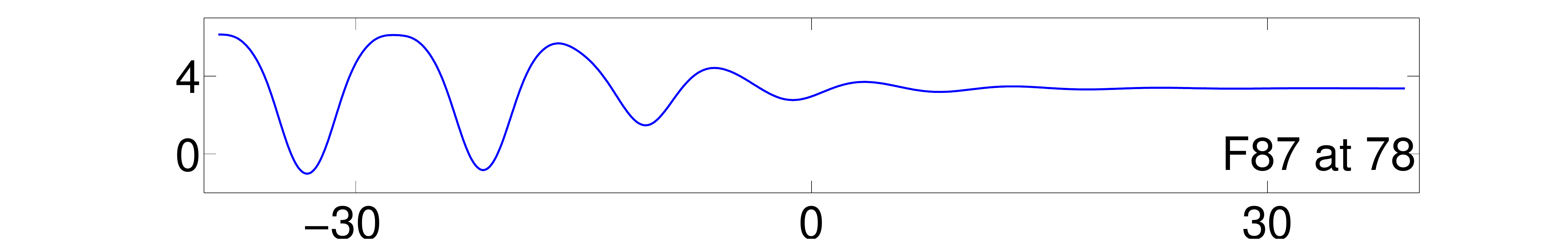}\\
     \ig[width=1\textwidth]{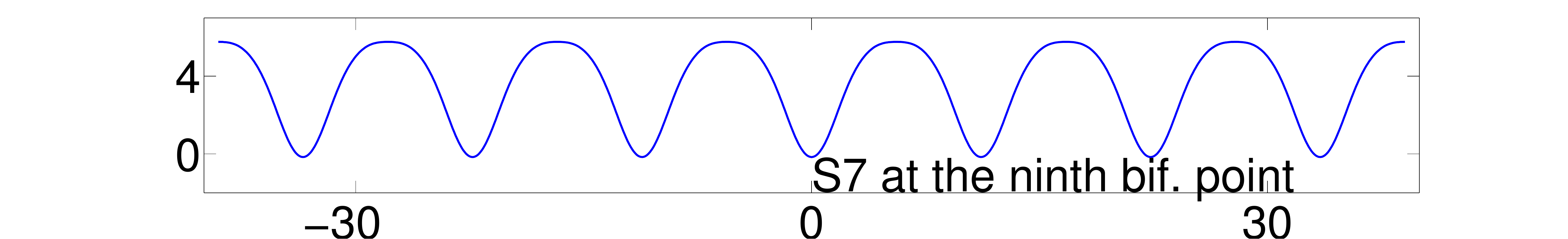}
  \end{minipage}
 \begin{minipage}{0.5\textwidth}
 (c)\\
    \ig[width=0.9\textwidth]{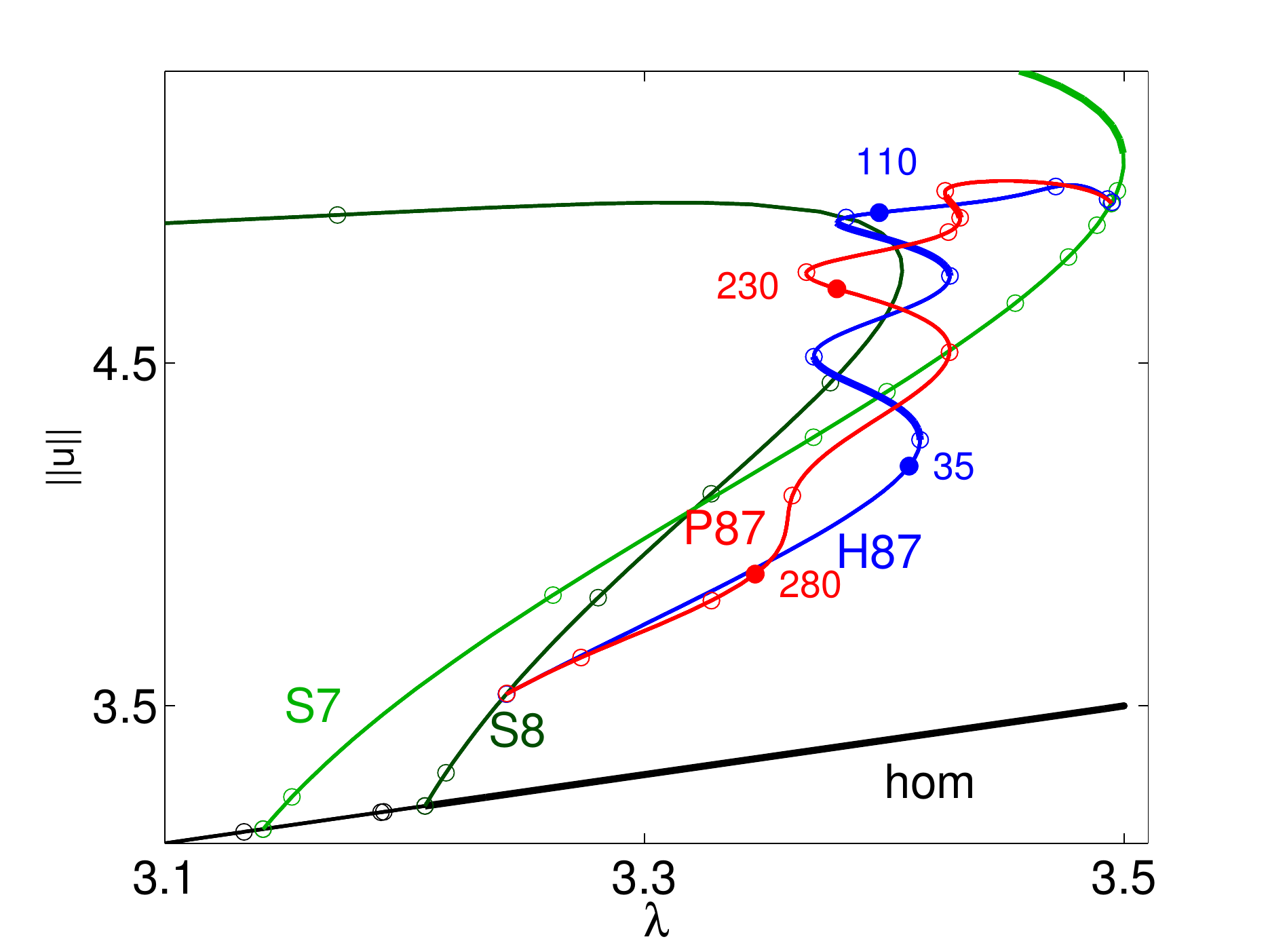}
  \end{minipage}
  \begin{minipage}{0.5\textwidth}
 (d)\\
    \ig[width=1\textwidth]{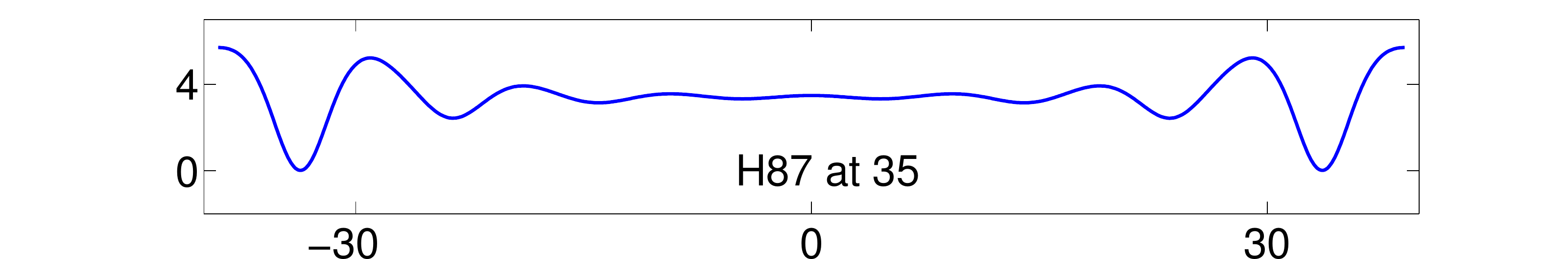}\\
    \ig[width=1\textwidth]{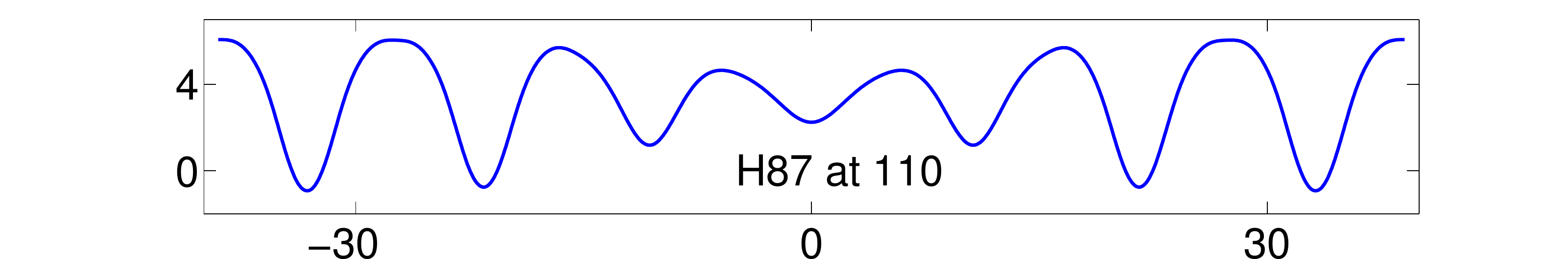}\\
    \ig[width=1\textwidth]{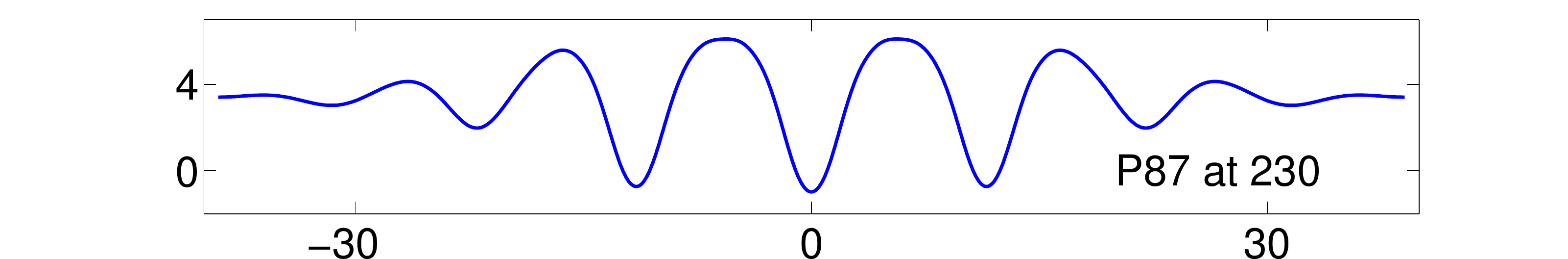}\\
    \ig[width=1\textwidth]{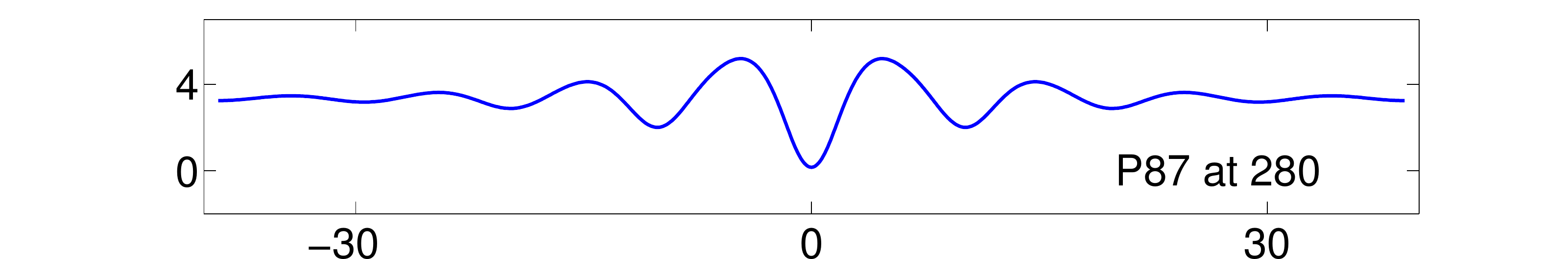}
  \end{minipage}
}
  \caption{{\small Two primary periodic branches S8 and S7, and branches 
connecting S8 and S7 for
    $\sigma{=}-0.6$ over the one dimensional domain 
$\Omega{=}(-8\pi/{k_c},
    {8\pi}/{k_c})$. (a),(b) front type; (c),(d) homoclinic type. The red dot 
in (a) marks a Hopf point, see Remark \ref{Hopf-rem}. 
See also \cite{smov} for a movie.} \label{F87}}
\end{figure}

The first interesting observation is that the S8 branch does 
not aquire stability in the fold near $\lam\approx 3.4$, but 
for smaller $\lam\approx 2.93$, see Fig.\ref{F87}(a), 
and the fold for the S7 branch is to the right of the S8 fold, 
and S7 becomes stable in its fold. 
This is somewhat similar to Fig.~9 in \cite{bbkm2008} 
where it is illustrated that also for the 1D quadratic-cubic SHe \reff{qcshe} 
the leftmost fold does not belong to the critical wave-number $k_c=1$ 
(in that case). 
Moreover, from the first bifurcation
point on S8 there bifurcates a snaking branch of stationary fronts F87 
which however does not connect back to S8, but to the 9th bifurcation 
point on S7, which is close to the fold. 
This agrees with the heuristics that branches of 
localized patterns connect patterns with bistability. 

Similarly, from the second bifurcation point on S8 we obtain a connecting 
branch to the 8th bifurcation point on S7, see 
Fig.~\ref{F87}(c),(d).  In fact, here we plot the full loop connecting these 
points, as the way from S8 to S7, which consists of solutions with 
patterns at the boundary and a {\em h}ole in the middle, 
is different from the way back, given by solutions with a (negative) 
{\em p}eak in the 
middle, and these two branches are not related by symmetry. 

\begin{figure}[h]
{\small 
\bce
  \begin{minipage}{0.28\textwidth}
(a)\\ \ig[width=50mm,height=50mm]{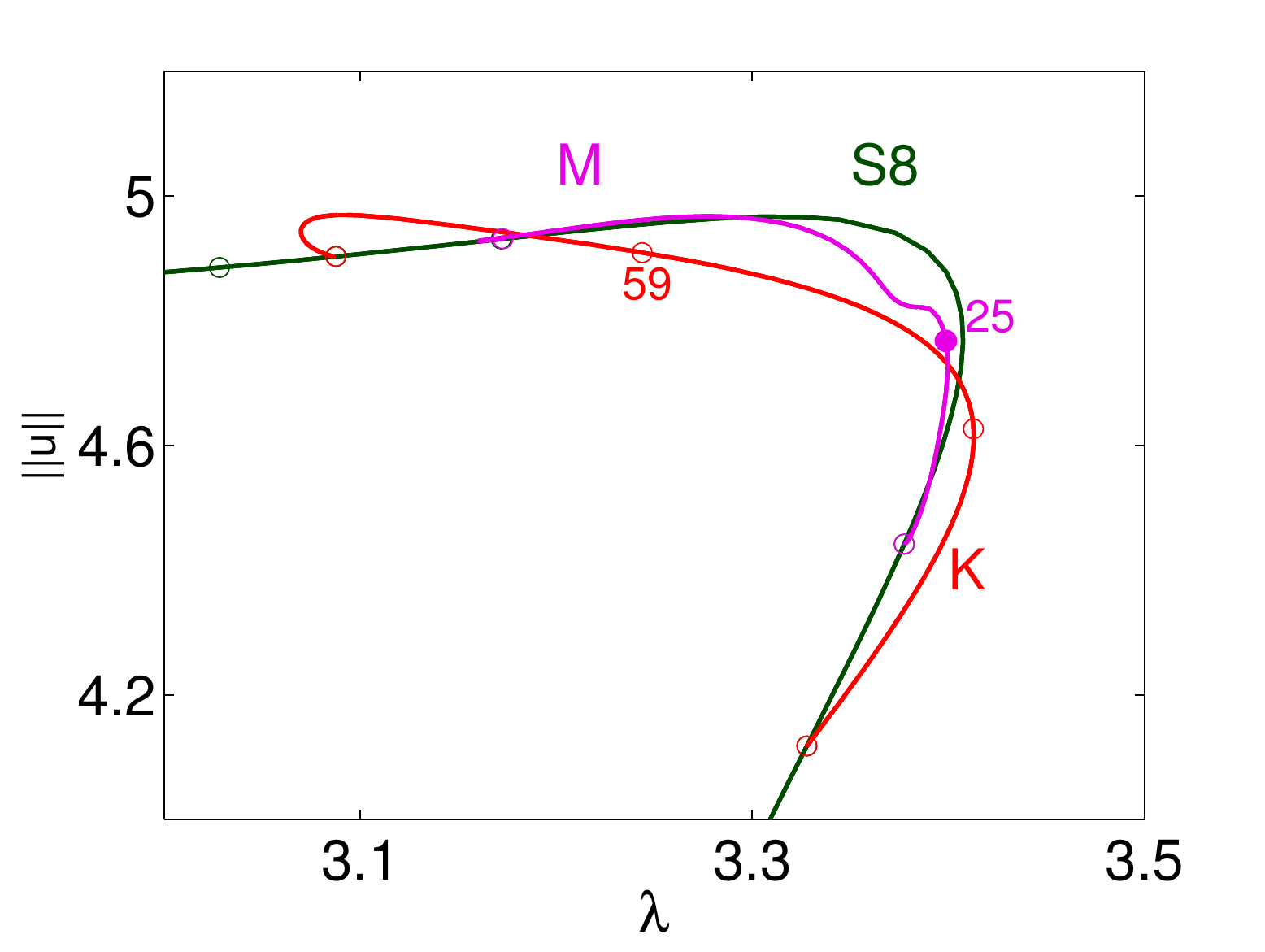}
\end{minipage}
\begin{minipage}{0.28\textwidth}
(b)\\ 
\ig[width=50mm,height=50mm]{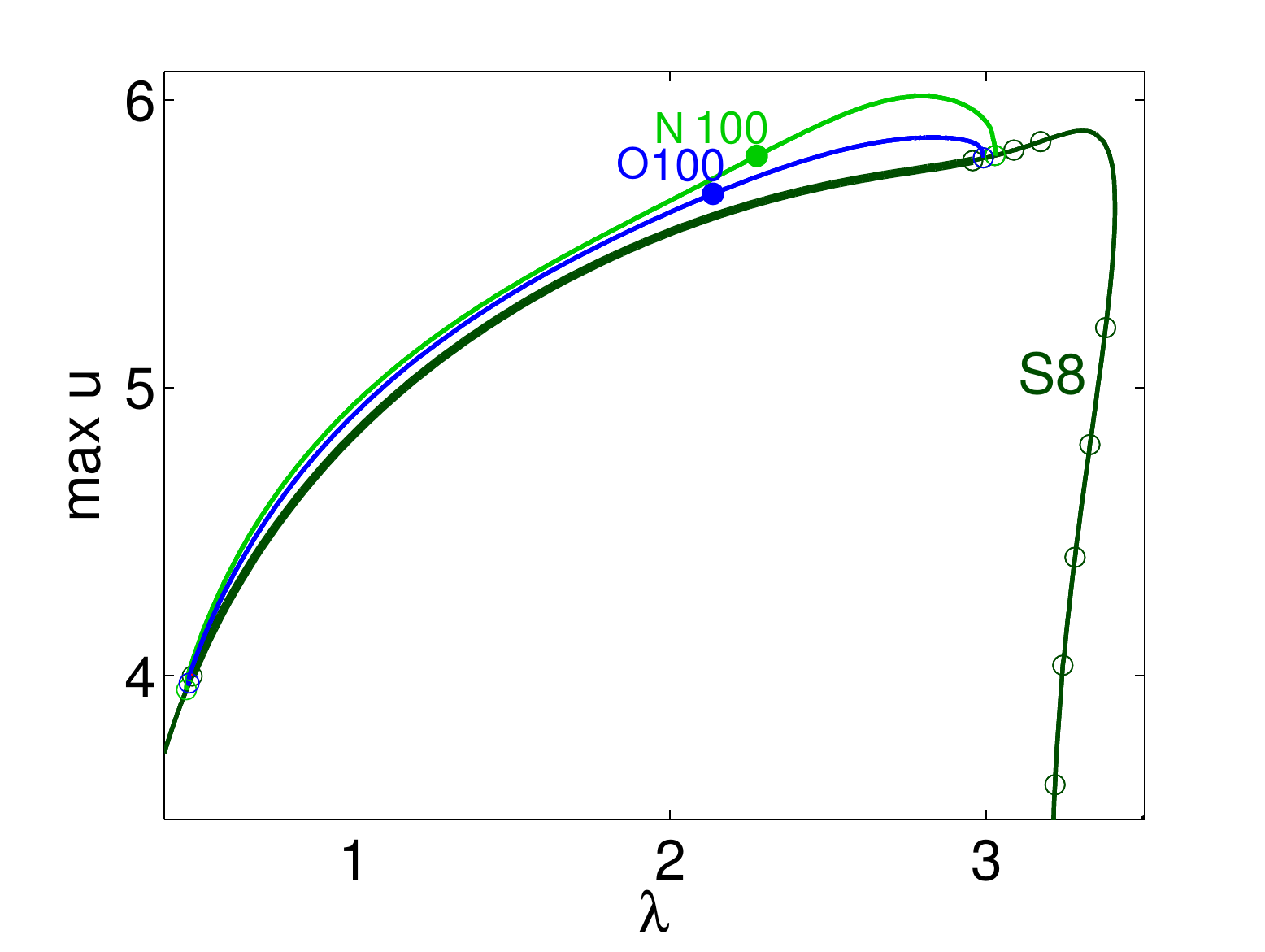}
\end{minipage}
\begin{minipage}{0.42\textwidth}
(c)\\
      \ig[width=1.1\textwidth]{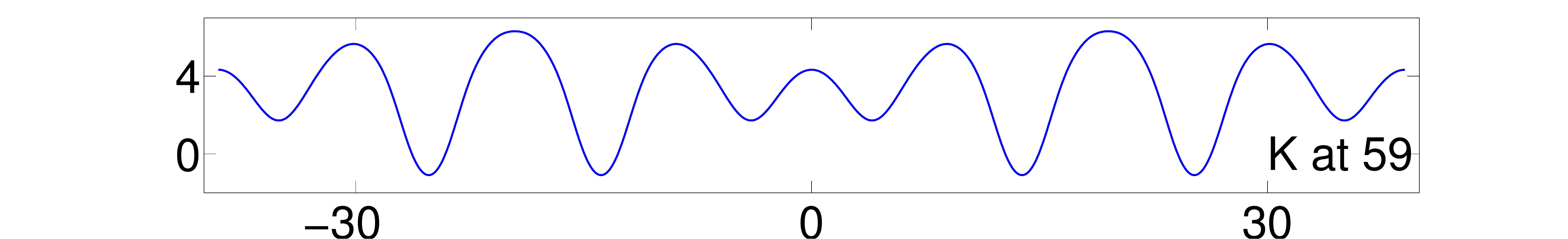}\\
      \ig[width=1.1\textwidth]{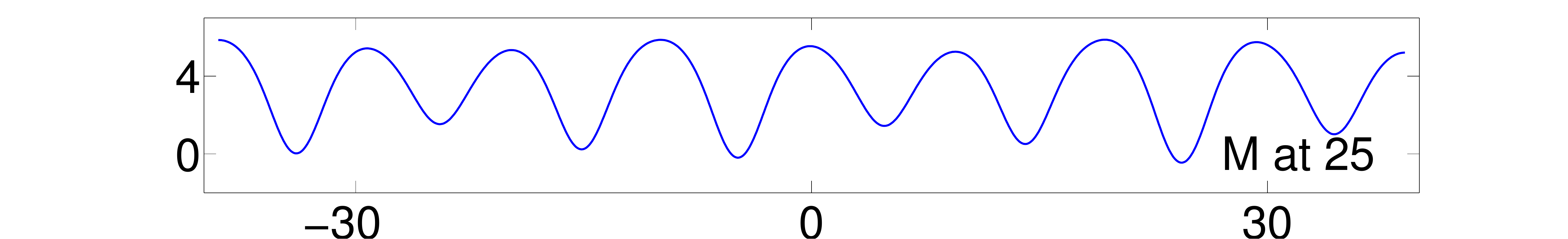}\\
  \ig[width=1.1\textwidth]{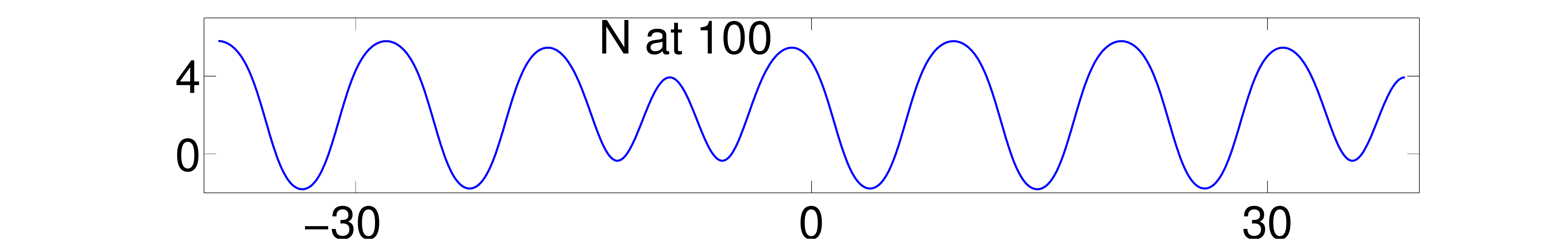}\\
    \ig[width=1.1\textwidth]{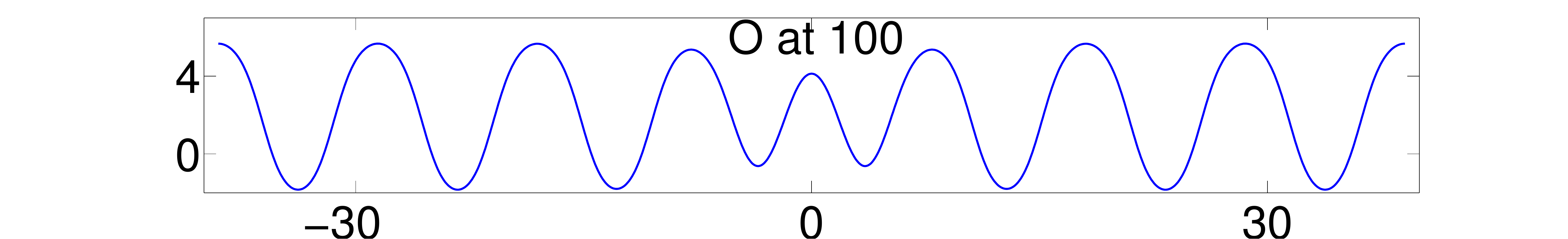}
\end{minipage}
  \caption{{\small Bifurcations from the remaining bifurcation points on S8. } \label{F89}}
\ece 
}
\end{figure}

From the third bifurcation point on S8 bifurcates a branch which 
connects to the S6.5 branch, i.e., a branch with 6.5 periods. 
The branches bifurcating from the remaining bifurcation points on S8 
reconnect to S8 as illustrated in Fig.~\ref{F89}, and similarly 
bifurcation points 1-6 on S7 connect pairwise as P1 to P6, P2 to P5, 
and P3 to P4, via branches of localized patterns without snaking.  
Finally, from 
the 7th and last unaccounted simple bifurcation point on S7 there bifurcates 
a branch which connects to the S8.5 branch. Thus, even on the small domain we have some rather complicated 
secondary bifurcations from the stripes, but of all the branches 
discussed above only S7 and the snakes connecting to S7 have some stable 
parts for $\lam>3$. 

\brem\label{Hopf-rem}{\rm Looking at the stability indices, i.e., the number 
of positive eigenvalues and the eigenvalues themselves, on the branches 
we also find a Hopf point  
near the fold of S7, see Fig.~\ref{F87}(a). This, however,  
will be studied elsewhere.}\eex\erem
\section{Outlook: additional numerical results}
\label{addsec}
\subsection{Fully localized spot patches over radial stripes}
\label{radsec}
Besides the straight (planar) stripes considered so far, we may 
expect so called radial stripes which only depend on the radius $r$ and 
are asymptotically periodic in $r$, see, e.g.,  \cite{scheel-ams2003}. Moreover, 
for the 2D quadratic--cubic and cubic--quintic SHe \reff{qcshe} and 
\reff{cqshe} also a number of localized radial stripe patterns 
are known,  
e.g., radial pulses (or spots) and rings \cite[Figure 1 and \S5]{hexsnake}, 
see also \cite{slloyd09}, and again these patterns can also be expected for RD systems 
with a Turing instability, cf.~\cite[Remark 1]{slloyd09}. 

Thus we may search for such solutions of \reff{rd} or \reff{mod1}, 
and may study their bifurcations to patterns 
with (dihedral) $D_n$ symmetry, including 
some fully localized patterns. Here we restrict to \reff{rd}, 
i.e., let $\sig=0$ again,  
over circular sectors with opening angle $\pi/3$, thus 
effectively restricting to patterns with $D_6$ symmetry. Again we choose 
some rather small domain 
$$
\Om_r=\{(x,y)=\rho(\cos\phi,\sin\phi): 0<\rho<r, 
0<\phi<\pi/3\},$$ 
with Neumann boundary conditions, and $r=12\pi/k_c$ 
(this  precise 
value of $r$ is rather arbitrary, even if we expect radial stripes 
with asymptotic period $2\pi/k_c$), 
which yields some interesting results at acceptable 
numerical costs. (Typically we need meshes of about $60.000$ triangles, 
locally refined near $(x,y)=(0,0)$, to avoid branch jumping.)

\begin{figure}[!ht]
{\small 
\bce
\begin{minipage}{\textwidth}
(a) bifurcation directions (first components) at selected points on 
the homogeneous branch. 

\ig[width=40mm]{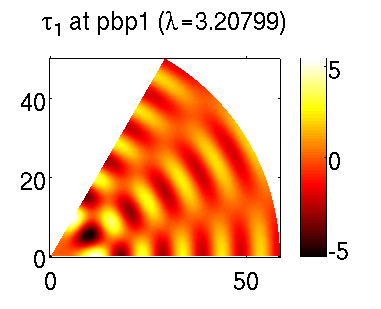}
\ig[width=40mm]{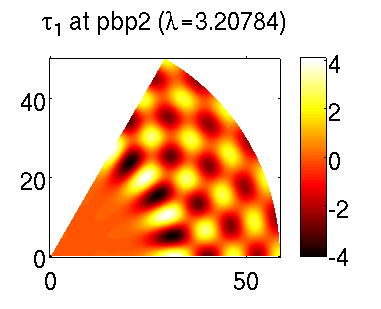}
\ig[width=40mm]{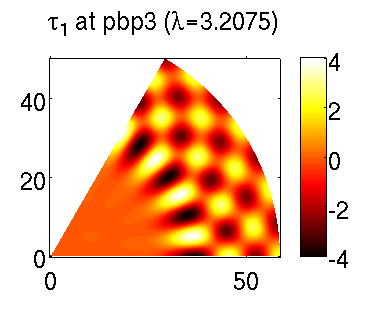}
\ig[width=40mm]{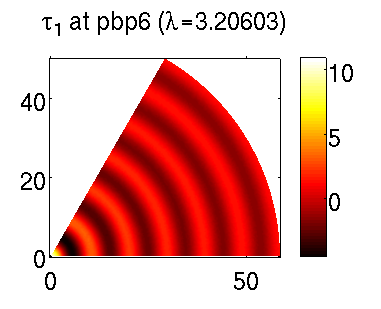}
\end{minipage}

\begin{minipage}{40mm}(b) bifurcation diagram of radial stripes\\
\ig[width=40mm]{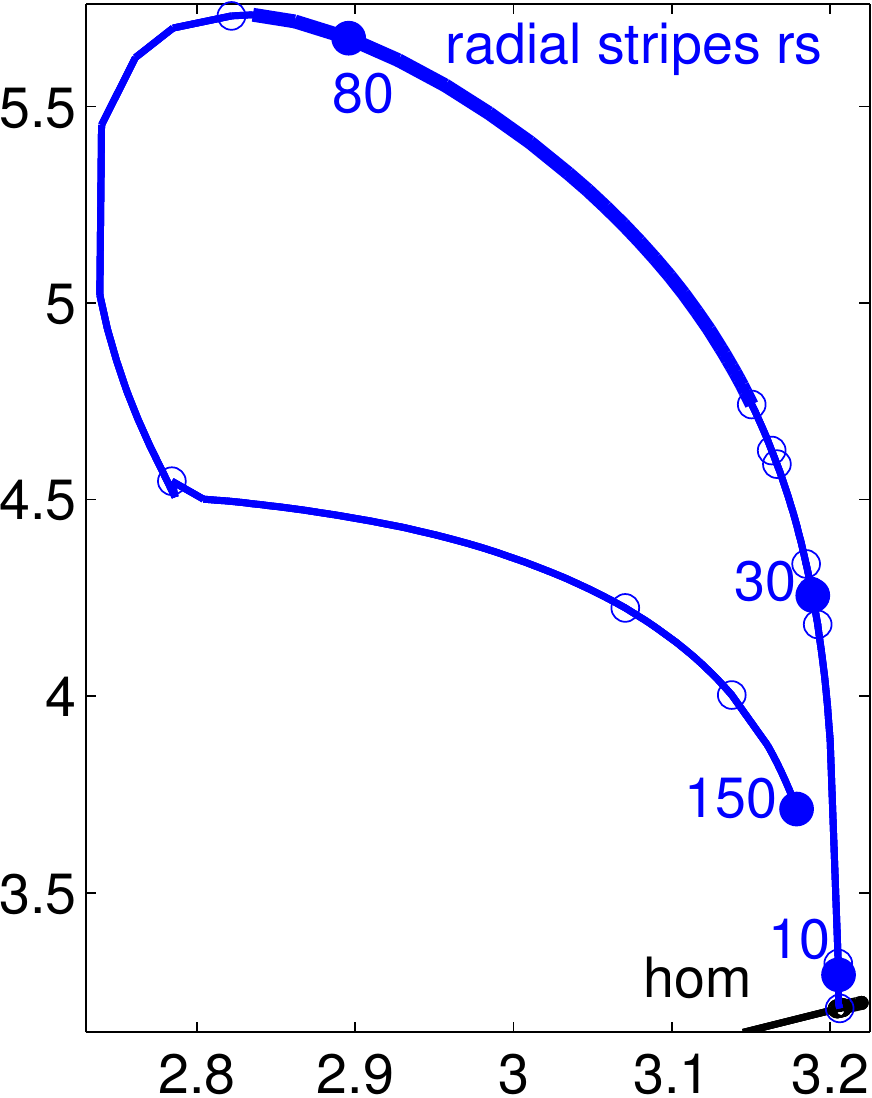}
\end{minipage}\quad
\begin{minipage}{40mm}(c) radial stripe at c80, \newline and bifurcation 
direction $\tau_1$ at the nearby bif.~point. \\
\ig[width=40mm]{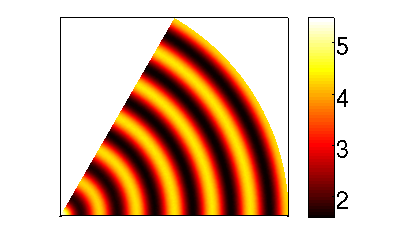}\\
\ig[width=40mm]{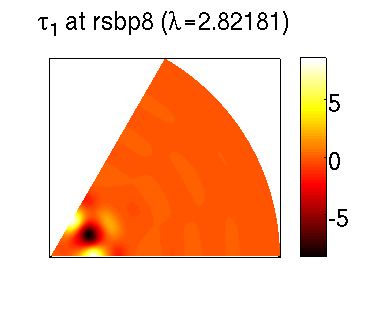}
\end{minipage}\quad
\begin{minipage}{70mm}(d) selected radial profiles\\
\begin{tabular}{p{36mm}p{34mm}} 
rs10 and rs30& rs80 and rs150\\
\ig[width=35mm]{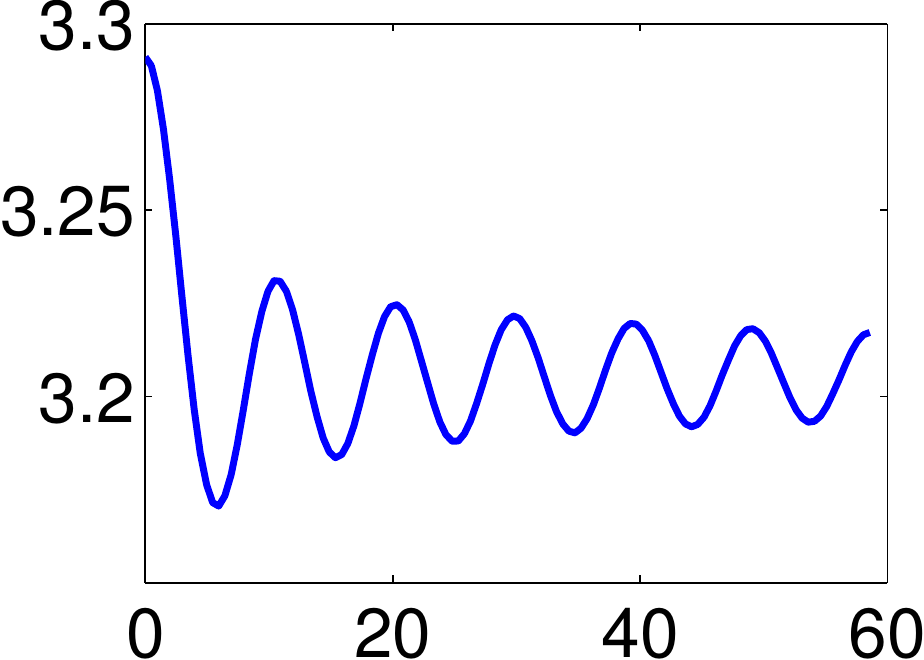}&\ig[width=32mm]{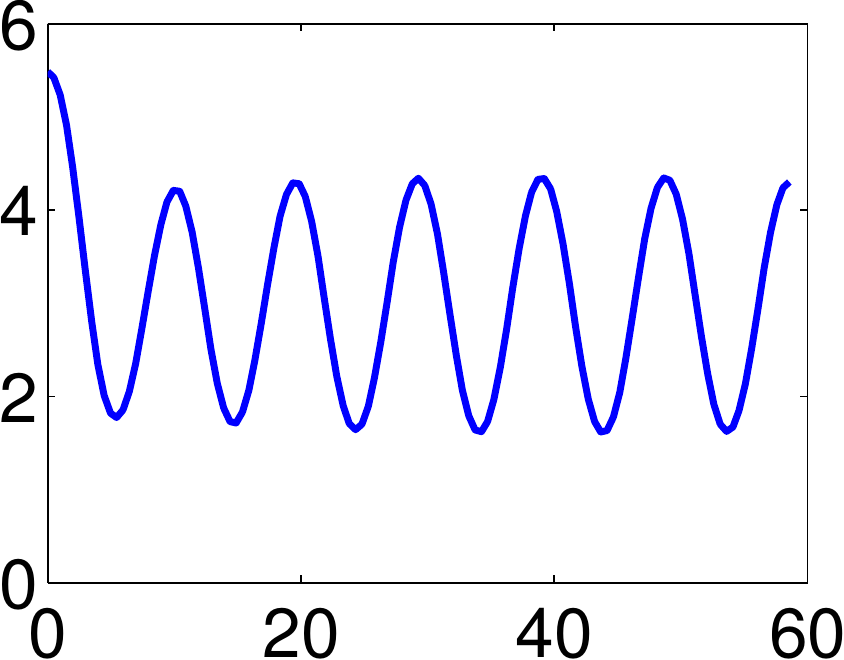}\\
\quad\ig[width=32mm]{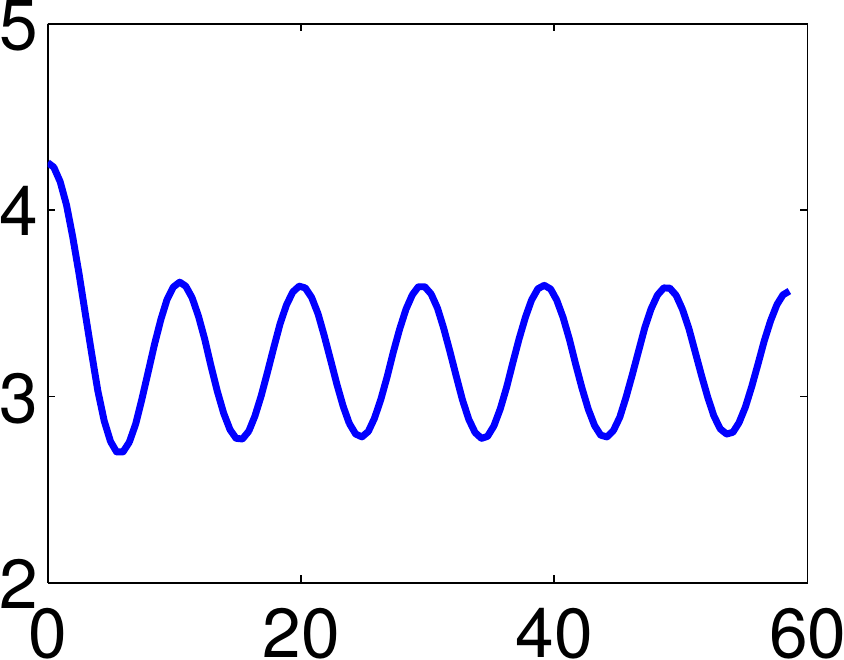}&\ig[width=32mm]{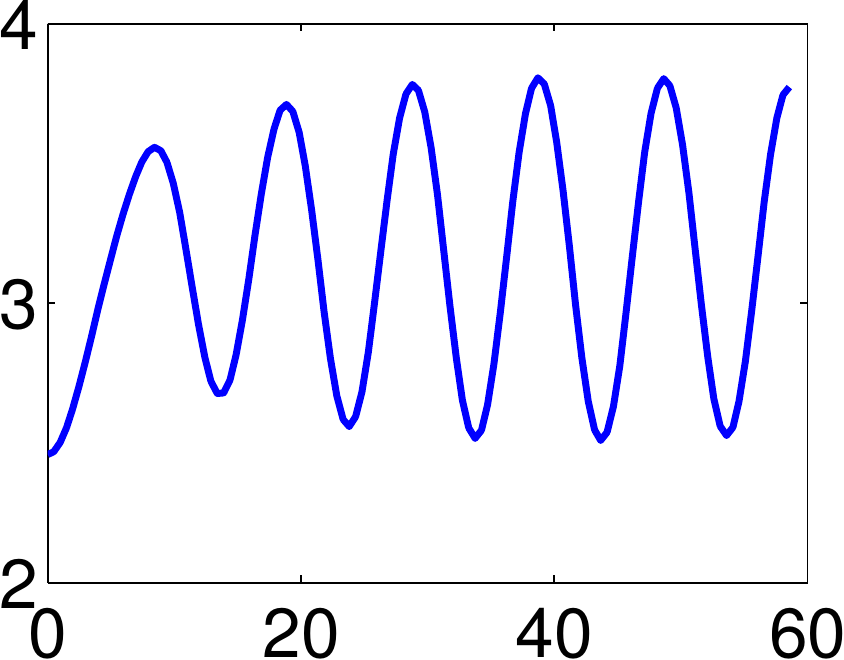}
\end{tabular}
\end{minipage}
\ece 
}

\vs{-5mm}
\caption{{\small Bifurcation of radial stripes from the homogeneous branch 
$w(\lam)=(\lam,1/\lam)$ over 
a sector with radius $r=12\pi/k_c$ and angle $\pi/3$. (a) Some bifurcation 
directions from $w$. (b)-(d) Bifurcation diagram and some example plots, 
including the bifurcation direction to the spot patch branch, see 
Fig.~\ref{cf2}.}  \label{cf1}}
\end{figure}
Already the bifurcations from the trivial branch (Fig.~\ref{cf1}(a)) 
show a number of interesting directions, but we focus on 
the branch {\tt rs} ({\tt r}adial {\tt s}tripes, see Remark \ref{radrem}(iii)
on that terminology) coming from branch point 6. The bifurcation diagram in (b) 
shows that similar to the straight stripes from 
\S\ref{numrel} the radial stripes bifurcate supercritically but unstable, 
gain stability at some $\lam_{rs}^b\approx 3.15$ and loose stability 
again at $\lam=\lam_{rs}^e\approx 2.82$. The example plots in (c),(d) 
show that {\tt rs} starts as a ``hot stripe'' with a peak in $0$, and after 
losing stability turns around to return to the homogeneous branch 
as a ``cold stripe'' with a slightly longer period. See also Remark 
\ref{radrem}(iv) on the ``cold stripe'' half of {\tt rs} branch.  

\begin{figure}[!ht]
{\small 
\bce
\begin{minipage}{33mm}(a) bifurcation of spot patch branch from radial stripes.\\
\ig[width=38mm]{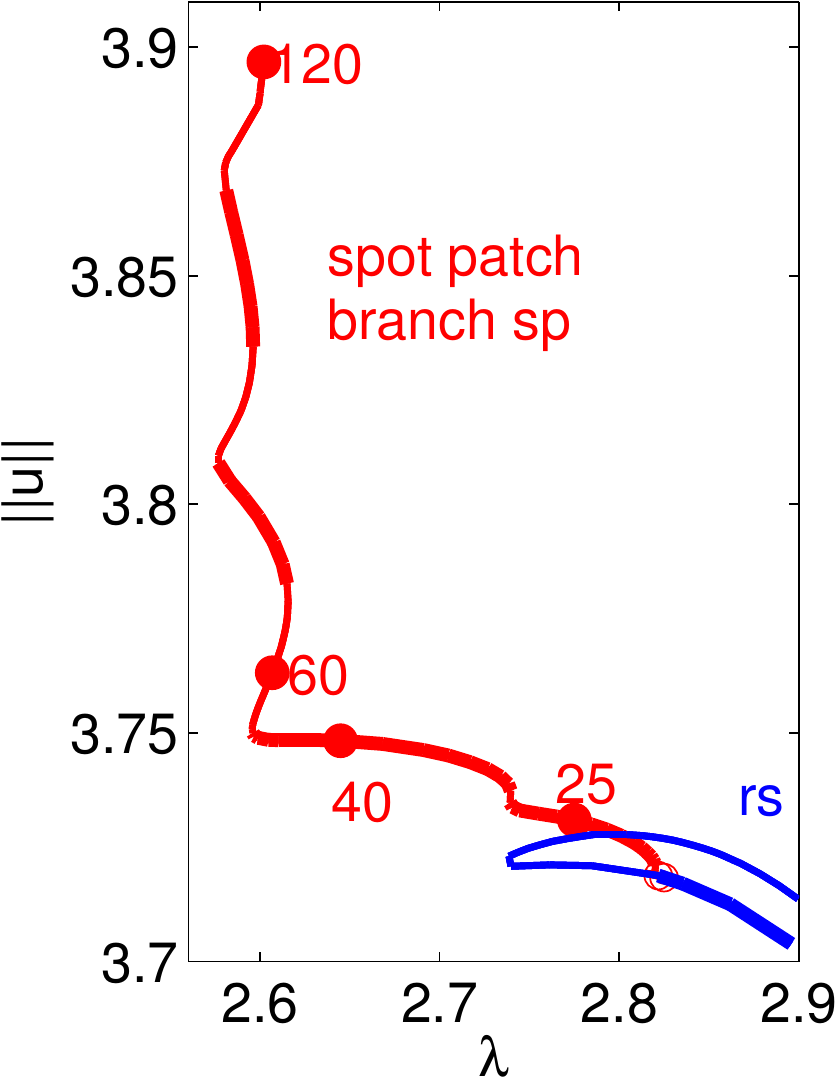}
\end{minipage}\quad
\begin{minipage}{70mm}\text{}\quad(b) example plots: sp25, sp40 (top) and\\
 \text{}\quad sp60,sp120 (bottom)\\
\begin{tabular}{p{34mm}p{34mm}}
\ig[width=35mm]{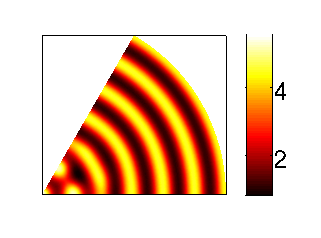}&\ig[width=35mm]{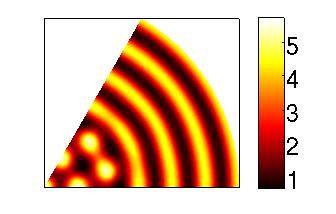}\\
\ig[width=35mm]{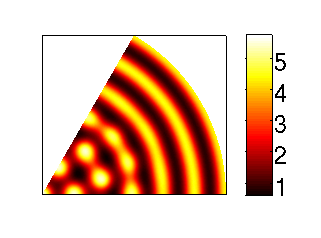}&\ig[width=35mm]{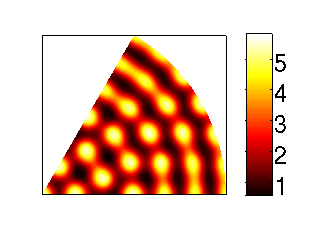}
\end{tabular}
\end{minipage}\quad
\begin{minipage}{54mm}\text{}\quad(c) sp60 with rotations by $j\frac \pi 3$.
\ig[width=54mm]{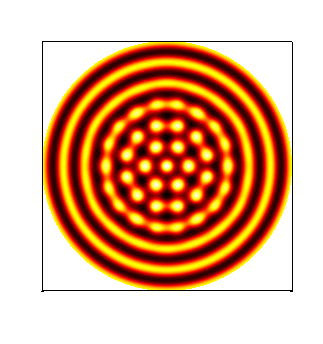}
\end{minipage}
\ece 
}

\vs{-5mm}
\caption{{\small Bifurcation of (hexagonal, in center) spot patch branch 
{\tt sp} from radial stripes.  See \cite{smov} for a movie.
}  \label{cf2}}
\end{figure}

In Fig.~\ref{cf2} we show the bifurcation of spot patches from the 
radial stripes at $\lam_{rs}^e$. The tangent at bifurcation Fig.~\ref{cf1}(c) 
shows that we may expect one sixt of a regular hexagon patch to appear 
near $0$, and this  
indeed happens on the {\tt sp} branch. During continuation the branch snakes  
as additional spots are added, with alternating stable and 
unstable parts of the branch. Except for the central hexagon, the 
added spots yield somewhat distorted 
hexagon patches, and as the spots approach the radial boundary of the domain 
the numerical continuation becomes difficult: the mirror symmetry 
along angle $\pi/6$ is lost due to numerical inaccuracy 
somewhere beyond point 120, the stepsize 
decreases quickly, and eventually continuation fails. 
 Figure \ref{cf2}(c) shows the five-fold rotation of {\tt sp60}, while 
 Fig.~\ref{ipic2}(d) in the Introduction is obtained from {\tt sp40} in the same way.

As this section is only intended as an outlook, and as we plan to study 
at least some of these patterns elsewhere, including some analytical 
approaches in the spirit of \cite{slloyd09}, here we close with the 
following remarks. 
\brem\label{radrem} 
\bcen[(i)]
\item Branch switching at some other bifurcation 
points on the homogeneous branch gives a number of interesting different 
patterns, see Fig.~\ref{cf1}(a) for some example bifurcation directions. 
However, in contrast to the radial stripes, all of these appear to depend 
strongly on the somewhat artificial choice of the bounded domain, in particular 
the radial Neumann boundary conditions. 
\item Some of the earlier bifurcation points on the {\tt rs} branch do give 
localized patterns, where, e.g., one of the ``middle'' 
stripes is replaced by a spot ring, 
but these are unstable. 
\item From the bifurcation direction in Fig.~\ref{cf1}(a), and from 
the radial profiles near bifurcation, e.g., {\tt rs10} in (d), 
it appears that the {\tt rs} branch 
should rather be classified as a radial spot branch, i.e., 
solutions show radially decaying oscillations. However, since the 
branch continuously turns into non-decaying oscillations we find 
the name radial stripes more appropriate. Presumably, one needs 
significantly larger domains to cleanly distinguish radial spots and 
stripes. The same applies to radial rings, 
which following \cite[Theorem 1, Remark 1]{slloyd09} we may expect 
for  \reff{mod1} with $\sig<\sig_c$ where we have a subcritical 
bifurcation of (straight) stripes. 
\item From the ``cold'' branch bifurcating from the homogeneous branch at 
the same point as the ``hot'' {\tt rs} branch but in the opposite 
direction we obtain a bifurcation of a ``cold'' spot branch.  
\item It appears natural to conjecture that also localized radial stripes 
(i.e., spots or rings) 
over a background of hexagons should exist. To find these we 
considered triangular domains as in Fig.~\ref{huf4}, which 
support hexagon patterns, and looked for radial stripes near $(x,y)=0$ by 
two methods: (a) via bifurcation as the hot hexagons gain stability 
roughly at $\lam_{hh}^b$ (the start of the bean branch over 
rectangular domains); 
(b) by using suitable initial guesses and time integration or a 
direct Newton loop. Both methods failed so far.
\ecen
\eex\erem

\subsection{An interface perpendicular to stripes}\label{pisec} 
\begin{figure}[!ht]
\bce
{\small
\begin{tabular}{ll}
(a) bifurcation diagram &(b) b1-50, b1-90, hh{-}10, f1-20, l1-80, and l1-180\\
\ig[width=48mm]{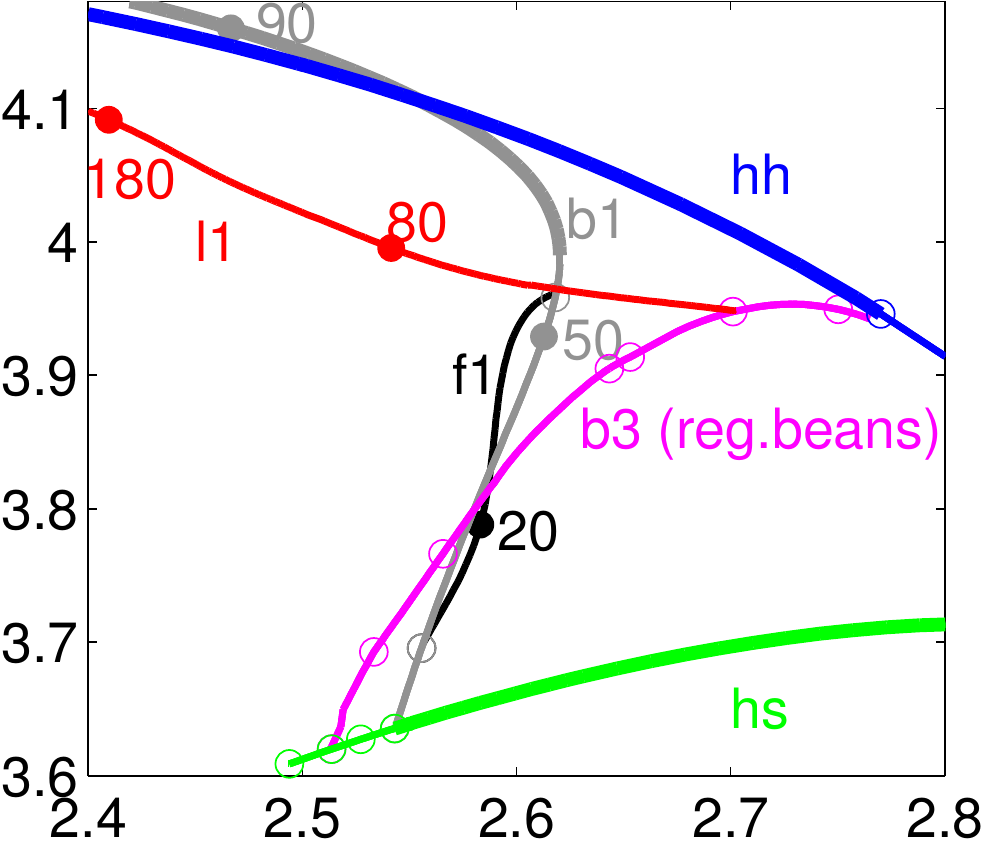}&
\text{}\hs{-4mm}\ig[width=27mm]{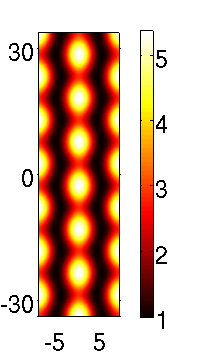}\hs{-4mm}
\ig[width=20mm]{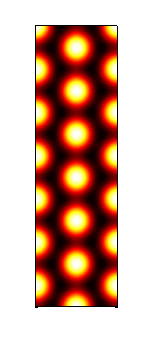}\hs{-4mm}
\ig[width=20mm]{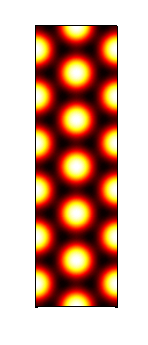}\hs{-4mm}
\ig[width=20mm]{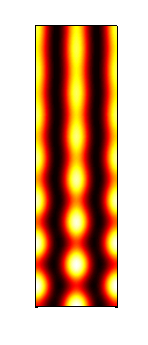}\hs{-4mm}
\ig[width=19.5mm]{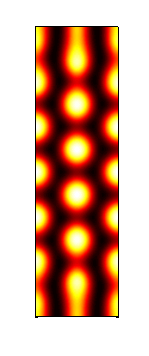}\hs{-4mm}
\ig[width=19.5mm]{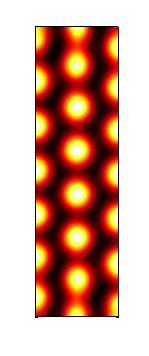}
\end{tabular}
}\ece

\vs{-5mm}
\caption{{\small Bifurcations from the stripes in the hot regime ($\sig=0$) 
over the $1\times 6$ domain.} \label{1x6f}}
\end{figure}

For all planar fronts so far 
we chose different lengths in $x$ but kept the $y$ dimension 
fixed and rather small. Together with the choice of vertical stripes 
this gave interfaces parallel to the stripes, and essentially this 
also holds true for the interfaces to radial stripes in \S\ref{radsec}.  
We now aim at more general interfaces, and thus first consider 
{\em long} vertical stripes, i.e., large $y$--lengths. Figure \ref{1x6f} 
shows some branches of solutions and example plots over the $1\times 6$ 
domain $\Om=(-\frac{2\pi}{k_c},\frac{2\pi}{k_c})\times (-\frac{12\pi}
{\sqrt{3}k_c},\frac{12\pi}{\sqrt{3}k_c})$. Decreasing $\lam$, 
the stripes no longer lose stability at $\lam_s^e\approx 2.51$ 
to the regular beans 
with $k_{2,3}=\frac {k_c} 2 (-1,\sqrt{3})$, but earlier to 
streched beans with wave vectors 
$k_{2,3}=\frac {k_c} (-\frac 1 2, \pm \frac {13}{12}\frac{\sqrt{3}}2)$, 
denoted by {\tt b1}, with 6.5 spots in $y$--direction. 
 This branch  {\em turns} into a streched hexagon branch, becoming 
stable in a fold.
The regular beans {\tt b3} bifurcate at the third 
bifurcation point after loss of stability of the stripes and contain a number 
of bifurcation points. Two of the bifurcation points on {\tt b3} 
are connected by a front in $x$--direction. However, on the 
example branch {\tt l1} the solutions first 
change wave-number near the top and the bottom and then become very similar to 
hexagons with 7 spots in $y$--direction, and similar effects occur on 
branches from the remaining bifurcation points on {\tt b3}. 

The {\tt b1} branch contains two bifurcation points which 
are connected by a branch containing a rather flat vertical front between 
stripes and 
(streched) hexagons, i.e., the interface is perpendicular to the stripes. 
We did not find similar branches bifurcating from other bean 
branches, in particular not from the regular beans {\tt b3}. We do not 
know, why this is so. 

The numerics are considerably more difficult on 
the $1\times 6$ domain compared to, e.g., the $4\times 2$ from 
Fig.~\ref{snakea}. For instance, in Fig.~\ref{1x6f} 
we had to refine the mesh to about 100000 points on the bifurcating 
branches to avoid some behaviour like increasingly 
tilted interfaces during continuation, which do not seem to be due 
to properties of the PDE, but due to finite size and discretization effects, 
i.e., poor meshes.

\subsection{Fully localized (sideband) hexagon 
patches over straight stripes}
\label{flsec}
In our final simulation we return to the modified system \reff{mod2}, with 
$\sig=-0.3$, mainly because the numerics for fully localized hexagons 
over stripes 
appear to be easier in that regime, than, e.g., in the hot regime from 
Fig.~\ref{1x6f}, where $y$--fronts are quite flat, while $x$--fronts are 
rather steep, cf.~\S\ref {lpsec}. Since over large domains we then again find 
that the selected wave vectors for patterns are 
not necessarily $k_1,k_2,k_3$ in a hexagonal lattice with $|k_j|=k_c$, 
here we choose 
the large square domain $\Omega=(-8 \pi
/ k_c,8 \pi / k_c )\times(-8 \pi / k_c,8 \pi / k_c ) $.  
Note the missing $1/\sqrt{3}$ in the $y$--length, such 
that this does not allow regular hexagons with wave-vectors 
$k_{2,3}=k_c(-\frac 1 2, \pm \frac {\sqrt{3}}2)$. Instead, similar to 
\S\ref{sidesec} and Fig.~\ref{1x6f} we expect sideband hexagons 
with $k_{2,3}$ near 
$k_c(-\frac 1 2, \pm \frac {\sqrt{3}}2)$, and the closest wave vectors 
allowed by the domain 
are $k_{2,3}=k_c(-\frac 1 2, \pm \frac 3 4)$. 

\begin{figure}[!ht]
\bce
{\small 
\begin{minipage}{0.4\textwidth}
(a) bifurcation diagram\\
\includegraphics[width=1\textwidth, height=60mm]{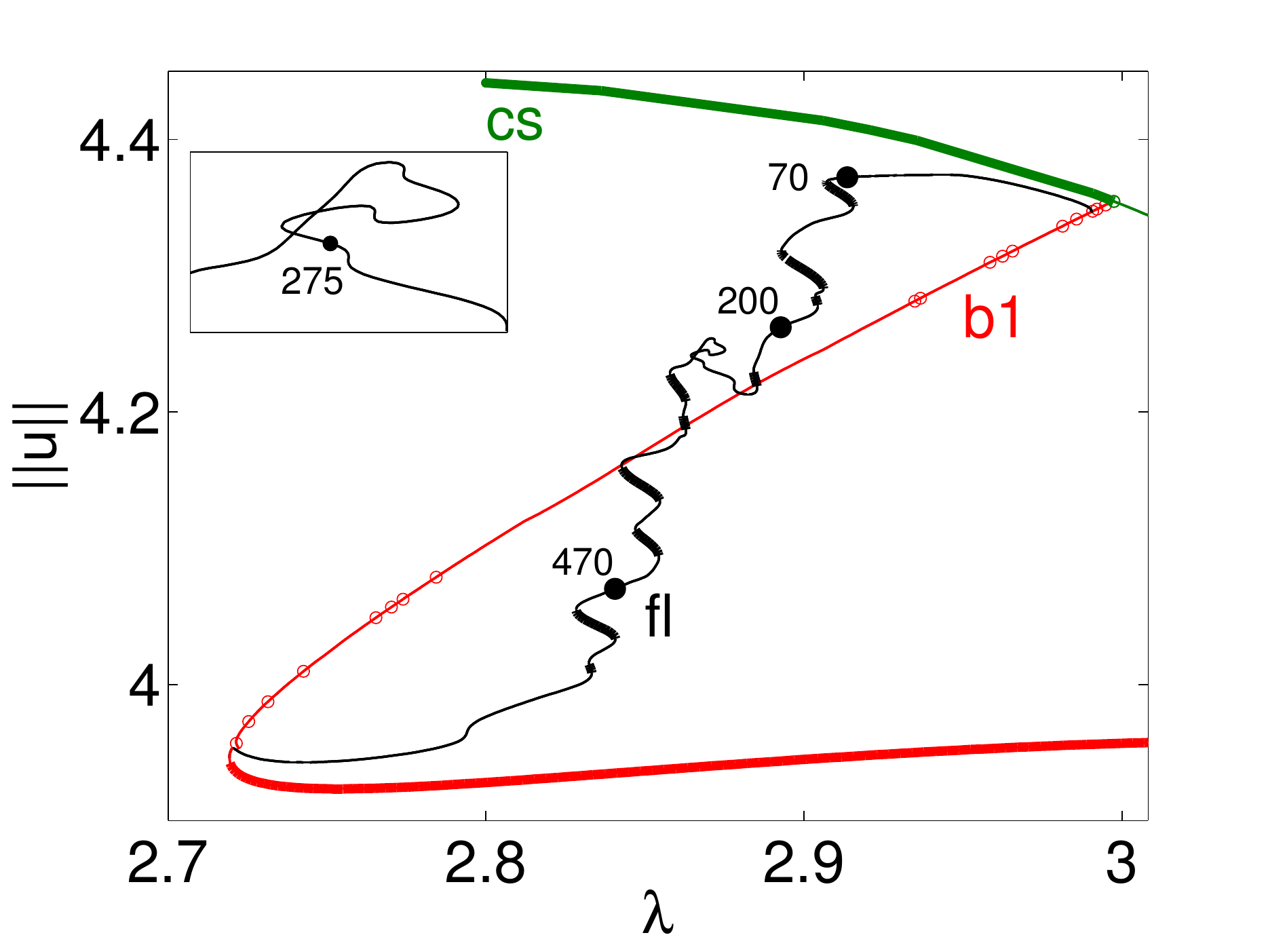}
\end{minipage}
\begin{minipage}{0.5\textwidth}
\begin{tabular}{l}
(b) fl70, fl200 (top), and fl275, fl470 (bottom)\\
\includegraphics[width=37mm]{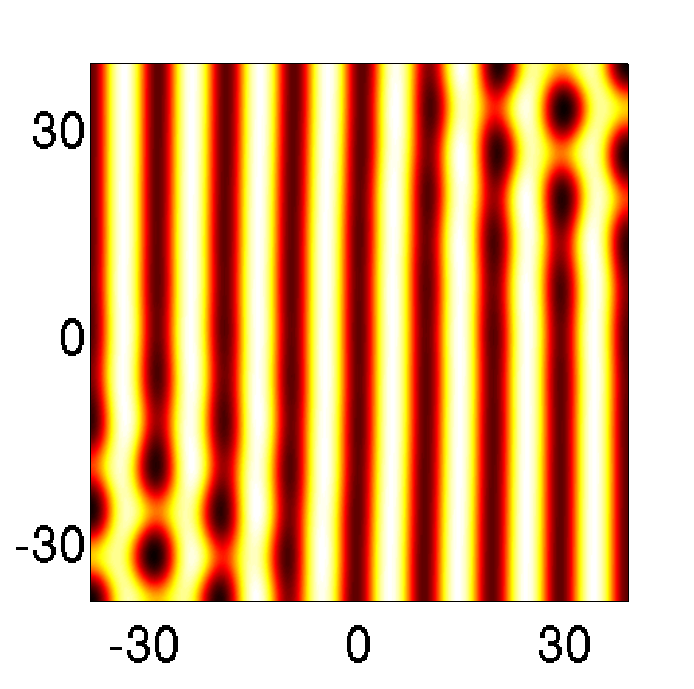}
\includegraphics[width=40mm]{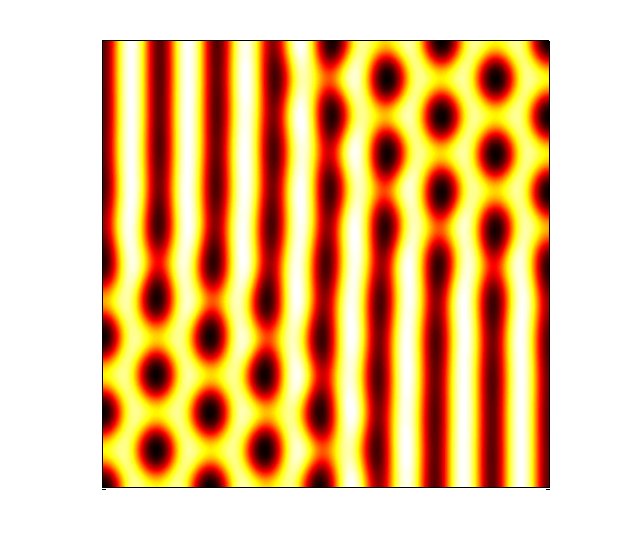}\\
\text{}\hs{-3mm}
\includegraphics[width=40mm]{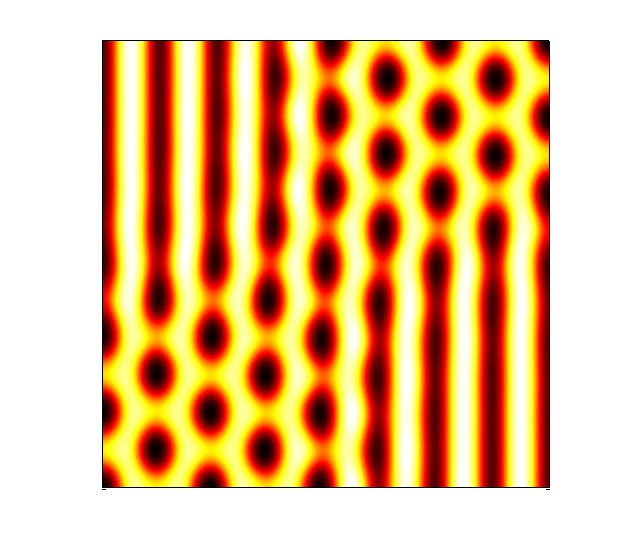}\hs{-2mm}
\includegraphics[width=40mm]{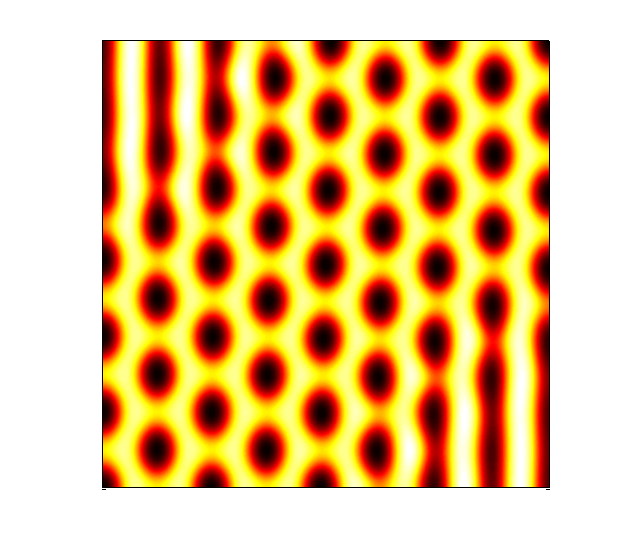}
\end{tabular}
\end{minipage}
}
\ece 

\vs{-6mm}
\caption{{\small Fully localized patches of stretched hexagons 
on a stripe background for $\sigma=-0.3$ over the square domain 
$\Omega=(-8 \pi / k_c,8 \pi / k_c )\times(-8 \pi / k_c,8 \pi / k_c) $. 
Colormap in all plots roughly between 0 and 5. 
See \cite{smov} for a movie.}\label{flplot}}
\end{figure}
We start at $\lam=2.8$ in the stable part of the vertical cold stripes
{\tt cs} for $\sigma=-0.3$ and follow the branch in positive $\lam$
direction. When the stripes lose their stability, there bifurcates a
bean-branch {\tt b1}, which {\em turns} into a branch of (streched) cold
hexagons and connects to the homogeneous branch at $\lam$  
very close to $\lam_c$.  
As expected, the patterns of {\tt b1} 
are stretched in vertical direction, i.e., $k_{2,3}=k_c(-1/2,\pm
3/4)$. From the first bifurcation point on {\tt b1} bifurcates a
(snaking) branch of horizontal planar fronts between the 
{\tt b1} hexagons and the stripes {\tt cs}, i.e., similar to Fig.~\ref{huf3}, 
or rather similar to Fig.~\ref{sigplot}(f), 
while from the second bifurcation
point bifurcates a (snaking) branch of planar fronts 
perpendicular to the stripe direction, similar to {\tt f1} 
in Fig.~\ref{1x6f}. 

Finally, from the third bifurcation point on {\tt b1} 
a branch {\tt fl} of fully localized hexagon patches bifurcates. 
At the beginning the solutions have a patch of (stretched) 
hexagons in two corners, see {\tt flp70} of Fig. \ref{flplot}(c). 
The branch snakes
and the patches grow towards the middle ({\tt flp200}). At {\tt flp275}
the patches get together and the branch makes a loop, such that afterwards 
we rather have two patches of stripes on a hexagon background, see, e.g., 
{\tt flp400}. 

Clearly, as for \S\ref{radsec}, many open questions remain for the patterns 
calculated in Fig.~\ref{1x6f} and \ref{flplot}. 

\section{Discussion, and Open Problems}\label{dsec}
As for localized patterns over homogeneous backgrounds, the main 
ingredients for localized patterns over a different background pattern 
is a bistable range, usually generated by a subcritical bifurcation 
(of mixed modes, in the case of two patterns). For planar interfaces 
the existence and approximate location of the pinning region can be predicted 
if there exists a heteroclinic connection between the corresponding 
fixed points in the associated Ginzburg--Landau system at 
some Maxwell point $\lam=\lam_m$. The pinning 
effect in the full system then yields the existence of heteroclinics in a 
parameter interval around $\lam_m$, and the same effect yields branches 
of localized patterns in a patterned background even if the
 Ginzburg--Landau system only has ``approximate'' homoclinic solutions. 

These results have already been suggested in 
\cite{pomeau}, and have been further worked out in 1D problems 
and for patterns over homogeneous 
background in various papers cited above, 
but \S\ref{numrel} of 
the present paper appears to be the first numerical illustration 
for different patterns in 2D. Moreover, in \S\ref{addsec} we gave a numerical 
outlook on some classes of fully localized patterns over patterns. 
All the numerics are somewhat delicate 
due to the very many solution branches that exist in the Turing unstable 
range, in particular over large domains, which require rather fine 
discretizations to avoid uncontrolled branch switching. 
Here we restricted to domains of still small to intermediate size. 

The results are certainly not special for the model problems \reff{rd}, 
\reff{mod1}, and we have for instance used the same 
method to predict and numerically find various snaking branches in other 
2D reaction diffusion systems \cite{dwfeu13}. 

Our work has been almost entirely numerical, and at many places we just give 
first steps towards understanding the new patterns found. Main Open Problems 
and ``things to do''  
include (in order of appearance): 
\bcen[(a)]
\item Clarify why the snaking in the hot regime is strongly slanted, 
see Fig.~\ref{nnf3}. 
\item\label{td1} Derive a formula like \reff{swform} giving 
exponential smallness 
of the snaking width in the subcriticality parameter $\eps$. For \reff{mod1} 
this should probably be done near the codimension 2 point $(\lam_c,\sig_c)$, 
cf.~Remark \ref{fcrem}(a). Alternatively, discuss this near a codimension 
2 point 
for the full Selkov--Schnakenberg model using (at least one of) 
the parameters $a,b$ discussed after \reff{rdt}. 
\item\label{td2} Following (\ref{td1}), or independently, calculate more 
numerical data 
for the relation of subcriticality and snaking width using fold continuation, 
cf.~Remark \ref{fcrem}(b). 
\item\label{td3} Detect and study Hopf bifurcations in \reff{rd} or \reff{mod1}, 
see Remark \ref{Hopf-rem} for just one example. 
\item Study, both numerically and analytically, fully localized patterns 
over patterns in more detail; in \S\ref{addsec} we gave just a very first 
outlook of these. 
\ecen 
Additionally, it might be interesting to study the stick--slip 
motion (Fig.~\ref{ssfig}) near a snake in more detail, 
including the derivation of equations of motion for interfaces. 
For (\ref{td2}) and (\ref{td3}), here we refrain from setting up ad hoc 
modifications of \pdep, but instead plan to include the 
needed upgrades (multi--parameter continuation and Hopf bifurcation) 
in a general way into the package.  

\renewcommand{\refname}{References}
\renewcommand{\arraystretch}{1.05}\renewcommand{\baselinestretch}{1}
\small
\newcommand{\etalchar}[1]{$^{#1}$}

\end{document}